\documentclass{report}

\usepackage{amsfonts, amsmath, amssymb, amscd, amsthm}
\usepackage{graphicx}

\newtheorem{thm}{Theorem}[chapter]
\newtheorem{cor}[thm]{Corollary}
\newtheorem{lem}[thm]{Lemma}
\newtheorem{prop}[thm]{Proposition}
\newtheorem{prob}[thm]{Problem}
\newtheorem{conj}[thm]{Conjecture}
\theoremstyle{definition}
\newtheorem{defn}[thm]{Definition}
\theoremstyle{remark}
\newtheorem{rem}[thm]{Remark}
\newtheorem{ex}[thm]{Example}
\numberwithin{section}{chapter} \numberwithin{equation}{chapter}
\numberwithin{figure}{chapter}

\newcommand{\G }{\Gamma (G, X\cup\mathcal H)}
\newcommand{\dxh }{dist_{X\cup\mathcal H}}
\newcommand{\dx }{dist_X}
\newcommand{\dxl }{dist_{X\cup \mathcal L}}

\newfont{\eufm}{eufm10}

\begin{document}






\title{Relatively hyperbolic groups:
Intrinsic geometry, algebraic properties, and algorithmic
problems}
\author{\Large \rm D. V. Osin}
\date{}
\maketitle



\tableofcontents


\chapter{Introduction}



\section{Preliminary remarks}


Originally, the notion of a relatively hyperbolic group was
proposed by Gromov \cite{Gr1} in order to generalize various
examples of algebraic and geometric nature such as fundamental
groups of finite--volume non--compact Riemannian manifolds of
pinched negative curvature, geometrically finite Kleinian groups,
word hyperbolic groups, small cancellation quotients of free
products, etc. Gromov's idea has been elaborated by Bowditch in
\cite{B1}. (An alternative approach was suggested by Farb
\cite{F}.) In the present paper we obtain a characterization of
relative hyperbolicity in terms of isoperimetric inequalities and
adopt techniques based on van Kampen diagrams to the study of
algebraic and algorithmic properties of relatively hyperbolic
groups. This allows to establish a background for the subsequent
paper \cite{Osin}, where we use relative hyperbolicity to prove
embedding theorems for countable groups.

Since the words 'relatively hyperbolic group' seem to mean
different things for different people, we briefly explain here our
terminology. There are two different approaches to the definition
of the relative hyperbolicity of a group $G$ with respect to a
collection of subgroups $\{ H_1, \ldots , H_m \} $. The first one
was suggested by Bowditch \cite{B1}. It is similar to the original
Gromov's concept and characterizes relative hyperbolicity in terms
of the dynamics of properly discontinuous isometric group actions
on hyperbolic spaces. (For exact definitions we refer to the
appendix).

In the paper \cite{F}, Farb formulated another definition in terms
of the coset graphs. In the simplest case of a group $G$ generated
by a finite set $S$ and one subgroup $H\le G$ it can be stated as
follows. $G$ is {\it hyperbolic relative to} $H$ if the graph
$\widehat \Gamma (G,S)$ obtained from the Cayley graph $\Gamma (G,
S)$ of $G$ by contracting each of the cosets $gH$, $g\in G$, to a
point is hyperbolic. In fact, the hyperbolicity of $\widehat
\Gamma (G,S)$ is independent on the choice of the finite
generating set $S$ in $G$.

The two definitions were compared in \cite{Sz}, where
Szczepa\'nski showed that if a group $G$ is hyperbolic with
respect to a collection of subgroups $\{ H_1, \ldots , H_m \} $ in
the sense of Bowditch, then $G$ is hyperbolic with respect to $\{
H_1, \ldots , H_m \} $ in the sense of Farb, but not conversely.
However, in \cite{F} Farb do not simply consider relatively
hyperbolic groups. He introduces an additional (and central in his
theory) condition, the so--called Bounded Coset Penetration
property (or BCP, for brevity). It turns out that the notion of
the relative hyperbolicity with BCP in the sense of Farb is
equivalent to the notion of the relative hyperbolicity in the
sense of Bowditch \cite{B1,Bum1,DahPhD}.

In order to define relative hyperbolicity of a group $G$ with
respect to a collection of subgroups $\{ H_1, \ldots , H_m \} $,
the approaches of Farb and Bowditch require $G$ to be finitely
generated as well as $ H_1, \ldots , H_m $ (although the last
assumption is rather technical). In the present paper we introduce
a more general definition which is based on relative isoperimetric
inequalities. This allows us to eliminate the assumption of the
existence of finite generation sets of $G$ and $ H_1, \ldots , H_m
$ as well as the assumption of the finiteness of the collection of
subgroups. Such a generalization is important in certain
applications and, in particular, allows to include the small
cancellation theory over free products (see \cite[Ch. V, Sec.
9]{LS}, \cite[Ch. 11]{Ols-book}) within the framework of the
theory of relatively hyperbolic groups. In case the group $G$ is
finitely generated, our notion of relative hyperbolicity is
equivalent to that of Bowditch and to that of Farb with the BCP
property.

Before stating the main theorems, let us survey certain motivating
examples for which, in particular, the results of our paper can be
applied.

(I) Let $M^n$ be a noncompact, complete, finite-volume Riemannian
manifold with (pinched) negative sectional curvature
$$-b^2\le K(M^n)\le -a^2<0.$$  Then $\pi _1(M^n)$
is hyperbolic in the sense of Bowditch with respect to the
collection of cusp subgroups \cite{F} (for the definition of cusp
subgroups we refer to \cite{E}). The examples of such a type
include, for instance, non--uniform lattices in real $\mathbb
R$--rank one simple Lie groups.

(II) Let $G$ be a $C^\prime (1/6)$--small cancellation quotient of
the free product of groups $X_1, \ldots , X_k$ (see \cite[Ch. V,
Sec. 9]{LS}). Then $G$ is hyperbolic relative to the natural
images of subgroups $X_i$ in $G$ in the sense of Bowditch. This
result follows directly from our characterization of relative
hyperbolicity in terms of relative isoperimetric inequality and
the Greendlinger Lemma \cite{LS} (see also
\cite{Pan1}--\cite{PanPhD}, where 'hyperbolic factorizations' of
free products are studied).

(III) Suppose that $H_1, H_2, \ldots , H_k$ are quasi--convex
subgroups of a word hyperbolic group $G$. Then $G$ is hyperbolic
with respect to the collection $H_1, H_2, \ldots , H_k$ in the
sense of Farb \cite{Ger}. If, in addition, $|H_i^g\cap H_j|<
\infty $ whenever $g\notin H_i$ or $i\ne j$, then $G$ is
hyperbolic relative to $H_1, H_2, \ldots , H_k$ in the sense of
Bowditch \cite[Theorem 7.11]{B1}.

(IV) Yaman \cite{Yaman} characterized relatively hyperbolic groups
as geometrically finite convergence groups (see Appendix for
definitions). More precisely, suppose that $M$ is a non--empty
perfect compact metric space and suppose that a group $G$ acts as
a geometrically finite convergence group on $M$. Suppose also that
the stabilizer of each bounded parabolic point is finitely
generated. Then $G$ is hyperbolic relative to the set of the
maximal parabolic subgroups in the sense of Bowditch and $M$ is
equivariantly homeomorphic to the boundary of $G$. In the case of
ordinary hyperbolic groups (i.e. in case the set of bounded
parabolic points of $M$ is empty) this result was obtained by
Bowditch \cite{Bow-98}.

(V) Recently Dahmani \cite{Dah2} proved combination theorem for
relatively hyperbolic groups. This allowed to show that the limit
groups introduced by Sela \cite{Sela} in his solution of the
Tarski problem are hyperbolic in the sense of Bowditch relative to
their maximal abelian non--cyclic subgroups (see also \cite{Ali}).

We also mention some examples of groups, which are hyperbolic in
the sense of Farb.

(VI)  Let $Mod (S)$ denote the mapping class group corresponding
to a surface $S$. Using the isometric action of  $Mod (S)$ on the
complex of curves introduced in \cite{Harv}, Masur and Minsky
proved that $Mod(S)$ is hyperbolic relative to a finite collection
of stabilizers of certain curves in the sense of Farb \cite{MM}.
However, in most cases it is not hyperbolic in the sense of
Bowditch. An alternative proof can be found in \cite{Bow03}

(VII) Applying a technique related to small cancellation theory,
Kapovich \cite{Kap} proved the relative hyperbolicity in the sense
of Farb of some Artin groups of extra large type with respect to
certain families of parabolic subgroups. Another result of this
type was obtained by Bahls \cite{Bah}. He showed that right-angled
Coxeter groups are relatively hyperbolic in the sense of Farb with
respect to natural collections of rank $2$ parabolic subgroups.

(VIII) Finally we mention two combination theorems. If $G$ is an
HNN--extension of a group $H$ (respectively an amalgamated product
of $H_1$ and $H_2$) with associated (respectively amalgamated)
subgroups $A$ and $B$, then $G$ is hyperbolic in the sense of Farb
relative to $H$ (respectively relative to $\{ H_1, H_2\} $.
Moreover, if $H$ is hyperbolic in the sense of Farb relative to
$\{ A,B\} $ (respectively $H_1$ is hyperbolic in the sense of Farb
relative to $A$ and $H_2$ is hyperbolic in the sense of Farb
relative to $B$), then $G$ is hyperbolic in the sense of Farb
relative to $A$ \cite{Wh}. In particular, this allows to construct
a finitely presented group $G$ which is hyperbolic in the sense of
Farb relative to a hyperbolic subgroup $H$ and has undecidable
word problem.

\newpage


\section{Main results}


In this section we discuss shortly the main results of our paper.
We assume the reader to be familiar with such notions as Cayley
graph, Dehn function, hyperbolic group, quasi--geodesic path,
etc., and refer to the next chapter for the precise definitions.

Let $G$ be a group, $X$ a subset of $G$, and $H$ an arbitrary
subgroup of $G$. For simplicity we consider here the case of a
single subgroup and refer to Section 2.1 for the general case. We
say that $X$ is a {\it relative generating set} of $G$ with
respect to $H$, if $G$ is generated by the set $H\cup X$.

In this situation there exists a canonical homomorphism
$$
\varepsilon : F=H \ast F(X) \to G,
$$
where $F(X)$ is the free group with the basis $X$. If $Ker \;
\varepsilon $ is a normal closure of a subset $\mathcal R\subseteq
N$ in $F$, we say that $G$ has the {\it relative presentation}
\begin{equation}
\langle X,H\; |\; R=1,\, R\in\mathcal R\rangle \label{Fintr}
\end{equation}
The relative presentation (\ref{Fintr}) is {\it finite } if the
sets $X$ and $\mathcal R$ are finite. $G$ is said to be {\it
relatively finitely presented} with respect to $H$, if it admits a
finite relative presentation. Note that $G$ and $H$ need not be
finitely presented or even finitely generated in the usual sense.

Similarly one can define the notion of a relatively finitely
generated and a relatively finitely presented group with respect
to arbitrary collection of subgroups. We begin with the theorem,
which shows some restrictions in case $G$ is finitely generated
(in the usual sense) and is finitely presented with respect to a
collection of subgroups $\{ H_\lambda \} _{\lambda \in \Lambda } $
(see Proposition \ref{Hfg} and Corollary \ref{Lambda}).

\begin{thm}\label{2}
Let $G$ be a finitely generated group, $\{ H_\lambda \} _{\lambda
\in \Lambda }$ a collection of subgroups of $G$. Suppose that $G$
is finitely presented with respect to $\{ H_\lambda \} _{\lambda
\in \Lambda }$. Then the following conditions hold.

1) The collection of subgroups is finite, i.e. $card\; \Lambda
<\infty $.

2) Each subgroup $H_\lambda $ is finitely generated.
\end{thm}

\begin{defn}
Given a finite relative presentation (\ref{Fintr}), we say that
$f: \mathbb N\to \mathbb N$ is {\it a relative isoperimetric
function} of (\ref{Fintr}) if for every word $W$ of length $\|
W\|\le n$ over the alphabet $X^{\pm 1}\cup H$ representing the
identity in $G$, there exists an expression
\begin{equation}
W=_F\prod\limits_{i=1}^k f_i^{-1}R_if_i \label{intrprod}
\end{equation}
(with the equality in the group $F$), where $R_i\in \mathcal R$,
$f_i\in F$, and $k\le f(n)$. The minimal relative isoperimetric
function of (\ref{Fintr}) is called the {\it relative Dehn
function} of $G$ (with respect to $H$). We denote it by $\delta
^{rel}_{G,H}$.
\end{defn}

For instance, any finitely presented group $G$ is relatively
finitely presented with respect to the trivial subgroup. In this
case the corresponding relative Dehn function coincides with the
ordinary Denh function of $G$.

\begin{ex} \label{wd}
We stress that the relative Dehn function is not always
well--defined, i.e., it can be infinite for certain values of the
argument, since the number of words of bounded relative length can
be infinite. Indeed consider the group
$$G=\langle a,b \; |\; [a,b]=1 \rangle\cong \mathbb Z\times
\mathbb Z $$ and the cyclic subgroup $H$ generated by $a$. Clearly
one can take $X=\{ b\} $. It is easy to see that the word
$W_n=[a^n,b]$ has length $4$ as a word over $H\cup \{ b\} $ for
every $n$. However, the minimal number of factors in the
expression of type (\ref{intrprod}) corresponding to $W_n$ growths
linearly as $n\to\infty $. Thus we have $\delta ^{rel}_{G,H}
(4)=\infty $ for this relative presentation.
\end{ex}

The definition of the relative Dehn function in case of arbitrary
collection of subgroups $\{ H_\lambda \} _{\lambda \in \Lambda }$
is similar (see Section 2.3). Analogously to the case of ordinary
Dehn functions, if the relative Dehn function of $G$ with respect
to $\{ H_\lambda \} _{\lambda \in \Lambda }$ is finite for each
value of the argument, it is independent of the choice of a finite
relative presentation up to some equivalence relation (Theorem
\ref{Df}). Thus we can speak about the relative Dehn function of
the pair $G, \{ H_\lambda \}_{\lambda \in \Lambda }$ by means of
the corresponding equivalence classes. The next result is obtained
in Section 2.3 and provides us with necessary conditions for the
relative Dehn function to be well--defined.

\begin{thm}\label{3}
Let $G$ be a group, $\{ H_\lambda \} _{\lambda \in \Lambda }$ a
collection of subgroups of $G$. Suppose that $G$ is finitely
presented with respect to $\{ H_\lambda \} _{\lambda \in \Lambda
}$ and the relative Denh function of $G$ with respect to $\{
H_\lambda \} _{\lambda \in \Lambda })$ is well--defined. Then the
following conditions hold.

1) For any $g_1,g_2\in G$, the intersection $H_\lambda ^{g_1} \cap
H_\mu ^{g_2} $ is finite whenever $\lambda \ne \mu $.

2) The intersection $H_\lambda ^g \cap H_\lambda $ is finite for
any $g\not\in H_\lambda $.
\end{thm}

The main reason for our study of relative Dehn functions is the
following characterization of relative hyperbolicity of finitely
generated groups.

\begin{thm}\label{MainTh}
Let $G$ be a finitely generated group,  $\{ H_1, H_2, \ldots , H_m
\}$ a collection of subgroups of $G$. Then the following
conditions are equivalent.

1) $G$ is finitely presented with respect to $\{ H_1, H_2, \ldots
, H_m \}$ and the corresponding relative Dehn function is linear.

2) $G$ is hyperbolic with respect to the collection $\{ H_1, H_2,
\ldots , H_m \}$ in the sense of Farb and satisfies the BCP
property (or, equivalently, $G$ is hyperbolic relative to $\{ H_1,
H_2, \ldots , H_m \}$ in the sense of Bowditch).
\end{thm}

Theorem \ref{MainTh} allows to consider the definition below as a
generalization of Bowditch's one.

\begin{defn}\label{O}
We say that a group $G$ is {\it hyperbolic relative to a
collection of subgroups} $\{ H_\lambda \} _{\lambda \in \Lambda }$
if $G$ is finitely presented with respect to $\{ H_\lambda \}
_{\lambda \in \Lambda }$ and the relative Dehn function of $G$
with respect to $\{ H_\lambda \} _{\lambda \in \Lambda }$ is
linear.
\end{defn}

Standard arguments show that if the relative Dehn function is
subquadratic, then, in fact, it is linear (Corollary \ref{Ghyp}).

Given a group $G$, a collection of subgroups $\{ H_\lambda \}
_{\lambda \in \Lambda }$, and a finite relative generating set $X$
of $G$ with respect to $\{ H_\lambda \} _{\lambda \in \Lambda }$,
we denote by $\G $ the Cayley graph of $G $ with respect to the
generating set $X\cup\left( \bigcup\limits_{\lambda\in \Lambda }
H_\lambda \right) $.  It is easy to prove that if the relative
Dehn function of $G$ with respect to $\{ H_\lambda \} _{\lambda
\in \Lambda }$ is linear, then $\G $ is hyperbolic. A partially
converse result is obtained in Section 2.5.

\begin{thm}
Suppose that a group $G$ is finitely presented with respect to a
collection of subgroups $\{ H_\lambda \} _{\lambda \in \Lambda }$
and the relative Dehn function of $G$ with respect to $\{
H_\lambda \} _{\lambda \in \Lambda }$ is well--defined. Then
the following conditions are equivalent.

1) The group $G$ is hyperbolic relative to $\{ H_\lambda \}
_{\lambda \in \Lambda }$.

2) The Cayley graph $\G $ is a hyperbolic metric space.
\end{thm}

If the group $G$ is generated by a finite set $X$ in the ordinary
non--relative sense, we can think of the Cayley graph of $G$ with
respect to $X$, $\Gamma (G,X)$, as a subgraph of the graph $\Gamma
(G, X\cup \mathcal H)$ defined above. Obviously these graphs have
the same set of vertices. Assuming the lengths of any edge of
$\Gamma (G,X)$ and $\G $ to have length $1$, we get two
combinatorial metrics $\dx $ and $\dxh $ on $\Gamma (G,X)$ and $\G
$ respectively. The proof of the Theorem \ref{MainTh} given in the
appendix is based on some results about the geometry of the
embedding of $\Gamma (G, X)$ into $\G $, which are proved in
Chapter 3. Although the proofs of these facts take a significant
part of our paper, the theorems seem to be too technical to
formulate them here.

We also introduce and study the notion of relative
quasi--convexity for subgroups of $G$ in case $G$ is finitely
generated.

\begin{defn}\label{qc}
Let $G$ be a group generated by a finite set $X$, $\{ H_1, \ldots
, H_m\} $ a collection of subgroups of $G$.  A subgroup $R$ of $G$
is called {\it relatively quasi--convex} with respect to $\{ H_1,
\ldots , H_m\} $ if there exists a constant $\sigma >0$ such that
the following condition holds. Let $f$, $g$ be two elements of
$R$, and $p$ an arbitrary geodesic path from $f$ to $g$ in $\G $.
Then for any vertex $v\in p$ there exists a vertex $w\in R$ such
that
$$\dx (u,w)\le \sigma .$$ In case $G$ is hyperbolic relative to $\{ H_1, \ldots
, H_m\} $, this notion is independent of the choice of the finite
generating set of $G$ (Proposition \ref{qcinv}). A subgroup $R$ is
called {\it strongly relatively quasi--convex} if, in addition,
the intersections $R\cap H_i^g$ are finite for all $i=1,\ldots
,m$, $g\in G$.
\end{defn}

The next three theorems are the relative analogues of well--known
facts about quasi--convex subgroups of ordinary hyperbolic groups
(see Section 4.2). In these theorems we suppose $G$ to be a
finitely generated group hyperbolic relative to a collection of
subgroups $\{ H_1, \ldots , H_m\}$.

\begin{thm}
Let $R$ a subgroup of $G$. Then the following conditions are
equivalent.

1) $R$ is strongly relatively quasi--convex.

2) $R$ is generated by a finite set $Y$ and the natural embedding
$(R, dist _Y) \to (G, \dxh )$ is a quasi--isometry.
\end{thm}

As a corollary, we obtain

\begin{thm}
If $R$ is a strongly relatively quasi--convex subgroup of $G$,
then $R$ is a hyperbolic group.
\end{thm}

Finally we show that the set of strongly relatively quasi--convex subgroups
is closed under intersections.

\begin{thm}\footnote{In the published version of the paper, this theorem was stated for quasi--convex subgroups instead of strongly quasi--convex. Although the statement remains true for quasi-convex subgroups (see \cite{Hru}), it is not proved here. We only prove it in the strongly relatively quasi-convex case.}  
Let $P$ and $R$ be two strongly relatively quasi--convex subgroups of $G$.
Then $P\cap R$ is strongly relatively quasi--convex.
\end{thm}

Let us mention some applications of the technique developed in
this paper to the study of algebraic and algorithmic properties of
relatively hyperbolic groups. We say that an element $g\in G$ is
{\it hyperbolic} if $g$ is not conjugate to an element of one of
the subgroups $ H_\lambda $, $\lambda \in \Lambda $. In the next
two theorems and their corollaries we assume that $G$ is an
arbitrary (not necessary finitely generated) group that is
relatively hyperbolic with respect to a collection of subgroups
$\{ H_\lambda \} _{\lambda \in \Lambda }$. The proofs can be found
in Sections 4.1, 4.3.

\begin{thm} \label{101}
There exist only finitely many conjugacy classes of hyperbolic
elements of finite order in $G$. In particular, the set of orders
of hyperbolic elements is finite.
\end{thm}

\begin{cor}
If $G$ is residually finite and all subgroups $H_\lambda $ are
torsion free, then $G$ is virtually torsion free, that is $G$
contains a torsion free subgroup of finite index.
\end{cor}

It is worth to notice that the assumption of residually finiteness
is essential in the last corollary. Indeed in Section 4.1 we
construct an example of a (finitely generated) group $G$
hyperbolic relative to a torsion--free subgroup $H$ such that $G$
is not virtually torsion--free.

\begin{thm}\label{10}
For any hyperbolic element $g\in G$ of infinite order, there exist
$\lambda >0$, $c\ge 0$ such that $$\dxh (1, g^n)>\lambda |n|-c$$
for any $n\in \mathbb N$.
\end{thm}

It follows from the proof of Theorem \ref{10} that every cyclic
subgroup generated by a hyperbolic element of infinite order has
finite index in its centralizer. In other terms this was first
proved by Tukia in \cite{Tuk}. We also show that the constant
$\lambda $ in Theorem \ref{10} can be chosen independently of $g$.

\begin{cor}
If $g$ is a hyperbolic element of $G$ of infinite order and $g^k$
is conjugate to $g^l$ for some $k,l\in \mathbb Z$, then $k=\pm l$
\end{cor}

We note that the existence of an action of $G$ on a hyperbolic
graph such that all edge stabilizers are trivial and every vertex
stabilizer is conjugate to $H_i$ for some $i$ (that is a 'weak'
form of hyperbolicity of $G$) is not sufficient for the fulfilment
of Theorems \ref{101} and \ref{10}. The counterexamples are
provided in Sections 4.1, 4.3.

It is known that the word and the conjugacy problems are decidable
in a (finitely generated) group $G$ hyperbolic with respect to $\{
H_1, \ldots , H_m\}$ provided these problems are decidable for
each $H_i$ \cite{Bum}, \cite{F}. Also if the conjugacy problem is
solvable in $H_1, \ldots , H_m$, then given $g\in G$ and $i\in \{
1, \ldots , m\} $, it is possible to decide whether $g$ is
conjugate to an element of $H_i$ \cite{Bum}. Some of these results
are generalized in Sections 5.1, 5.2.

The application of our approach to the study of other algorithmic
problems leads to the following theorem.

\begin{thm}
Let $G$ be a group hyperbolic with respect to a collection of
recursively presented subgroups $\{ H_1, \ldots , H_m\}$.  Then
each of the algorithmic problems listed below is solvable in $G$
whenever the word problem is solvable in each of the subgroups
$H_1, \ldots , H_m$.

1) The conjugacy problem for hyperbolic elements, that is, given
two hyperbolic elements $f,g\in G$, to decide whether $f$ and $g$
are conjugate.

2) The order problem for hyperbolic elements, that is to calculate
the order of a given hyperbolic element $g\in G$.

3) The root problem for hyperbolic elements, that is, given an
hyperbolic element $g\in G$, to decide whether there exists a
nontrivial root of $g$ in $G$.

4) The power conjugacy problem for hyperbolic elements, that is to
decide, for any two hyperbolic elements $g,f\in G$, whether or not
there exist two hyperbolic conjugate powers of $f$ and $g$.
\end{thm}

\paragraph{Acknowledgements.}
The paper was written in part during the author's visits in
Courant Institute of Mathematical Sciences, in spring 2002, and
Vanderbilt University, in fall 2002. I would like to thank Brian
Bowditch, Benson Farb, Mihail Gromov, Alexander Ol'shanskii, and
Mark Sapir for useful conversations and remarks.


\chapter{Relative isoperimetric inequalities}



\section{Relative presentations and length functions}


We begin by introducing relative generating sets of groups with
respect to fixed collection of subgroups.

\begin{defn}
Let $G$ be a group, $\{ H_\lambda \} _{\lambda\in \Lambda} $ a
collection of subgroups of $G$, $X$ a subset of $G$. We say that
$G$ is {\it generated by $X$ with respect to} $\{ H_\lambda \}
_{\lambda\in \Lambda } $ (or, equivalently, $X$ is a {\it relative
generating set of $G$ with respect to $\{ H_\lambda \}
_{\lambda\in \Lambda }$}) if $G$ is generated by the set
$\left(\bigcup\limits_{\lambda \in \Lambda} H_\lambda \right) \cup
X$. We will always assume the set $X$ to be symmetrized, i.e.,
$X=X^{-1}$.
\end{defn}

In the above situation the group $G$ can be regarded as the
quotient group of the free product
\begin{equation}
F=\left( \ast _{\lambda\in \Lambda } \widetilde H_\lambda  \right)
\ast F(X), \label{F}
\end{equation}
where the groups  $\widetilde H_\lambda $ are isomorphic copies of
$H_\lambda$, and $F(X)$ is the free group with the basis $X$. Let
us denote by $\mathcal H$ the disjoint union
$$ \mathcal H=\bigsqcup\limits_{\lambda\in \Lambda} (\widetilde
H_\lambda \setminus \{ 1\} ) .$$ It is easy to see that $F$ is
generated by $X\cup \mathcal H$.

{\bf Conventions and notation.} By $(X\cup \mathcal H)^\ast $, we
denote the free monoid generated by $X\cup \mathcal H$. Given
$W\in (X\cup \mathcal H)^\ast $, $\| W\| $ denotes the length of
the word $W$, $\overline{W}$ denotes the element of $G$
represented by $W$. Throughout this paper we write $U\equiv V$ to
express letter--for--letter equality of two words $U,V\in (X\cup
\mathcal H)^\ast $ and write $U=_FV$ when $U$ and $V$ represent
the same elements of the group $F$. To simplify our notation we
identify words from $(X\cup \mathcal H)^\ast $ and elements of the
group $F$ represented by them. We also write $x^t$ for $t^{-1}xt$
and $[x,t]$ for $x^{-1}t^{-1}xt$.

For every $\lambda \in \Lambda $, we denote by $\mathcal S_\lambda
$ the set of all words over the alphabet $\widetilde H_\lambda
\setminus \{ 1\} $ that represent the identity in $F$. Thus the
group $F$ can be defined by the presentation
\begin{equation}
\langle X, \mathcal H \; | \; S=1, S\in \bigcup\limits_{\lambda
\in \Lambda } \mathcal S_\lambda \rangle .
\end{equation}
The isomorphisms $\widetilde H_\lambda $ and the identity map on
$X$ can be uniquely extended to a homomorphism $\epsilon :F\to G$.
We denote its kernel by $N$.

\begin{defn}[Relative presentation]
We say that the group $G$ has the {\it relative presentation}
\begin{equation}
\langle X, \mathcal H\; | \; S=1, S\in \bigcup\limits_{\lambda \in
\Lambda } \mathcal S_\lambda , R=1, R\in \mathcal R\rangle,
\label{G1}
\end{equation}
with respect to the collection of subgroups $\{ H_\lambda\}
_{\lambda \in \Lambda } $, where $\mathcal R\subseteq (X\cup
\mathcal H)^\ast $, if $N$ is the normal closure of the set
$\mathcal R$ in the group $F$. It is convenient to assume that
$\mathcal R$ is symmetrized that is for every $R\in \mathcal R$,
the set $\mathcal R$ contains all cyclic shifts of $R$ and
$R^{-1}$.
\end{defn}

For brevity, we use the following reduced record for the
presentation (\ref{G1})
\begin{equation}
\langle X, H_\lambda, \lambda\in \Lambda \; | \;  R=1, R\in
\mathcal R\rangle, \label{G2}
\end{equation}

\begin{defn}
The relative presentation (\ref{G2}) is called {\it finite } if
both the sets $\mathcal R$ and $X$ are finite. If there exists a
finite relative presentation of a group $G$ with respect to a
collection of subgroups $\{ H_\lambda \} _{\lambda\in \Lambda }$,
we say that $G$ is {\it finitely presented} relative to $\{
H_\lambda \} _{\lambda\in \Lambda }$.
\end{defn}

\begin{ex}
Consider the amalgamated product $$G=H_1 \ast _{K=L}H_2$$ of two
arbitrary groups $H_1$, $H_2$ associated to an isomorphism $\alpha
: K\to L$ between subgroups $K\le H_1 $ and $L\le H_2$. Then $G$
has the relative presentation $$\langle H_1, H_2\; | \; k=\alpha
(k), \; k\in K\rangle $$ with respect to $\{ H_1, H_2\} $. If $K$
is finitely generated, one can construct a finite relative
presentation for $G$ with respect to $\{ H_1, H_2\}$. Indeed it is
sufficient to impose the relations of type $k=\alpha (k)$ for all
generators of $K$.
\end{ex}

\begin{ex}
If $H$ is a group and $\alpha :A\to B$ is an isomorphism between
two subgroups $A, B\le H$, then the corresponding HNN--extension
has the relative presentation $$\langle H\; | \; a^t=\alpha (a),
a\in A\rangle . $$ As above this relative presentation can be made
finite in case $A$ is finitely generated.
\end{ex}

\begin{defn}
To each element $g\in G$, we assign its {\it relative length with
respect to the collection of subgroups} $\{ H_\lambda \} _{\lambda
\in \Lambda }$ of $G$ (or simply {\it relative length})
$|g|_{X\cup \mathcal H}$ that is the length of a shortest word
from $(X\cup \mathcal H)^\ast $ representing $g$ in $G$. This
length function induces the left--invariant {\it relative distance
function} on $G\times G$ by the rule $$dist _{X\cup \mathcal H}
(g,h)=|g^{-1}h|_{X\cup \mathcal H}.$$
\end{defn}

We note that the group $G$ endowed with the relative metric is not
always a proper metric space, i.e., the closed balls can be
infinite. It is easy to see that $G$ is proper if and only if
$card\, \Lambda <\infty $ and all subgroups $H_\lambda $ are
finite.

It is clear that this distance strongly depends on the choice of
the set $X$ and the collection of subgroups. However, if the
collection of subgroups $\{ H_\lambda \} _{\lambda \in \Lambda }$
is fixed and $G$ is finitely generated with respect to $\{
H_\lambda \} _{\lambda \in \Lambda }$, then the asymptotic
behavior of the distance function is essentially independent of
the choice of a finite relative generating set. This can be
expressed in the following way.

\begin{defn}
Two metrics $dist_1$ and $dist_2$ on the same space are called
{\it Lipschitz equivalent} if the ratios $dist_1/dist_2$ and
$dist_2/dist_1$ are bounded when they are considered as functions
on the Cartesian square of the space minus the diagonal.
\end{defn}

The proposition below is the relative analogue of the well known
property of ordinary word metrics on finitely generated groups.

\begin{prop}\label{Lip}
Let $G$ be a group, $\{ H_\lambda \} _{\lambda \in \Lambda }$ a
collection of subgroups of $G$. Suppose that $X$ and $Y$ are two
finite relative generating sets of $G$ with respect to $\{
H_\lambda \} _\Lambda $. Then the corresponding distance functions
$ dist _{X\cup \mathcal H}$ and $dist _{Y\cup \mathcal H}$ are
Lipschitz equivalent.
\end{prop}

\begin{proof}
Let us take an arbitrary element $g\in G$. Suppose that $W\in
(Y\cup \mathcal H)^\ast $ is a shortest word representing $g$. For
every $y\in Y$, we fix a word $V_y\in (X\cup \mathcal H)^\ast $
that represents $y$. Put $$ M=\max\limits_{y\in Y} \| V_y\| .$$ As
the set $Y$ is finite, we have $M< \infty $.

If we consider the word $U\in (X\cup \mathcal H)^\ast $ that is
obtained from $W$ by replacing all letters $y\in Y$ with the
corresponding words $V_y$, then $$\| U\| \le M\| W\| = M
|g|_{Y\cup \mathcal H}.$$ Finally we have $$|g|_{X\cup \mathcal
H}\le \| U\| \le M |g|_{Y\cup \mathcal H}.$$ The reverse
inequality can be obtained in the analogous way. The proposition
is proved.
\end{proof}


\section{Geometry of van Kampen diagrams over relative
presentations}


Our study of properties of relatively hyperbolic groups is based
on the combinatorial geometry of van Kampen diagrams over relative
representations. In this section we collect definitions and some
technical facts about the diagrams which are used in what follows.

\begin{defn}
Recall that a planar {\it map} $\Delta $ over a presentation
\begin{equation}
\langle Z\; | \; \mathcal P\rangle \label{ZP}
\end{equation}
is a finite oriented connected simply--connected 2--complex
endowed with a labelling function $\phi : E(\Delta )\to Z^{\pm
1}$, where $E(\Delta ) $ denotes the set of oriented edges of
$\Delta $, such that $\phi (e^{-1})=(\phi (e))^{-1}$. The label of
a path $p=e_1\ldots e_n$ is, by definition, the word
$\phi(e_1)\ldots \phi (e_n)$. By {\it length} $l(p)$ we mean the
number of edges in $p$. Given a cell $\Pi $ of $\Delta $, we
denote by $\partial \Pi$ the boundary of $\Pi $; similarly,
$\partial \Delta $ denotes the boundary of $\Delta $. The labels
of $\partial \Pi $ and $\partial \Delta $ are defined up to a
cyclic permutation.
\end{defn}

\begin{defn}
A map $\Delta $ over a presentation (\ref{ZP}) is called a {\it
van Kampen diagram} over (\ref{ZP}) if the following holds. For
any cell $\Pi $ of $\Delta $, the boundary label $\phi (\partial
\Pi)$ is equal to a cyclic permutation of a word $P^{\pm 1}$,
where $P\in \mathcal P$. Sometimes it is convenient to use the
notion of $0$--refinement in order to assume diagrams to be
homeomorphic to a disc. We do not explain here this notion and
refer the interested reader to \cite[Ch. 4]{Ols-book}.
\end{defn}

The van Kampen lemma states that a word $W$ over the alphabet
$Z^{\pm 1}$ represents the identity in the group given by
(\ref{ZP}) if and only if there exists a simply--connected planar
diagram $\Delta $ over (\ref{ZP}) such that $\phi (\partial \Delta
)\equiv W$ \cite{LS}.

Dealing with diagrams over relative presentations we will divide
the set of 2--cells into two parts as follows.

\begin{defn}
Let $\Delta $ be a van Kampen diagram over the relative
presentation (\ref{G1}).  A cell $\Pi $ of $\Delta $ is called an
{\it $\mathcal S_\lambda $--cell} if it corresponds to a relator
from $\mathcal S_\lambda $. A cell $\Pi $ is called an $\mathcal
S$--cell if it is an $\mathcal S_\lambda $--cell for some
$\lambda\in \Lambda $. Obviously for any $\mathcal S$--cell, the
label of its boundary represents $1$ in the group $F$. Similarly
we call $\Pi $ an $\mathcal R$--cell if it corresponds to a
relator from $\mathcal R$.
\end{defn}

\begin{defn}
Given a van Kampen diagram $\Delta $ over (\ref{G1}), we denote by
$N_\mathcal R (\Delta )$ (respectively by $N_\mathcal S (\Delta
)$) the number of $\mathcal R$--cells (respectively $\mathcal
S$--cells) of $\Delta $. We define the area of $\Delta $ by the
formula $$ Area (\Delta )= N_\mathcal R(\Delta )+N_\mathcal
S(\Delta ). $$
\end{defn}

\begin{defn}
The {\it type} of the diagram $\Delta $ is the triple $\tau
(\Delta )= (N_\mathcal R (\Delta ), N_\mathcal S (\Delta ), card\,
E(\Delta ))$. We fix the lexicographic order on the set of all
triples, that is,
$$
\begin{array}{cc}
(a_1,b_1,c_1)\le (a_2,b_2,c_2)  \Longleftrightarrow \\ \\
(a_1< a_2)\; \vee\; (a_1 =a_2\, \wedge\, b_1<b_2)\; \vee \;
(a_1=a_2\, \wedge\, b_1=b_2\, \wedge \, c_1\le c_2).
\end{array}
$$
\end{defn}

\begin{figure}
\label{T}
\begin{picture}(120,25)(0,-3)
\put(40,10){\oval(30,20)} \put(80,10){\oval(20,20)}
\put(55,10){\line(1,0){15}} \put(35,0){\line(0,1){20}}
\put(21,9){$e_2$} \put(45,1){$e_3$} \put(45,17){$e_1$}
\put(36,9){$d$} \put(62,11){$e_5$} \put(91,9){$e_4$}
\end{picture}
\caption{Types of edges in diagrams.}
\end{figure}
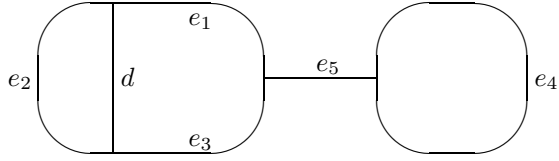

We will exploit the following classification of edges of van
Kampen diagrams over group presentations.

\begin{defn}\label{edges}
A (non--oriented) edge $e$ of a diagram $\Delta $ is called {\it
internal} if it is a common edge of two cells of $\Delta $;
otherwise $e$ is called {\it external}. Thus the union of all
external edges forms the boundary $\partial \Delta $ of the
diagram $\Delta $. Further an external edge $e$ is called an {\it
external edge of the first type} if $e$ belongs to the boundary of
some cell of the diagram; otherwise $e$ is called an {\it external
edge} of the second type.
\end{defn}

The notions introduced in Definition \ref{edges} are illustrated
on Fig. 2.1. The edges $e_1, \ldots , e_5$ are external, the edge
$d$ is internal. More precisely, $e_1, \ldots , e_4$ are external
edges of the first type and $e_5$ is external of the second type.

\begin{lem}\label{DiagMinType}
Let $\Delta $ be a van Kampen diagram over (\ref{G1}). Suppose
that $\Delta $ has the minimal type among all van Kampen diagrams
over (\ref{G1}) with the same boundary label. Then every internal
edge of $\Delta $ belongs to the boundary of some $\mathcal
R$--cell of $\Delta $.
\end{lem}

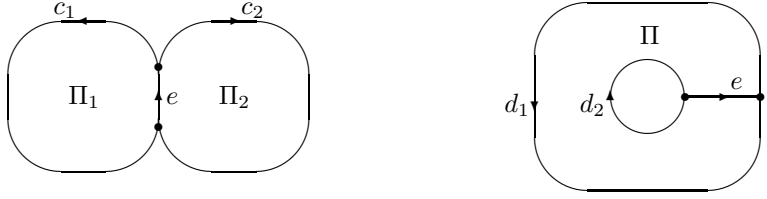
\begin{figure}

\begin{picture}(120,30)(-10,0)
\multiput(10,15)(20,0){2}{\oval(20,20)}
\multiput(20,11)(0,8){2}{\circle*{1}}

\put(8,14){$\Pi _1$} \put(28,14){$\Pi _2$}

\put(9,25){\vector(-1,0){0}} \put(31,25){\vector(1,0){0}}
\put(6,26){$c_1$} \put(31,26){$c_2$} \put(21,14){$e$}
\put(20,16){\vector(0,1){0}}

\put(85,15){\circle{10}} \put(85,15){\oval(30,25)}
\put(90,15){\circle*{1}} \put(100,15){\circle*{1}}
\put(90,15){\line(1,0){10}} \put(96,15){\vector(1,0){0}}

\put(70,13){\vector(0,-1){0}} \put(80,16){\vector(0,1){0}}
\put(66,13){$d_1$} \put(76,13){$d_2$} \put(96,16){$e$}
\put(84,22){$\Pi $}
\end{picture}

\caption{Two cases in the proof of Lemma \ref{DiagMinType}.}
\label{DMTfig}
\end{figure}

\begin{proof}
Let $e\in E(\Delta )$ be an internal edge. Suppose that $e$ does
not belong to the boundary of any $\mathcal R$--cell of $\Delta $.
Then there are two possibilities (see Fig. \ref{DMTfig}): either
$e$ is a common edge of two distinct $S_\lambda $--cells $\Pi _1 $
and $\Pi _2$ or there is an $S_\lambda $--cell $\Pi $ in $\Delta $
such that $$\partial \Pi =ed_1e^{-1}d_2, $$ where $d_1$ and $d_2$
are cycles in $\Delta $ (one of them may be trivial) and one of
the cycles $d_1, d_2$, say $d_2$, is contained in the part of the
diagram bounded by the other one. In the first case let
$$\partial \Pi _1 =ec_1, \;\;\;\;\; \partial \Pi _2 =
c_2^{-1}e^{-1}.$$

Note that the label of the path $c_1c_2^{-1}$ (respectively $d_1$)
consists of letters from $\widetilde H_\lambda $ and represents
$1$ in $F$. Since $\mathcal S_\lambda $ contains all words over
$\widetilde H_\lambda $ that represents the identity in the group
$F$, we can replace the cells $\Pi _1$, $\Pi _2$ (respectively the
cell $\Pi $) with the cell $\Sigma $ with the boundary
$c_1c_2^{-1}$ (respectively $d_1$). In both cases the type of the
diagram decreases and we get a contradiction.
\end{proof}

\begin{cor} \label{sums}
Let $\Delta $ be a van Kampen diagram over (\ref{G1}). By
$S(\Delta )$ we denote the set of all $\mathcal S$--cells of
$\Delta $. Suppose that (\ref{G1}) is relatively finite and
$\Delta $ has the minimal type among all van Kampen diagrams over
(\ref{G1}) with the same boundary label. Then we have
\begin{equation} \label{sums1}
\sum\limits_{\Pi \in S(\Delta )} l (\partial \Pi )\le M
N_{\mathcal R}(\Delta ) + l(\partial \Delta),
\end{equation}
where $M=\max\limits_{R\in \mathcal R} \| R\| $.
\end{cor}

\begin{proof}
Let $R(\Delta )$ denote the set of all $\mathcal R$--cells of
$\Delta $. According to Lemma \ref{DiagMinType}, the number of all
(non--oriented) edges of $\Delta $ satisfies the inequality $$
card\; E(\Delta ) \le \sum\limits_{\Xi \in R(\Delta )} l(\partial
\Xi ) +l(\partial \Delta )\le M N_{\mathcal R}(\Delta ) +
l(\partial \Delta).$$ This implies (\ref{sums1}) as every edge of
$\Delta $ belongs to the boundary of at most one $\mathcal
S$--cell.
\end{proof}

\begin{defn}
The {\it Cayley graph} $\Gamma = \Gamma (G, Z)$ of a group $G$
generated by a set $Z$ is an oriented labelled 1--complex with the
vertex set $V(\Gamma )=G$ and the edge set $E(\Gamma )=G\times Z$.
An edge $e=(g,s)\in E(\Gamma )$ goes from the vertex $g$ to the
vertex $gs$ and has the label $\phi (e)=s$. As usual, we denote
the origin and the terminus of the edge $e$, i.e., the vertices
$g$ and $gs$, by $e_-$ and $e_+$ respectively.

Given a combinatorial path $p=e_1e_2\ldots e_k$ in $\Gamma $,
where $e_1, e_2, \ldots , e_k\in E(\Gamma )$, we denote by $\phi
(p)$ its label. As in the case of diagrams, $\phi
(p)=\phi(e_1)\phi(e_2)\ldots \phi (e_k).$ By $p_-=(e_1)_-$ and
$p_+=(e_k)_+$ we denote the origin and the terminus of $p$
respectively. A path $p$ is called {\it irreducible} if it
contains no subpaths of type $ee^{-1}$ for $e\in E(\Gamma )$. A
{\it subpath} $q$ of $p=e_1e_2\ldots e_k$ is a path of type
$e_ie_{i+1}\ldots e_j$ for some $1\le i\le j\le k$. (So we always
assume $p$ and $q$ to have compatible orientations, that is,
starting from $p_-$ and passing along $p$, we first meet the
vertex $q_-$ and then $q_+$.) The graph $\Gamma $ can be regarded
as a metric space if we endow it with a combinatorial metric. This
means that the length of every edge of $\Gamma $ is assumed to be
equal to 1.
\end{defn}

In this paper we denote by $\G $ the Cayley graph of the group $G$
given by (\ref{G1}) with respect to the generating set $X\cup
\mathcal H$. Note that every van Kampen diagram $\Delta $ over
(\ref{G1}) can be mapped to $\G $ in such a way that the mapping
preserves labels and orientation. Taking into account this remark,
we will often consider configurations in the Cayley graph $\G $
instead of the treating the corresponding van Kampen diagrams.

Now we are going to introduce an auxiliary terminology, which
plays an important role in our paper.

\begin{defn}[$H_\lambda $--subwords]
Given a word $W\in (X\cup \mathcal H)^\ast $, we say that a
subword $V$ of $W$ is an {\it $H_\lambda $--subword} if $V$
consists of letters from $\widetilde H_\lambda $. An $H_\lambda
$--subword of $W$ is called an {\it $H_\lambda $--syllable} if it
is maximal, i.e., it is not contained in a bigger $H_\lambda
$--subword of $W$.

As usual, by a cyclic word $W$ we mean the set of all cyclic
shifts of $W$. As in the case of ordinary words, we say that
subword $V$ of a cyclic word $W$ is a $H_\lambda $--subword if it
is an $H_\lambda $--subword of a certain cyclic shift of $W$. A
maximal $H_\lambda $--subword of a cyclic word $W$ is called an
$H_\lambda $--syllable.
\end{defn}

\begin{defn}[$H_\lambda $--components]
Let $q$ be a path (respectively cyclic path) in $\G $. A subpath
$p$ of $q$ is called an {\it $H_\lambda $--subpath}, if the label
of $p$ is an $H_\lambda $--subword of the word $\phi (q)$
(respectively cyclic word $\phi (q)$). A {\it component} (or more
precisely an {\it $H_\lambda $--component}) of $q$ is an
$H_\lambda $--subpath $p$ such that the label of $p$ is an
$H_\lambda $--syllable of the the word $\phi (q)$ (respectively
cyclic word $\phi (q)$).
\end{defn}

\begin{defn}[Connected components]
Two $H_\lambda $--components $p_1, p_2$ of a path $q$ (cyclic or
not) in $\G $ are called {\it connected} if there exists a path
$c$ in $\G $ that connects some vertex of $p_1$ to some vertex of
$p_2$ and ${\phi (c)}$ is a word consisting of letters from $
\widetilde H_\lambda $. The path $c$ is called an {\it $H_\lambda
$--connector}. Note that this is equivalent to the requirement
that for any two vertices $v_1$ and $v_2$ of $p_1$ and $p_2$
respectively there exists a connector $c$ such that $c_-=v_1,
c_+=v_2$. (In algebraic terms this means that these two vertices
belong to the same coset $gH_\lambda $.) Clearly we can always
assume that $c$ consists of a single edge, as every element of
$H_\lambda $ is included in the set of generators.

Sometimes we will also speak about connected $H_\lambda
$--syllables in a word $W\in (X\cup \mathcal H)^\ast $ (cyclic or
not). By these we mean two $H_\lambda $--syllables $U,V$ of $W$
such that the corresponding components of some (or, equivalently,
of any) path in $\G $ labelled $W$ are connected.
\end{defn}

\begin{defn}[Isolated components]
An $H_\lambda $--component $p$ of a path $q$ (cyclic or not) is
called {\it isolated } if no (distinct) $H_\lambda $--component is
connected to $p$. The notion of an isolated $H_\lambda $--syllable
of a word $W\in (X\cup \mathcal H)^\ast $ is defined in the
obvious way.
\end{defn}

\begin{figure}
\begin{picture}(120, 20)(-10,-5)
\multiput(10,0)(10,0){8}{\circle*{1}}
\multiput(37,10)(6,0){2}{\circle*{1}} \thicklines{
\put(10,0){\vector(1,0){6}} \put(16,0){\vector(1,0){10}}
\put(26,0){\line(1,0){4}} \put(50,0){\vector(1,0){6}}
\multiput(56,0)(10,0){2}{\vector(1,0){10}}
\put(76,0){\line(1,0){4}} \put(40,0){\oval(20,20)[t]}}

\put(41,10){\vector(1,0){0}}
\multiput(30,1)(20,0){2}{\vector(0,-1){0}}

\thinlines \put(36,0){\vector(1,0){0}} \put(46,0){\vector(1,0){0}}

\qbezier[30](30,0)(40,0)(50,0)

\multiput(14,-3)(10,0){7}{$a$} \put(39,12){$a$}
\multiput(27,5)(25,0){2}{$t$}
\end{picture}
\caption{The path corresponding to the word $W\equiv a^2a^ta^3$.}
\label{bs}
\end{figure}
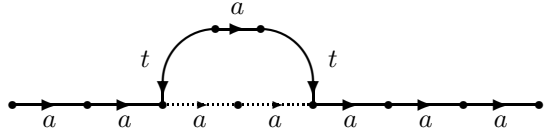

\begin{ex}\label{comp}
Let us consider the Baumslag--Solitar group $$BS(1,2)=\langle a,t
\; | \; a^t=a^2\rangle .$$ Set $H=\langle a\rangle $ and consider
the word $W\equiv a^2a^{t}a^3$. Then the $H$--syllables $a^3$ and
$a^2$ are connected since $a^{t}$ represents the same element of
$G$ as $a^2$ (see Fig. \ref{bs}). The $H$-syllable $a$ of $W$ is
isolated.
\end{ex}

In the next sections we will often use the following result
without any references. The proof follows from the definitions in
the obvious way.

\begin{lem}
Let $p$ be a geodesic path in the Cayley graph $\G $. Then every
component of $p$ is isolated.
\end{lem}

\begin{defn}\label{redpres}
A relative representation (\ref{G1}) of a group $G$ with respect
to a collection of subgroups $\{ H_\lambda \} _{\lambda \in
\Lambda }$ is called {\it reduced } if each relator  $R\in
\mathcal R$ has minimal length among all words from $(X\cup
\mathcal H)^\ast $ representing the same element of the group $F$.
In particular this means that for any $\lambda \in \Lambda $ and
any $R\in \mathcal R$, every $H_\lambda $--syllable of $R$
consists of a single letter. Without loss of generality we may
assume all relative presentations under consideration to be
reduced.
\end{defn}

The next lemma shows that, without loss of generality, we can
assume finite relative presentations to be reduced.

\begin{defn} \label{defofOmega}
For every $\lambda \in \Lambda $, we denote by $\Omega _\lambda $
the subset of all elements $g\in H_\lambda $ such that there
exists a relator $R\in \mathcal R$, and an $H_\lambda $--syllable
$V$ of $R$ such that $V$ represents $g$ in $G$. We also put
$$\Omega =\bigcup\limits_{\lambda \in \Lambda} \Omega_\lambda .$$
It is important that the set $\Omega $ is finite, whenever the set
$\mathcal R$ is finite.
\end{defn}

\begin{defn}[Relative area of a cycle in the Cayley graph]
Let $q$ be a cycle in the Cayley graph $\Gamma (G, X\cup\mathcal
H)$. We define its {\it relative area}, $Area^{rel}(q)$, with
respect to the presentation (\ref{G1}) as the minimal number
$N_\mathcal R(\Delta )$ of $\mathcal R$--cells among all van
Kampen diagrams $\Delta $ over (\ref{G1}) with the boundary label
$\phi (\partial \Delta )\equiv \phi (q)$.
\end{defn}

Now we are ready to formulate the main lemma of this section.
Although Lemma \ref{Omega} is quite trivial, it allows to obtain
some important results on groups given by finite relative
presentations (see, for example, Proposition \ref{Hfg} and Theorem
\ref{malnorm} below). Recall that for a word $W\in (X\cup\mathcal
H)^\ast $, $\overline{W} $ means the element of $G$ represented by
$W$.

\begin{lem}\label{Omega}
Suppose that a group $G$ is given by the reduced finite relative
presentation (\ref{G1}) with respect to a collection of subgroups
$\{ H_\lambda \} _{\lambda \in \Lambda }$.  Let $q$ be a cycle in
$\G $, $p_1, \ldots , p_k$ a certain set of isolated $H_\lambda
$--components of $q$. Then
\begin{equation}\label{om1}
\overline{\phi (p_i)} \in \langle \Omega _\lambda \rangle
\end{equation}
for any $i=1, \ldots , k$. Moreover, the lengths of the elements
$\overline{\phi (p_1)} , \ldots , \overline{\phi (p_k)} $ with
respect to the generating set $\Omega _\lambda $ of the subgroup
$\langle \Omega _\lambda \rangle $ satisfy the inequality
\begin{equation}\label{om2}
\sum\limits_{i=1}^k|\overline{\phi (p_i)} |_{\Omega _\lambda }\le
M\cdot Area^{rel} (q),
\end{equation}
where $$M=\max\limits_{R\in \mathcal R} ||R||.$$
\end{lem}

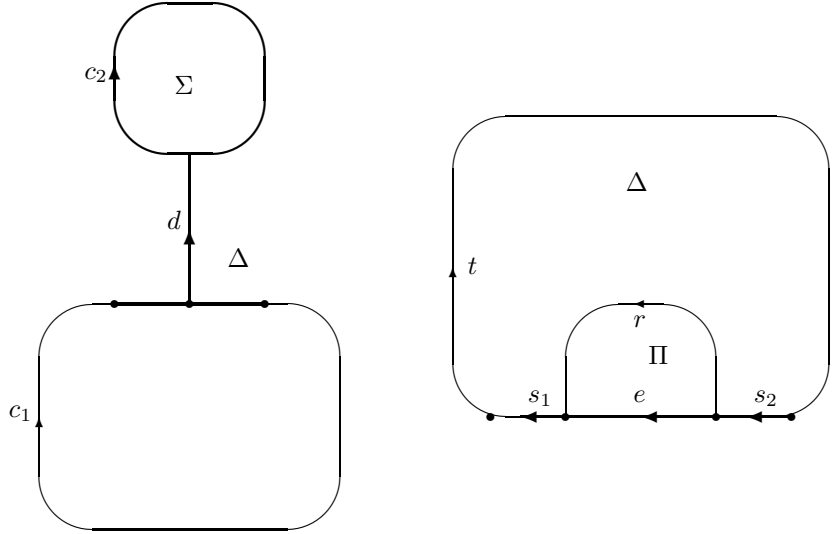
\begin{figure}
\begin{picture}(120,80)(0,-2)

\put(30,20){\oval(40,30)} \put(10,20){\vector(0,1){0}}
\put(6,20){$c_1$} \put(16,65){$c_2$} \put(27,45){$d$}
\put(28,63){$\Sigma $} \put(35,40){$\Delta $}

\thicklines { \put(30,65){\oval(20,20)}
\put(30,35){\line(0,1){20}} \put(30,35){\circle*{1}}
\put(20,35){\circle*{1}} \put(40,35){\circle*{1}}
\put(30,45){\vector(0,1){0}} \put(20,67){\vector(0,1){0}}
\put(20,35){\line(1,0){20}} }

\thinlines

\put(90,40){\oval(50,40)}  \put(90,20){\oval(20,30)[t]}
\put(88,50){$\Delta $} \put(65,40){\vector(0,1){0}}
\put(89,35){\vector(-1,0){0}}

\thicklines {\put(80,20){\vector(-1,0){6}}
\put(110,20){\vector(-1,0){6}} \put(76,20){\line(1,0){34}}
\put(90,20){\vector(-1,0){0}} }

\put(70,20){\circle*{1}} \put(80,20){\circle*{1}}
\put(100,20){\circle*{1}} \put(110,20){\circle*{1}}

\put(89,22){$e$} \put(75,22){$s_1$} \put(105,22){$s_2$}
\put(89,32){$r$} \put(67,39){$t$} \put(91,27){$\Pi $}
\end{picture}
\caption{Two cases in the proof of Lemma \ref{Omega}; the
component $p_i^\prime $ is marked by the thick line.}
\label{indstep}
\end{figure}

\begin{proof}
Changing the enumeration of the components $p_1, p_2, \ldots ,
p_k$ if necessarily, we may assume that $q=p_1q_1p_2q_2\ldots
p_kq_k $ for some paths $q_1, \ldots , q_k$ in $\G $.  By
$\mathcal Q$ we denote the set of all irreducible cycles $q^\prime
$ with the same initial point as $q$ in $\G $ such that $q^\prime
$ can be represented as
\begin{equation}
q^\prime =p_1^\prime q_1p_2^\prime q_2\ldots p_k^\prime q_k,
\end{equation}
where $p_1^\prime , p_2^\prime , \ldots , p_n^\prime $ are some
$H_\lambda $--components of $q^\prime $ satisfying the conditions
$(p_i^\prime )_-=(p_i)_-$ and $(p_i^\prime )_+=(p_i)_+$ for
$i=1,2, \ldots , k$. In particular, we have $\overline{\phi
(p^\prime _i)}=\overline{\phi (p_i)}$ for any $i=1, \ldots , k$.
Note that for any $q^\prime \in \mathcal Q$, the subpaths
$p_1^\prime , p_2^\prime , \ldots , p_k^\prime $ are isolated
$H_\lambda $--components of $q^\prime $. By $\mathcal D$ we denote
the set of all diagrams $\Delta $ over (\ref{G1}) such that $\phi
(\partial \Delta )\equiv \phi (q^\prime ) $ for some $q^\prime \in
\mathcal Q$.

Now let $\Delta $ be a diagram of minimal type in $\mathcal D$. To
simplify our notation we will identify the boundary of $\Delta $
with the corresponding cycle $q^\prime =p_1^\prime q_1p_2^\prime
q_2\ldots p_k^\prime q_k$ in $\G $. We are going to show that
every edge of subpaths $p_1^\prime , p_2^\prime , \ldots ,
p_k^\prime $ of $\partial \Delta $ belongs to the boundary of some
$\mathcal R$--cell.

Suppose this is not true; then there are two possibilities: either
a certain $p_i^\prime $ contains an external edge of the second
type, or at least one $\mathcal S$--cell has a common edge with
$p_i^\prime $.

{\it Case 1.} First assume that a certain component $p_i^\prime $,
contains an external edge $d$ of the second type. Then $\partial
\Delta =dc_2d^{-1}c_1$. Since $p_i^\prime $ is isolated, the cycle
$dc_2d^{-1}$ is a subpath of $p_i^\prime $ (see Fig. 2.4). As
$q^\prime $ is irreducible, the subdiagram $\Sigma  $ bounded by
the cycle $dc_2d^{-1}$ contains at least one cell. However, this
contradicts to our choice of $\Delta $ as we can decrease the type
of the diagram by eliminating the subdiagram $\Sigma $.

{\it Case 2.} Now assume that some $\mathcal S$--cell $\Pi $ has a
common edge with $p_i^\prime $ for some $i$. Obviously $\Pi $ is
an $S_\lambda $--cell. Let $\partial \Pi =e^{-1}r$, where $e$ is
the largest common subpath of $p_i^\prime $ and $\partial \Pi $,
and let $\partial \Delta =tp_i^\prime $, where $p_i ^\prime =
s_2es_1 $ (see Fig. 2.4). Then we can decrease the type of the
diagram by passing to the subdiagram bounded by the paths
$ts_2rs_1$, which obviously belongs to $\mathcal Q$. Thus we get a
contradiction again.

Therefore, each edge of any subpath $p_i^\prime $ of $\partial
\Delta $ belongs to the boundary of some $\mathcal R$--cell of
$\Delta $. Since we assumed (\ref{G1}) to be reduced, the length
of any $H_\lambda $--component of any $R\in \mathcal R$ equals $1$
and so the element represented by the label of any edge of any
$p_i^\prime $ belongs to $\Omega _\lambda $. This gives
(\ref{om1}). Finally, by our choice of $\Delta $, we have
$N_\mathcal R (\Delta)\le Area^{rel} (q) $. Hence the total number
of edges of $p_1^\prime , p_2^\prime , \ldots , p_k^\prime $
satisfies the inequality
$$
\sum\limits_{i=1}^kl(p_i^\prime )\le MN_\mathcal R (\Delta)\le M
Area^{rel} (q),
$$
which yields (\ref{om2}). The lemma is proved.
\end{proof}

Note that for any path $q$ in $\G $, every $H_\lambda $--component
$p$ of $q$ is contained in the unique maximal system of connected
$H_\lambda $--components of $q$. In the next section we will use
the following corollary of Lemma \ref{Omega}.

\begin{cor}\label{cut}
Suppose that a group $G$ is given by the reduced finite relative
presentation (\ref{G1}). Let $q=pr$ be a cycle in the Cayley graph
$\G $, where $p$ an $H_\lambda $--component of $q$ for some
$\lambda\in\Lambda $. Let $p, p_1, p_2,\ldots , p_l$ be the
corresponding maximal system of connected $H_\lambda $--components
of $q$. Then
\begin{equation} \label{cut1}
\overline{\phi (p)}\in \langle \overline{\phi (p_1)},
\overline{\phi (p_2)}, \ldots , \overline{\phi (p_l)}, \Omega
_\lambda \rangle .
\end{equation}
\end{cor}

\begin{proof}
Without loss of generality we may assume that $$ r=r_1p_1\ldots
r_lp_lr_{l+1}$$ for some paths $r_1, \ldots , r_l+1$. Let $c_1,
\ldots , c_{l+1}$ be edges in $\G $ labelled by elements of
$\widetilde H_\lambda $ such that
$$(c_1)_-=p_+, \;\;\; (c_{l+1})_+=p_-,$$ and
$$(c_i)_+=(p_i)_-, \;\;\; (c_{i+1})_-=(p_i)_+ $$ for $i=1, \ldots
, l$ (see Fig. \ref{cutfig}).

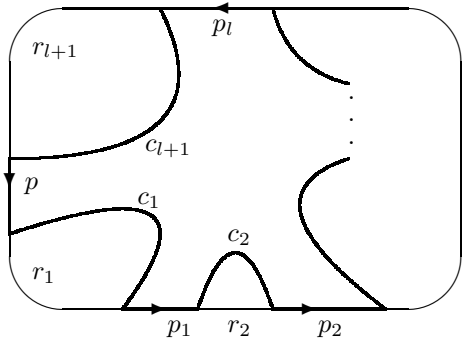
\begin{figure}
\begin{picture}(120,50)(-30,-5)
\put(30,20){\oval(60,40)}

\thicklines{ \put(0,10){\line(0,1){10}} \put(15,0){\line(1,0){10}}
\put(35,0){\line(1,0){15}} \put(20,40){\line(1,0){15}}
\put(0,16){\vector(0,-1){0}} \put(21,0){\vector(1,0){0}}
\put(41,0){\vector(1,0){0}} \put(27,40){\vector(-1,0){0}} }

\put(2,16){$p$} \put(21,-3){$p_1$} \put(41,-3){$p_2$}
\put(27,37){$p_l$}

\put(3,34){$r_{l+1}$} \put(3,4){$r_1$} \put(29,-3){$r_2$}

\put(18,21){$c_{l+1}$} \put(17,14){$c_1$} \put(29,9){$c_2$}

\qbezier(20,40)(30,20)(0,20) \qbezier(15,0)(30,20)(0,10)
\qbezier(25,0)(30,15)(35,0) \qbezier(50,0)(30,15)(45,20)
\qbezier(45,30)(37,32)(35,40)

\put(45,22){.} \put(45,25){.} \put(45,28){.}
\end{picture}

\caption{The decomposition of the cycle $q$ in the proof of
Corollary \ref{cut}.} \label{cutfig}
\end{figure}

Let us consider cycles $q_j=c_jr_j^{-1}$. Since the system $p,
p_1, \ldots , p_l$ is maximal, it follows that $c_j$ is an
isolated $H_i$--component of $q_j$ for every $j$. Applying Lemma
\ref{Omega} for cycles $q_j$, $j=1, \ldots , l+1$, we obtain
$$
\overline{\phi (c_j)} \in \langle \Omega _\lambda \rangle .
$$
This immediately implies (\ref{cut1}) since $p=(c_1p_1\ldots
c_lp_lc_{l+1})^{-1}$.
\end{proof}

Applying Lemma \ref{Omega}, we can prove the following.

\begin{prop}\label{Hfg}
Let $G$ be a group, $\{ H_\lambda \} _{\lambda \in \Lambda }$ a
collection of subgroups of $G$. Suppose that $G$ is generated by a
finite set $X$ in the ordinary non--relative sense and is finitely
presented with respect to $\{ H_\lambda \} _{\lambda \in \Lambda
}$. Then each subgroup $H_\lambda $ is generated by $\Omega
_\lambda $. In particular, $H_\lambda $ is finitely generated.
\end{prop}

\begin{proof}
Suppose that $G$ is given by a reduced finite relative
presentation (\ref{G1}). Let us fix $\lambda \in \Lambda $. For
every nontrivial element $h\in H_\lambda $, there is a word $W$
over $X$ such that $\overline{W}=h$ in $G$. Let us consider a
cycle $q=pr$ in $\G $ such that $\overline{\phi (p)}=h$ and $\phi
(r)\equiv W^{-1}$. Clearly $p$ is an isolated $H_\lambda
$--component of $q$ as $r$ contains no edges labelled by elements
of $\widetilde{H} _\lambda $. Applying Lemma \ref{Omega} for
$k=1$, $p_1=p$, we obtain $h\in \langle \Omega _\lambda \rangle $.
Therefore $\Omega _\lambda $ generates $H_\lambda $. As the
presentation is finite, the set $\Omega _\lambda $ is finite and
the proposition is proved.
\end{proof}


\section{Relative Dehn functions}


In the case of a finitely generated group $G$ and a finitely
generated subgroup $H\le G$, the relative Dehn function
corresponding to the pair $(G, H)$ was introduced in \cite{BC} as
a special case of the Howie function of a complex of groups (see
\cite{BC1} and \cite{Co2} for details). However, the
straightforward extension of the definition from \cite{BC} does
not work in case the number of subgroups is greater that one.

Indeed, the definition from \cite{BC} involves the following
construction. Let $(L,K) $ be a geometric realization of the pair
$(G,H)$, i.e., a pair of connected combinatorial 2--complexes with
finite 1--skeletons such that $K$ is a subcomplex of $L$ and there
are finitely many cells in $L\setminus K$. Moreover, suppose that
for any choice of base point in $K$, there exist isomorphisms $\pi
_1(L)\to G$, $\pi _1(K)\to H$ for which the following diagram is
commutative
$$
\begin{CD}
\pi _1(L) @>>> G\\ @AAA @AAA\\ \pi _1(K) @>>> H
\end{CD}
$$
(The vertical maps are the natural embeddings.) Then one can form
the combinatorial quotient $Q=L/K$ by contracting $K$ to a point
(see \cite{Ge1}). The relative Dehn function of the pair $(G,H)$
is defined as the Dehn function of the universal cover $E$ of $Q$
related to the associated complex of groups (\cite{BC},
\cite{Co2}).

However, if $L=K_1\cup K_2$, we obtain a point by contracting both
$K_1$ and $K_2$. This situation arises when we consider a free
product $G=H_1\ast H_2$ of two finitely presented groups $H_1$ and
$H_2$, and the corresponding two--complexes $L,$ $K_1$, $K_2$
canonically associated with the finite presentations of $G$,
$H_1$, and $H_2$ respectively.

In the present paper we exploit another approach based on the
notion of relative presentation. Our treatment is inspired by the
paper \cite{Pan1} and is similar to that in \cite{BC}. However we
use a language, which is rather combinatorial in contract to
geometric one in \cite{BC}. This allows to define the notion of
the relative Dehn function with respect to an arbitrary set of
arbitrary, not necessarily finitely generated, subgroups. In the
case of a single subgroup, our definition is equivalent to that of
Brick and Corson.

Suppose that $G$ is given by the relative presentation \ref{G1}.
For a word $W\in (X\cup \mathcal H)^\ast $ such that $W$
represents $1$ in $G$, there exists an expression
\begin{equation}
W=_F\prod\limits_{i=1}^k f_i^{-1}R_if_i \label{prod}
\end{equation}
with the equality in the group $F$ given by \ref{F}, where $R_i\in
\mathcal R$ and $f_i\in F$ for any $i$.

\begin{defn}\label{IP}
We say that a function $f:\mathbb N\to \mathbb N$ is a {\it
relative isoperimetric function} of the presentation (\ref{G1})
with respect to subgroups $\{H_\lambda \} _{\lambda \in \Lambda }$
if for any $n\in \mathbb N$ and any word $W\in (X\cup \mathcal
H)^\ast $ of length $\| W\| \le n$ representing the identity in
the group $G$, one can write $W$ as in (\ref{prod}) with $k\le
f(n).$
\end{defn}

Notice that the above definition coincides with the definition of
the ordinary Dehn function in case $\mathcal H=\emptyset  $. The
standard arguments show that Definition \ref{IP} is equivalent to
the following.

\begin{defn}
A function $f:\mathbb N\to \mathbb N$ is a {\it relative
isoperimetric function} of the presentation (\ref{G1}) with
respect to subgroups $\{H_\lambda \} _{\lambda \in \Lambda }$ if
for any $n\in \mathbb N$ and any cycle $q$ in the Cayley graph $\G
$ of length $l( q)\le n$, the relative area of $q$ with respect to
the presentation (\ref{G1}) satisfies $$ Area^{rel} (q)\le f(n).$$
\end{defn}

\begin{defn}
The smallest relative isoperimetric function of (\ref{G1}) is
called the {\it relative Dehn function} of (\ref{G1}) with respect
to subgroups $\{H_\lambda \} _{\lambda \in \Lambda }$. As was
notified in the introduction, it can happen that (\ref{G1}) does
not posses any finite relative isoperimetric function. In this
case we say that the relative Dehn function is not well--defined.
\end{defn}

As in the usual (non--relative) case, we consider the relative
isoperimetric functions up to the following equivalence relation.

\begin{defn}\label{equiv}
For two functions $f,g: \mathbb N\to \mathbb N$, we say that $f$
is asymptotically less than $g$ and write $f\preceq g$ if there
exist constants $C,K,L\in \mathbb N$ such that $$ f(n)\le
Cg(Kn)+Ln.$$ Further, we say that $f$ is asymptotically equivalent
to $g$ and write $f\sim g$ if $f\preceq g$ and $g\preceq f$.
\end{defn}

\begin{thm} \label{Df}
Let
\begin{equation}\label{P1}
\langle X_1, H_\lambda, \lambda\in \Lambda \; | \;  R=1, R\in
\mathcal R_1\rangle
\end{equation}
and
\begin{equation}\label{P2}
\langle X_2, H_\lambda, \lambda\in \Lambda \; | \;  R=1, R\in
\mathcal R_2\rangle
\end{equation}
be two finite relative presentations of the same group $G$ with
respect to a fixed collection of subgroups $\{ H_\lambda \}
_{\lambda \in \Lambda }$, $\delta _1$ and $\delta _2$ the
corresponding relative Dehn functions. Suppose that $\delta _1$ is
well--defined, i.e., $\delta _1(n)$ is finite for every $n$. Then
$\delta _2$ is well--defined and $\delta_1\sim \delta _2$.
\end{thm}

\begin{proof}
Although the presentations are infinite (in the usual,
non--relative sense), passing from (\ref{P1}) to (\ref{P2}) we
change only a finite part of the presentation as $X_i$ and
$\mathcal R_i$ are finite for $i=1,2$. Thus in order to prove our
theorem we have to repeat word--for--word the proof of its
non--relative analogue replacing the word "length" with "relative
length" everywhere. Proposition \ref{Lip} will play the same role
as the fact that the word metric on a finitely generated group is
independent of the choice of a finite generating set up to the
Lipschitz equivalence. This is straightforward and we leave
details to the reader.
\end{proof}

In what follows, speaking about the relative Dehn function of $G$
with respect to $\{ H_\lambda \} _{\lambda \in \Lambda }$ we
always mean the corresponding equivalence class.

\begin{defn}
Let $G$ be a group, $\{ H_\lambda \} _{\lambda \in \Lambda }$ a
collection of subgroups of $G$. We say that $G$ is {\it hyperbolic
relative to} $\{ H_\lambda \} _{\lambda \in \Lambda }$, if $G$ is
finitely presented with respect to $\{ H_\lambda \} _{\lambda \in
\Lambda }$ and the relative Denh function of $G$ with respect to
$\{ H_\lambda \} _{\lambda \in \Lambda })$ is linear.
\end{defn}

\begin{prop}\label{malnorm}
Let $G$ be a group, $\{ H_\lambda \} _{\lambda \in \Lambda }$ a
collection of subgroups of $G$. Suppose that $G$ is finitely
presented with respect to $\{ H_\lambda \} _{\lambda \in \Lambda
}$ and the relative Denh function $\delta ^{rel}$ of $G$ with
respect to $\{ H_\lambda \} _{\lambda \in \Lambda })$ is
well--defined. Then the following conditions hold.

1) The intersection $H_\lambda ^{g_1} \cap H_\mu ^{g_2}$ is finite
whenever $\lambda \ne \mu $.

2) The intersection $H_\lambda ^g \cap H_\lambda $ is finite for
any $g\ne H_\lambda $.
\end{prop}

\begin{proof}
First suppose that $\lambda \ne \mu$. It is sufficient to check
that $H_\lambda ^g\cap H_\mu $ is finite for every $g\in G$ since
$ H_\lambda ^{g_1} \cap H_\mu ^{g_2} = (H_\lambda ^g\cap H_\mu
)^{g_2}$ for $g=g_1g_2^{-1}$.

Consider a word $W\in (X\cup \mathcal H)^\ast $ that represents
$g$ and has length $$\| W\| =|g|_{X\cup\mathcal H}.$$ Assume that
$W=W_1W_2$, where $W_1$ is the maximal (may be empty) prefix of
$W$ consisting of letters from $\widetilde H_\lambda $. Denote by
$f$ the element of $G$ represented by $W_2$. It is clear that
$H^g_\lambda =H^{f}_\lambda $. Thus it suffices to show that
$H^{f}_\lambda \cap H_\mu $ is finite.

Taking into account this remark, we can always assume that if
$g\ne 1$, then the first letter of the shortest word $W\in (X\cup
\mathcal H)^\ast $ representing the element $g$ does not belong to
$H_\lambda $.

Let us take an arbitrary element $h\in H_\lambda ^g \cap H_\mu $
and denote by $h_1$, $h_2$ the letters from $\widetilde H_\lambda
$ and $\widetilde H_\mu $ that represent elements $h^{g^{-1}}\in
H_\lambda $ and $h\in H_\mu $ respectively. Since
$$
\overline{W^{-1}h_1W}=\overline{h_2},
$$
there is a cycle $q$ in $\G $ having the label
$$
\phi (q) \equiv W^{-1}h_1Wh_2^{-1},
$$
of relative area
\begin{equation}\label{mal1}
Area ^{rel} (q)\le \delta ^{rel}(l(q))\le \delta ^{rel}(2\| W\|
+2)=\delta ^{rel}(2|g|_{X\cup\mathcal H} +2).
\end{equation}
Let
\begin{equation}\label{q1}
q=v_1p_1v_2p_2,
\end{equation}
where
\begin{equation}\label{q2}
\phi (v_1)\equiv W^{-1}, \;\;\; \phi (v_2)\equiv W,
\end{equation}
and
\begin{equation}\label{q3}
\phi (p_1)=h_1, \;\;\; \phi (p_2)=h_2^{-1}.
\end{equation}

Note that the subpath $p_1$ is an isolated $H_\lambda $--component
of $q$. Indeed, since the first letter of $W$ does not belong to
$H_\lambda $, $p_1$ is an $H_\lambda $--component. Suppose that
there is another $H_\lambda $--component $p$ of $q$ that is
connected to $p_1$. Since $\lambda \ne \mu$, $p$ is a subpath of
$v_1$ or $v_2$. For definiteness, assume that $$v_1=v_1^\prime
pv_2^{\prime \prime },$$ where $|v_1^{\prime \prime }|\ne 0$, and
there exists a path $c$ in $\G $ that connects some vertex of $p$
to a vertex of $p_1$ and has label consisting of letters from
$\widetilde H_\lambda $ (see Fig. \ref{malfig}). Then $\phi
(v_1^{\prime \prime })$ represents an element of $H_\lambda $ in
$G$. Hence we can take a letter $k\in \widetilde H_\lambda $ such
that
$$\overline{k}=\overline{\phi (p)}\overline{\phi (v_2^{\prime
\prime })} $$ and consider the word $\phi (v_1^\prime )k\in (X\cup
\mathcal H)^\ast $, which is shorter than $W$ and represents the
same element $g$ of $G$. A contradiction.

\begin{figure}
\begin{picture}(120,35)(0,-2)
\multiput(40,0)(0,30){2}{\line(1,0){40}}
\multiput(40,0)(40,0){2}{\line(0,1){30}}

\multiput(40,6)(0,20){2}{\vector(0,1){0}}
\put(80,14){\vector(0,-1){0}}

\multiput(40,0)(0,10){4}{\circle*{1}}
\multiput(80,0)(0,30){2}{\circle*{1}}

\put(36,4){$v_1^\prime $} \put(36,24){$v_2^{\prime \prime }$}
\put(36,14){$p$} \put(82,14){$v_2$} \put(59,2){$p_2$}
\put(59,32){$p_1$}

\thicklines \qbezier(40,20)(60,20)(50,30) \put(51,18){$c$}
\put(40,16){\vector(0,1){0}} \put(40,10){\line(0,1){10}}
\put(59,30){\vector(1,0){0}}
\multiput(40,0)(0,30){2}{\line(1,0){40}}
\put(59,0){\vector(-1,0){0}}
\end{picture}
\caption{The cycle in $\G $ associated with $h\in H^g_\lambda \cap
H_\mu $.} \label{malfig}
\end{figure}

Thus $p_1$ is an isolated $H_\lambda $--component of $q$ and we
can apply Lemma \ref{Omega}. Using (\ref{mal1}), we obtain
\begin{equation}\label{mal3}
|h_1|_{\Omega _\lambda }=|\overline{\phi (p_1)}|_{\Omega _\lambda
}\le M\cdot Area^{rel}(q)\le M \delta ^{rel}(2|g|_{X\cup \mathcal
H} +2).
\end{equation}
If $\delta ^{rel}(n)$ is finite for any $n\in \mathbb N$, this
means that the length of every element $h_1\in H^g_\lambda \cup
H_\mu $ with respect to $\Omega _\lambda $ is bounded by a
constant which is independent of $h_1$. As $\Omega _\lambda $ is
finite, we have $|H^g_\lambda \cap H_\mu |< \infty $.

Let us prove the second assertion of the theorem. Suppose that
$g\notin H_\lambda $ and $W$ is as above. Arguing in the analogous
way, we can assume that the first letter of the word $W\in (X\cup
\mathcal H)^\ast $ does not belong to $H_\lambda $. As in the
previous case, we construct a cycle $q$ in $\G $ satisfying
(\ref{q1}), (\ref{q3}) and note that $p_1$ is an isolated
$H_\lambda $--component of $q$. The only additional argument we
have to use in this case is that $p_1$ can not be connected to the
$H_\lambda $--component of $q$ containing $p_2$. Indeed, for
otherwise we have $g\in H_\lambda $ that contradicts to our
assumption. The rest of the proof coincides with that in the
previous case.
\end{proof}

In case of torsion free groups we have immediately

\begin{cor}
Suppose that $G$ is a torsion free group, $G$ is finitely
presented with respect to a collection of subgroups $\{ H_\lambda
\} _{\lambda \in \Lambda }$, and the relative Denh function of $G$
with respect to $\{ H_\lambda \} _{\lambda \in \Lambda })$ is
well--defined. Then each subgroup $H_\lambda $ is malnormal, i.e.,
$H_\lambda ^g\cap H_\lambda =\{ 1\} $ whenever $g\notin H_\lambda
$.
\end{cor}

Let us mention one more corollary of Proposition \ref{malnorm}.

\begin{cor}\label{gninH}
Suppose that a group $G$ is finitely presented with respect to a
collection of subgroups $\{ H_\lambda \} _{\lambda \in \Lambda }$
and the relative Denh function of $G$ with respect to $\{
H_\lambda \} _{\lambda \in \Lambda })$ is well--defined. Suppose
that $g\in G$ is an element of infinite order and $g^n\in
H_\lambda $ for some $n\in \mathbb Z\setminus \{ 0\} $ and some
$\lambda \in \Lambda $. Then $g\in H_\lambda $.
\end{cor}

\begin{proof}
If $g^n\in H_\lambda $, then the intersection $H^g_\lambda \cup
H_\lambda $ contains $\langle g^n\rangle $. Since $\langle
g^n\rangle $ is infinite, $g\in H_\lambda $ by Proposition
\ref{malnorm}.
\end{proof}

To formulate our next results, we need an auxiliary notion.

\begin{defn}
A function ${f : {\mathbb N}\to {\mathbb N}}$ is said to be {\it
subnegative} if $$f(a+b)\ge f(a)+f(b)$$ for any $a,b\in {\mathbb
N}$. Given an arbitrary function  ${f : {\mathbb N}\to {\mathbb
N}}$, the {\it subnegative closure} of $f$ is defined to be
\begin{equation}
\bar f(n)=\max_{i=1,\dots , n}\left(\max_{a_1+\dots + a_i=n,\
a_i\in {\mathbb N}} \left(f(a_1)+\dots +f(a_i)\right)\right)
\end{equation}
In fact, $\bar f$ is the smallest subnegative function such that
$\bar f(n)\ge f(n)$ for all $n$.
\end{defn}

Below we will speak about subnegative closure of (relative) Dehn
functions. It is easy to see that even the ordinary Dehn function
of a finitely presented group is not necessarily subnegative. For
example, for the group presentation $\langle a\; | \; a=1, a^2=1
\rangle $ we have $\delta (1)=1$, $\delta (2)=1$. The question
whether or not every Dehn function of a finitely presented group
is equivalent to some subnegative function is more complicated. We
only note that this question is still open and refer the reader to
\cite{GubS} for more comprehensive discussion.

\begin{thm} \label{oDf}
Let $G$ be a group, $\{H_1, \ldots H_m\} $, $\{ K_1, \ldots ,
K_l\} $ two finite collections of subgroups of $G$. Assume that
$G$ is finitely presented with respect to $\{H_1, \ldots H_m\}
\cup \{ K_1, \ldots , K_l\} $ and each subgroup $H_i$ is finitely
presented itself. Then $G$ is finitely presented with respect to
$\{ K_1, \ldots , K_l\} $.

Moreover, if $\delta _1, \ldots , \delta _m$ are the ordinary Dehn
functions of $H_1, \ldots , H_m$ respectively and the relative
Dehn function $\delta $ of $G$ with respect to $\{H_1, \ldots
H_m\} \cup \{ K_1, \ldots , K_l\} $ is finite for each value of
the argument, then the relative Dehn function $\gamma $ of $G$
with respect to $\{ K_1, \ldots , K_l\} $ is well--defined and
satisfies the inequality
\begin{equation}\label{D0}
\gamma (n)\preceq \bar f\circ \delta (n),
\end{equation}
where $\bar f$ is the subnegative closure of the function $$
f(n)=\max\limits_{i=1, \ldots , m} \delta _i(n).$$
\end{thm}

\begin{proof}
Before proving the theorem we have to make a few remarks. In the
particular case when $G$ is a fundamental group of a complex of
groups with finite edge groups, $H_1, \ldots , H_m$ are vertex
groups, and $l=0$ this theorem is equivalent to the main result of
the paper \cite{BC1}. In the case $l=0$ and the group $G$ is
hyperbolic with respect to $\{H_1, \ldots H_m\} $ in the sense of
Farb with the BCP property the analogue of Theorem \ref{oDf} can
be found in \cite{F}. (However the statement of the theorem in
\cite{F} does non involve the subnegative closure, which is
required in order to make the proof correct.) Although the proof
in the general case exploit similar ideas, we provide it here for
convenience of the reader.

By the assumptions of the theorem, the group $G$ is a quotient of
the group $$F=F(X)\ast \widetilde H_1\ast \ldots \ast \widetilde
H_m\ast \widetilde K_1\ast \ldots \ast \widetilde K_l,$$ where
$\widetilde H_i\cong H_i$ and $\widetilde K_j\cong K_j $. Let
$$
\mathcal H=\bigsqcup\limits_{i=1}^m \left( \widetilde H_i
\setminus\{ 1\}\right)
$$
and
$$
\mathcal K=\bigsqcup\limits_{j=1}^l \left( \widetilde K_j
\setminus\{ 1\}\right) .
$$

We start with a finite reduced relative presentation
\begin{equation}\label{D1}
G=\left\langle X\cup \mathcal H\cup\mathcal K\; \left|  \; S=1,\;
S\in \bigcup\limits_{i=1}^m \mathcal S_i, \; P=1,\; P\in
\bigcup\limits_{i=1}^m \mathcal P_i, \; R=1,\; R\in \mathcal R
\right.\right\rangle ,
\end{equation}
of the group $G$ with respect to $\{ H_1 , \ldots H_m\} \cup \{
K_1, \ldots , K_l\} $, where $\mathcal S_i$ (respectively
$\mathcal P_i$) is the set of all words over $\widetilde H_i$
(respectively $\widetilde K_i$) representing the identity in $G$.
Suppose that the groups $H_1 , \ldots H_m$ have finite
presentations
\begin{equation}\label{D2}
H_i=\langle Y_i\; |\; T=1,\; T\in \mathcal T_i \rangle , \;\; i=1,
\dots , m.
\end{equation}
Let $O_i$ denote the set of all letters from $\widetilde H_i$ that
appear in words $R\in \mathcal R$. Since $O_i$ is finite for every
$i$, without loss of generality we can assume that the (finite,
symmetrized) generating set $Y_i$ of $H_i$ contains all elements
of $H_i $ that are represented by letters from $O_i$. Thus we can
regard words from $\mathcal R$ as words in the alphabet
$$Z=X\bigcup \left(\bigcup\limits_{i=1}^m Y_i\right)\bigcup
\mathcal K.$$ We set
$$ \mathcal T =\bigcup\limits_{i=1}^m \mathcal T_i .$$

Evidently $G$ can be defined by the finite relative presentation
\begin{equation}\label{D3}
G=\left\langle X\cup \left(\bigcup\limits_{i=1}^mY_i\right) \cup
\mathcal K\; \left| \; P=1,\; P\in \bigcup\limits_{i=1}^l \mathcal
P_i,\; T=1,\; T\in \mathcal T, \; R=1,\; R\in \mathcal R
\right.\right\rangle
\end{equation}
with respect to $\{ K_1, \ldots , K_l\} $. Assume that $\delta _i,
i=1, \ldots m$, are the Dehn function of the presentations
(\ref{D2}), and $\delta $, $\gamma $ are the relative Dehn
function of the presentations (\ref{D1}) and (\ref{D3})
respectively.

Let us fix an arbitrary $n\in \mathbb N$. Consider a word $W$ over
the alphabet $Z$ such that
\begin{equation}\label{D4}
\| W\|\le n,
\end{equation}
and $W$ represents the identity in $G$. Since $Z\subseteq X\cup
\mathcal H\cup \mathcal K$, one can regard $W$  as the word over
$X\cup \mathcal H\cup \mathcal K$. We take a van Kampen diagram
$\Delta $ over (\ref{D1}) such that
\begin{equation}\label{D5}
\phi (\partial \Delta )\equiv W
\end{equation}
and assume that $\Delta $ has minimal type among all van Kampen
diagrams satisfying (\ref{D5}). Thus
\begin{equation}\label{D6}
N_\mathcal R(\Delta )\le \delta (n).
\end{equation}

Suppose that $\Pi $ is an $S_i$--cell in $\Delta $. If $e\in
\partial \Pi $ is an internal edge in $\Delta $, then $e$ belongs
to the boundary of some $\mathcal R$--cell by Lemma
\ref{DiagMinType} and therefore $\phi (e)\in Y_i$ according to our
choice of $Y_i$. If $e\in \partial \Pi $ is external, then $\phi
(e)\in Y_i$ since $W$ is a word over $Z$. Thus for any $S_i$--cell
$\Pi $ of $\Delta $, $\phi (\partial \Pi )$ is a word over $Y_i$.
This observation allows to transform $\Delta $ into a van Kampen
diagram $\Theta $ over (\ref{D3}) in the following way. For every
$S_i$--cell $\Pi $ of $\Delta $, we consider a van Kampen diagram
$\Sigma (\Pi )$ over (\ref{D2}) whose boundary label coincides
with the boundary label of $\Pi $ and replace $\Pi $ with $\Sigma
(\Pi )$ in $\Delta $. Doing this for all $S_i$-cells $\Pi $ in
$\Delta $, $i=1, \ldots , m$, we get a diagram $\Theta $ over
(\ref{D3}).

We want to estimate the number of cells in $\Theta $ assuming that
for every $\mathcal S$--cell $\Pi $ in $\Delta $, $\Sigma (\Pi )$
has minimal possible number of cells. Let $S(\Delta )$ denote the
set of all $\mathcal S$--cells of the diagram $\Delta $. We also
recall that $Area\; (\Sigma (\Pi) )$ denotes the number of all
cells in the diagram $\Sigma (\Pi)$.

Using Corollary \ref{sums}, we obtain
\begin{equation}\label{D7}
\sum\limits_{\Pi \in S(\Delta )} l(\partial \Pi ) \le M N_\mathcal
R(\Delta )+l(\partial \Delta ),
\end{equation}
where $M$ is the maximum of lengths of the relators from $\mathcal
R$. Combining (\ref{D4})--(\ref{D7}), we have $$
\begin{array}{ll}
\sum\limits_{\Pi \in S(\Delta )} Area\; (\Sigma (\Pi )) & \le
\sum\limits_{\Pi \in S(\Delta )} f(l(\partial \Pi ))\le \bar f
\left( \sum\limits_{\Pi \in S(\Delta )} l(\partial \Pi )\right)
\\ & \\ & \le \bar f (MN_\mathcal R(\Delta )+\| W\| ) \le \bar f (M\delta
(n)+n).
\end{array}
$$
Finally, we have the following estimate for the sum of the numbers
$N_\mathcal R (\Theta )$ and $N_\mathcal T(\Theta )$ of $\mathcal
R$ and $\mathcal T$--cells in $\Theta $ respectively:
$$
\begin{array}{ll}
N_\mathcal R (\Theta )+ N_\mathcal T(\Theta )& \le N_\mathcal
R(\Delta )+\sum\limits_{\Pi \in S(\Delta )} Area\; (\Sigma (\Pi ))
\\ &\\ & \le \delta (n)+ \bar f (M\delta (n)+n).
\end{array}
$$ This yields (\ref{D0}).
\end{proof}

\begin{cor}
Suppose that $G$ is a group hyperbolic relative to a finite
collection of finitely presented subgroups $\{ H_1, \ldots H_m\}$
and $f(n)$ is an isoperimetric function of $H_i$ for any $i=1,
\ldots , m$. Then $G$ is finitely presented itself and $\bar f$ is
an isoperimetric function of $G$. In particular, if each of the
subgroups $H_i$ is hyperbolic, then $G$ is hyperbolic.
\end{cor}

In the particular case of hyperbolic products of groups, this
corollary was obtained by Pankrat'ev \cite{Pan1}.


\section{Splitting Theorem for relatively finitely presented groups}


Although the definition of relative Dehn functions has been given
in the general situation, by technical reasons it is more
convenient to deal with finitely generated groups. In this section
we prove the Splitting Theorem for groups given by finite relative
presentations, which allows to reduce some questions concerning
algebraic properties to the case when the group is finitely
generated in the ordinary non--relative sense.

We begin with some basic notions of the Bass--Serre theory of
groups acting on trees.

\begin{defn}
A {\it finite graph of groups} $\mathcal G$ consists of the
following data.

1) A finite connected oriented graph $\mathcal G$; we denote by
$E(\mathcal G)$ and $V(\mathcal G)$ its set of edges and set of
vertices respectively.

2) For every vertex $v\in V(\mathcal G)$, one associates a group
$G_v$; the groups $G_v$, $v\in V(\mathcal G)$, are called {\it
vertex groups}.

3) For every edge $e\in E(\mathcal G)$, one associates a group
$G_e$ together with monomorphisms $\alpha _e:G_e\to G_{e_-}$,
$\omega _e: G_e\to G_{e_+}$, where $e_-$ and $e_+$ are the origin
and the terminus of the edge $e$ respectively. The groups $G_e$,
$e\in E(\mathcal G)$, are called {\it edge groups.}
\end{defn}

\begin{defn}
Let $\Theta $ be a maximal tree in $\mathcal G$. The fundamental
group $\pi _1(\mathcal G, \Theta )$ of the finite graph of groups
$\mathcal G$ at $\Theta $ is the group generated by the groups
$G_v$, $v\in V(\mathcal G)$, and elements $t_e$, $e\in E(\mathcal
G)$, subject to the relations $$ t_e^{-1}\alpha _e(g)t_e=\omega
_e(g), \;\;\; g\in G_e,\; e\in E(\mathcal G),$$ $$
t_{e^{-1}}=t_{e}^{-1}, \; e\in E(\mathcal G),$$ and
$$ t_e=1,\;\;\; {\rm if }\; e\in \Theta .$$
The group $\pi _1(\mathcal G, \Theta )$ is independent up to
isomorphism of the choice of the maximal tree $\Theta $
\cite[Prop. 20]{Trees}.

\end{defn}

In particular, if $\mathcal G$ is a tree itself, then the
fundamental group has the presentation
\begin{equation} \label{fgtg}
\left\langle G_v, \; v\in V(\mathcal G)\; |\; \alpha _e(g)=\omega
_e(g), \; g\in G_e,\; e\in E(\mathcal G)\right\rangle .
\end{equation}

Now we formulate the main result of this section.

\begin{figure}
\begin{picture}(120,50)
\put(60,10){\circle*{1}} \put(45,40){\circle*{1}}
\put(60,10){\line(-1,2){15}} \put(60,10){\vector(-1,2){9}}
\put(60,10){\line(-2,1){30}} \put(60,10){\line(2,1){30}}
\put(60,10){\vector(2,1){18}} \put(60,10){\vector(-2,1){18}}
\put(30,25){\circle*{1}} \put(90,25){\circle*{1}}

\put(59,17){.}\put(62,17){.}\put(65,15){.} \put(59,6){$Q$}

\put(23,27){$H_{\lambda _1}$} \put(41,42){$H_{\lambda _2}$}
\put(86,27){$H_{\lambda _m}$}

\put(36,16){$L_1$} \put(52,29){$L_2$} \put(78,15 ){$L_m$}
\end{picture}
\caption{The tree of groups $\mathcal G$.} \label{star}
\end{figure}
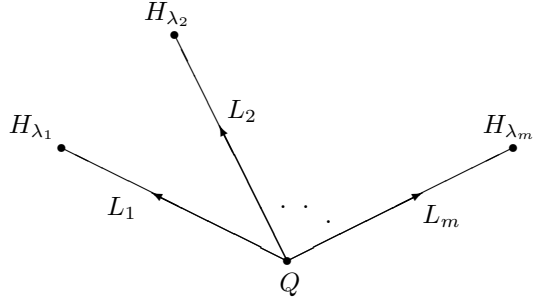

\begin{thm}\label{st}
Let $G$ be a group, $\{ H_\lambda \} _{\lambda \in \Lambda }$ a
collection of subgroups of $G$. Suppose that $G$ is finitely
presented with respect to $\{ H_\lambda \} _{\lambda \in \Lambda
}$.

\noindent 1) There exist a finite subset $\Lambda _0=\{ \lambda
_1, \ldots , \lambda _m\} \subseteq \Lambda $ such that $G$ splits
as the free product
\begin{equation} \label{splG}
G=\left(\ast _{\lambda \in \Lambda\setminus\Lambda _0} H_\lambda
\right) \ast G_0,
\end{equation}
where  $G_0$ is the subgroup of $G$ generated by $H_{\lambda _1},
\ldots , H_{\lambda _m}$ and $X$.

\noindent 2) The groups $G_0$ is isomorphic to the fundamental
group of the tree of groups $\mathcal G$ drawn on Fig. \ref{star}
for some finitely generated groups $L_1, \ldots , L_m$, $Q$, and
the following conditions hold.

a) The group $Q$ is finitely presented with respect to the
collection of subgroups $\{ L_1, \ldots , L_m\} $.

b) If, in addition, the Denh function $\delta ^{rel} _G$ of $G$
with respect to $\{ H_\lambda \} _{\lambda \in \Lambda }$ is
well--defined, then the relative Dehn function $\delta ^{rel} _Q$
of $Q$ with respect to $ \{ L_1, \ldots , L_m\} $ is well--defined
and
\begin{equation} \label{dQ}
\delta ^{rel}_Q \preceq \delta ^{rel}_G\preceq \bar \delta
_Q^{rel},
\end{equation}
where $\bar \delta _Q^{rel}$ is the subnegative closure of $\delta
_Q^{rel}$
\end{thm}

\begin{proof}
Suppose that $G$ is given by some finite reduced relative
presentation (\ref{G1}) with respect to $\{ H_\lambda \} _{\lambda
\in \Lambda }$. We keep our notation from Section 2.2, namely we
will use the sets $\Omega _\lambda $ and $\Omega $. For
simplicity, we will identify the set $\Omega _\lambda $, which is
a subset of elements of $G$, with the subset of $\widetilde
H_\lambda $ consisting of letters representing elements of $\Omega
_\lambda $.

As $G$ is finitely presented relative to $\{ H_\lambda \}
_{\lambda \in \Lambda }$, there is a finite set of indices
$\Lambda _0= \{ \lambda _1 , \ldots , \lambda _m\} $ such that no
relators from $\mathcal R$ involve letters from $\widetilde
H_\lambda $ for $\lambda \in \Lambda \setminus \Lambda _0$.
Evidently we have (\ref{splG}), where
\begin{equation}\label{sT0}
G_0=\langle X, \; \mathcal H\; | S=1, S\in \bigcup\limits_{i=1}^m
\mathcal S_{\lambda _i},\; R=1, \; R\in \mathcal R\rangle .
\end{equation}

We divide the rest of the proof into a few lemmas. For every $i=1,
\ldots , m$, we consider the group $L_i=\langle \Omega_i\rangle $.
Let $$F_Q=F(X)\ast \widetilde L_1\ast \ldots \ast \widetilde
L_m,$$ where $\widetilde L_i$ are isomorphic copies od $L_i$. Set
$$\mathcal L=\bigsqcup\limits_{i=1}^m \left( \widetilde L_i\setminus \{ 1\} \right) .$$
By $\mathcal T_i$ we denote the set of all words over the alphabet
$\widetilde L_i\setminus \{ 1\}$ representing the identity in the
group $G$. Let
$$ \mathcal T=\bigcup\limits_{i=1}^m \mathcal T_{\lambda _i} .$$

\begin{lem}\label{Qlem}
Let $Q$ be the subgroup of $G$ generated by the set $\Omega \cup
X$. Then $Q$ has the finite relative presentation
\begin{equation}
\langle X, \mathcal L \; |\; T=1,\; T\in \mathcal T,\; R=1,\; R\in
\mathcal R\rangle . \label{Q}
\end{equation}
with respect to $\{ L_1, \ldots , L_m\} $. Moreover, the relative
Dehn function of (\ref{Q}) satisfies the inequality $$ \delta
^{rel} _Q\preceq \delta ^{rel}_G.$$
\end{lem}

\begin{proof}
Indeed, let $W$ be a word over $X\cup \mathcal L$ of length $n$
that represents $1$ in the group $G$. To prove the lemma we have
to show that there exists a van Kampen diagram over (\ref{Q}) with
boundary label $W$ and number of cells at most
$\delta^{rel}_G(n)$.

To this end we note that any word over $X\cup \mathcal L$ can be
regarded as a word over $X\cup\mathcal H$. We consider a diagram
$\Theta $ of minimal type over (\ref{G1}) such that $\phi
(\partial \Delta )\equiv W$. Suppose that $e$ is an edge of
$\Theta $ labelled by a letter from $\mathcal H$. If $e$ is
internal, $e$ belongs to the boundary of some $\mathcal R$--cell
by Lemma \ref{DiagMinType}. Hence $\phi (e)$ represents an element
$w$ of $\Omega _{\lambda _i}\subseteq L_i$ for some $i$ and we can
regard $\phi (e)$ as an element of $\widetilde L_{i}\setminus \{
1\} $. Thus if $\Pi $ is an $S$--cell in $\Theta $, then $\partial
\Pi$ is labelled by a word in the alphabet $\widetilde
L_{i}\setminus \{ 1\} $ for a certain $i\in \{ 1, \ldots , m\} $.
As $\mathcal T_i$ contains all words over $\widetilde L_i\setminus
\{ 1\} $ representing the identity in $G$, one can think of $\Pi $
as a cell corresponding to a relator from $\mathcal T_i$. Thus
$\Theta $ can be regarded as a diagram over (\ref{Q}). This
completes the proof.
\end{proof}

Let $\alpha _i: L_i \to Q$ and $\omega _i: L_i\to H_{\lambda _i}$
be the natural embeddings. Then the fundamental group of the graph
of groups $\mathcal G$ can be represented as
\begin{equation}\label{G0fund}
\left\langle Q, H_{\lambda _i}, i=1, \ldots , m\; |\; \alpha_i
(l)=\omega _i(l), l\in L_i, i=1, \ldots , m\right\rangle .
\end{equation}
It is easy to see that using (\ref{Q}), we can obtain (\ref{sT0})
from (\ref{G0fund}) by applying Tietze transformations. This
proves that $G_0$ splits at the fundamental group of the graph of
groups $\mathcal G$.

\begin{lem}
We have $\delta _G^{rel}\preceq \bar \delta _Q^{rel}.$
\end{lem}

\begin{proof}
Let $q$ be a cycle in $\G $ of length $l(q)\le n$. We repeat here
the trick used in the proof of Lemma \ref{Omega}. Let
$q=q_1r_1\ldots q_nr_n$ be the decomposition of $q$ into the
product of components $q_1, \ldots , q_n$ and subpath $r_1, \ldots
, r_n$ labelled by words in $X$ ($r_i$ may be trivial for some
$i$). By $\mathcal Q$ we denote the set of all irreducible cycles
$q^\prime $ in $\G $ with the same initial point as $q$ and such
that $q^\prime $ can be represented as
\begin{equation}
q^\prime =q_1^\prime r_1 \ldots q_n^\prime r_n,
\end{equation}
where  $q_i^\prime $ is an $H_\mu $--component of $q^\prime $ for
the same $\mu $ as $q_i$ and $(q_i^\prime )_-=(q_i)_-$,
$(q_i^\prime )_+=(q_i)_+$ for $i=1,2, \ldots , n$. Let $\mathcal
D$ be the set of all diagrams over (\ref{G1}) with boundary label
$\phi (q^\prime ) $ for some $q^\prime \in \mathcal Q$.

Take the diagram $\Delta $ of minimal type in $\mathcal D$. Then
$\Delta $ can be represented as the union of subdiagrams $\Delta
_1, \Delta _2, \ldots , \Delta _k$, each $\Delta _i$ is
homeomorphic to a disk, connected by some paths entirely
consisting of external edges of the second type (see Fig.
\ref{Ddec}).

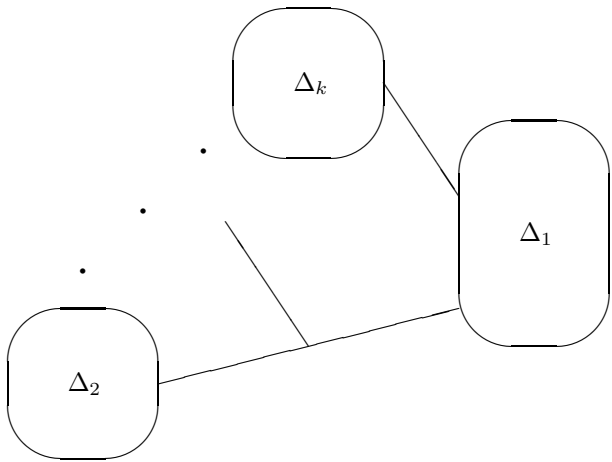
\begin{figure}

\begin{picture}(120,60)(-10,0)
\put(30,10){\oval(20,20)} \put(60,50){\oval(20,20)}
\put(90,30){\oval(20,30)} \put(40,10){\line(4,1){40}}
\put(70,50){\line(2,-3){10}} \put(60,15){\line(-2,3){11}}
\multiput(30,25)(8,8){3}{\circle*{0.8}}

\put(28,9){$\Delta _2$} \put(58,49){$\Delta _k$}
\put(88,29){$\Delta _1$}
\end{picture}

\caption{The decomposition of $\Delta $.} \label{Ddec}
\end{figure}

Repeating the arguments from the proof of Lemma \ref{Omega} and
using our assumptions about $\Delta $, one can easily show that no
$\mathcal S$--cell of $\Delta $ has a common edge with $\partial
\Delta _i$ for some $i$. Therefore, for any $i$, every edge $e$ of
$\Delta _i$ belongs to the boundary of some $\mathcal R$--cell.
Hence $\phi (e)$ can be regarded as a letter from $X\cup \mathcal
L $. As $\mathcal T_i$ contains all words over $\widetilde
L_i\setminus \{ 1\} $ representing the identity in $G$, this
allows us to consider $\Delta _i$, $i=1, \ldots , k$, as diagrams
over (\ref{Q}). We obtain
$$
Area^{rel}(q)=\sum\limits_{i=1}^k N_\mathcal R(\Delta _i)\le
\sum\limits_{i=1}^k \delta^{rel}_Q(l(\partial \Delta _i)) \le
\bar\delta ^{rel}_Q \left( \sum\limits_{i=1}^k l(\partial \Delta
_i) \right) \le \bar\delta ^{rel}_Q (n).
$$
\end{proof}

\end{proof}

We notice one corollary of Theorem \ref{st}.

\begin{cor} If the group $G$ in the Theorem \ref{st} is hyperbolic
with respect to the collection $\{ H_\lambda \} _{\lambda \in
\Lambda }$, then $Q$ is hyperbolic with respect to the collection
$\{ L_1, \ldots , L_m\} $.
\end{cor}

The Splitting Theorem together with the Gru\v sko--Neumann theorem
\cite[Proposition 3.7]{LS} also imply

\begin{cor}\label{Lambda}
Let $G$ be a finitely generated group, $\{ H_\lambda \} _{\lambda
\in \Lambda }$ a collection of nontrivial subgroups of $G$.
Suppose that $G$ is finitely presented with respect to $\{
H_\lambda \} _{\lambda \in \Lambda }$. Then the collection $\{
H_\lambda \} _{\lambda \in \Lambda }$ is finite, i.e., $card\;
\Lambda < \infty $.
\end{cor}

\begin{proof}
Given a finitely generated group $K$, by $rank\, K$ we denote the
minimal number of generators of $K$. By the Gru\v sko--Neumann
theorem, we have
\begin{equation}
rank\, G_0 +\sum\limits_{\lambda \in \Lambda\setminus \Lambda _0}
rank\, H_\lambda \le rank\, G < \infty \label{rank1} .
\end{equation}
As $\Lambda _0$ is finite, the corollary follows immediately from
(\ref{rank1}).
\end{proof}

We conclude this section with one more result describing the
geometry of the embedding $Q\to G$. This result will be used in
Section 4.3 (as well as the Splitting Theorem) to study cyclic
subgroups of relatively hyperbolic group. Besides it seems to be
of independent interest.

\begin{prop}\label{QisomG}
Let $G$ be a group, $\{ H_\lambda \} _{\lambda \in \Lambda }$ a
collection of subgroups of $G$. Suppose that $G$ is finitely
presented with respect to $\{ H_\lambda \} _{\lambda \in \Lambda
}$. Let $Q$ and $L_i$, $i=1, \ldots , m$ be the groups provided by
Theorem \ref{st}. We denote by $\dxl $ the relative metric on $Q$
with respect to the collection $\{ L_1, \ldots , L_m\} $. Then the
mapping of metric spaces $(Q, \dxl )\to (G, \dxh )$ induced by the
natural embedding $Q\to G$ is an isometry.
\end{prop}

\begin{proof}
It is easy to see that $|g|_{X\cup \mathcal L}\ge |g|_{X\cup
\mathcal H}  $ for any $g\in G$. Let us prove the converse
inequality.

Suppose that $g\in Q$ is an element such that the relative length
of $g$ in $Q$ satisfies $|g|_{X\cup \mathcal L}> |g|_{X\cup
\mathcal H} $. Let $W_1$ and $W_2$ be shortest words over $X\cup
\mathcal L$ and $X\cup \mathcal H$ respectively representing the
element $g$. In what follows we regard $\mathcal L$ as a subset of
$\mathcal H$.

Consider the cycle $pq^{-1}$ in $\G $ such that $$\phi (p)\equiv
W_1 , \;\;\;\;\; \phi (q)\equiv W_2.$$ Let $s$ be an $H_\lambda
$--component of $q^{-1}$ for some $\lambda $. We are going to show
that $\lambda \in \Lambda _0$ and $\overline{\phi (s)}\in L_i$ for
some $i\in \{ 1, \ldots , m\} $.

If $s$ is an isolated $H_\lambda $--component of $pq^{-1}$,
$\overline{\phi (s)}\in \langle \Omega_\lambda\rangle $ by Lemma
\ref{Omega}. Note that $\Omega _\lambda $ is non--empty only if
$\lambda \in \Lambda _0$. Thus $\overline{\phi (s)}\in L_i$ for
some $i\in \{ 1, \ldots , m\} $. Further suppose that $s$ is not
isolated in $pq^{-1}$. Let $t_1, \ldots , t_k, s$ be the maximal
connected system of $H_\lambda $--components of $pq^{-1}$
containing $s$. Since the word $W_2$ is a shortest words over
$X\cup \mathcal H$ representing the element $g$, the path $q$ is
geodesic in $\G $. Hence any component of $q$ is isolated in $q$.
It follows that $t_1, \ldots , t_k$ are subpaths of $p$.
Therefore, $\phi (t_1), \ldots , \phi (t_k)$ are $H_\lambda
$--syllables in $W_1$. Thus $\lambda =\lambda _i$ for some $i$ and
$\overline{\phi (t_j)}\in L_i$ for $j=1, \ldots , k$. Finally,
using Corollary \ref{cut}, we obtain
$$
\overline{\phi (s)}\in \langle \Omega _{\lambda _i},
\overline{\phi (t_1)}, \ldots , \overline{\phi (t_k)} \rangle
=L_i.
$$

We have proved that any $H_\lambda $--syllable of $W_2$ represents
an element of $L_i$ for some $i$. This means that one can think of
$W_2$ an a word over $X\cup \mathcal L$. However, $$\| W_2 \|
=|g|_{X\cup \mathcal L}< |g|_{X\cup \mathcal H}= \| W_1\| $$ that
contradicts to the minimality of $W_1$.
\end{proof}


\section{Isoperimetric functions of Cayley graphs.}


In order to apply various well--known results about hyperbolic
spaces to  a pair $(G, \{ H_\lambda \} _{\lambda \in \Lambda })$
with linear relative Dehn function, we have to study more
extensively the geometry of the corresponding Cayley graph  $\G $.
In the present section we establish the relation between the
linearity of the Dehn function of $G$ with respect to $\{
H_\lambda \} _{\lambda \in \Lambda } $ and hyperbolicity of $\G $.

We begin with a definition of the notion of area of a cycle in a
metric graph, which is a particular case of a more general concept
(see \cite[Sec. 5.F]{MG1}, \cite[Ch. III.H, Sec.2]{BriH} or
\cite{Ale57}, \cite{BerN93} for a refined version of area).

First we have to introduce an auxiliary terminology. Let $D^n$
denote the unit disk in the $n$--dimensional Euclidean space.
Recall that a {\it combinatorial map} between $CW$ complexes is a
map that sends open cells homeomorphically onto open cells, and a
{\it combinatorial complex} is a $CW$ complex $K$ such that for
any $(n+1)$--cell $e$ of $K$, the attaching map $\phi _e:
\partial D^{n+1}\to Sk^{(n)}K$ , where $Sk^{(n)}K$ denotes the
$n$--skeleton of $K$, is combinatorial with respect to some
combinatorial structure on $\partial D^{n+1}={\mathbb S}^{n}$.
(This definition involves the induction on dimension.) In what
follows we will also work with a lager category of maps. One
defines a {\it singular combinatorial map} $f: L\to K$ between
$CW$ complexes to be a continuous map such that for every open
$n$--cell $e$ in $L$, either $f| _e$ is a homeomorphism onto an
open cell of $K$ or else $f(e)$ is contained in the
$(n-1)$--skeleton of $L$ \cite{Bri}.

\begin{defn} ($k$--filling, $A^{(k)}$, and Isoperimetric
Inequality). Let $k\ge 3$. A {\it $k$--partition} of $D^2$ is a
homeomorphism $P$ from $D^2$ to a combinatorial $2$--complex in
which every $2$--cell is an $l$--gone for some $3\le l\le k$. We
endow $D^2$ with the induced cell structure and refer to the
preimages under $P$ of $0$--cells, $1$--cells, and $2$--cells as,
respectively, vertices, edges, and faces of $D^2$. We also denote
by $Sk^{(1)}D^2$ the $1$--skeleton of $D^2$ with the induced cell
structure.

Let $\Sigma $ be a graph equipped with the combinatorial metric,
that is, each edge of $\Sigma $ has length $1$. Let $c$ be a
combinatorial cycle in $\Sigma $. A {\it $k$--filling} of $c$
consists of a $k$--partition  $P$ of $D^2$ and a singular
combinatorial map $\Phi :Sk^{(1)}D^2\to \Sigma $ such that $\Phi
(\partial D^2)=c$. In this situation we write $|\Phi |$ to denote
the number of faces of $D^2$ with respect to the cell structure
induced by $P$. The {\it $k$--area } of $c$ is defined to be
$$
A^{(k)} (c)=\min \{ |\Phi |, \; \Phi \; {\rm is\; a \; } k{\rm
-filling\; of \; c} \}
$$
(if there is no $k$--filling of $c$, we put $A^{(k)}(c)=\infty $).

We also set $$ f^{(k)}_\Sigma (n)=\sup\limits_{l(c)\le n}
A^{(k)}(c),
$$ where the supremum is taken over all cycles of combinatorial
length at most $n$.
\end{defn}

It is easy to check that if for some $k\ge 3$ any cycle in $\Sigma
$ admits a $k$--filling and the corresponding function $f_\Sigma
^{(k)}$ is finite for each value of the argument, then for any
$k^\prime \ge k$, $f_\Sigma ^{(k^\prime )}$ is finite for each
value of the argument and, moreover, $f_\Sigma ^{(k^\prime )}\sim
f_\Sigma ^{(k)}$. (The proof of this fact in a much more general
situation can be found in \cite[Ch. III.H, Sec.2]{BriH}.) In this
case $f_\Sigma ^{(k)}$ (regarded up to equivalence) is called the
{\it Dehn function } of $\Sigma $ and is denoted by $f_\Sigma $.
As above, we consider Dehn functions of metric graphs up to the
usual equivalence.

\begin{ex}
In case $\Sigma $ is the Cayley graph of a group $Q$ given by a
finite presentation $\langle Y\; |\; \mathcal Q\rangle $, any
cycle in $\Sigma $ admits a $k$--filling for $k=\max\limits_{Q\in
\mathcal Q} \| Q\| $ and $f^{(k)} (n)< \infty $ for any $n$.
Moreover, $f^{(k)}$ is equivalent to the ordinary Dehn function of
$G$. The assumption of finite presentability of $G$ is essential,
since if $G$ is not finitely presented, then for any $k\in \mathbb
N $, there exists a cycle $c$ is $\Sigma $ that admits no
$k$--filling.
\end{ex}

If we consider a group $G$ which is finitely presented with
respect to a collection of the subgroups $\{ H_\lambda \}
_{\lambda \in \Lambda }$, then the relation between the relative
Dehn function $\delta _G^{rel} $ of $G$ with respect to $\{
H_\lambda \} _{\lambda \in \Lambda }$ and the Dehn function $f_{\G
}$ of $\G $ can be expressed as follows.

\begin{lem} \label{Ak}
Suppose $G$ is a group with a reduced finite relative presentation
(\ref{G1}) with respect to a collection of the subgroups $\{
H_\lambda \} _{\lambda \in \Lambda }$. Set
\begin{equation} \label{k}
k=\max\{ \max\limits_{R\in \mathcal R} \| R\| , 4\} .
\end{equation}
Then every cycle $q$ in $\G $ admits a $k$--filling and we have
$$ A^{(k)}(q)\le (k+1) Area^{rel} (q)+2l(q).$$
\end{lem}

\begin{proof}
We consider a van Kampen diagram $\Delta $ of minimal type over
(\ref{G1}) such that $\phi (\partial \Delta )\equiv \phi (q)$. In
particular, this means that
\begin{equation}\label{d-f0}
N_{\mathcal R} (\Delta )= Area ^{rel} (q).
\end{equation}
Using Corollary \ref{sums}, we obtain
\begin{equation}\label{d-ff1}
\sum\limits_{\Pi \in S(\Delta)} l(\partial \Pi )\le
\max\limits_{R\in \mathcal R} \| R\|Area^{rel}(q)+l(q)\le  k
Area^{rel}(q)+l(q),
\end{equation}
where $S(\Delta )$ is the set of all $\mathcal S$--cells of
$\Delta $. Now we obtain a diagram $\Psi $ from $\Delta $ as
follows. Let us take a cell $\Pi \in S(\Delta )$. If $l(\partial
\Pi )>3$, we can "triangulate" $\Pi $ by considering a van Kampen
diagram $\Upsilon (\Pi )$ over (\ref{G1}) with the same boundary
label as $\Pi $ such that every cell of $\Upsilon (\Pi )$ is an
$\mathcal S$--cell whose boundary has length $3$. Obviously the
minimal number of cells in such a diagram satisfies
\begin{equation} \label{d-ff2}
Area (\Upsilon (\Pi ))\le l(\partial \Pi )-2.
\end{equation}
Then we replace the cell $\Pi $ of $\Delta $ with the diagram
$\Upsilon (\Pi )$. Doing this for all $\Pi \in S(\Delta )$, we
obtain a new diagram $\Psi $ over (\ref{G1}) in which the boundary
of every $S$--cell has length at most $3$. Combining (\ref{d-f0}),
(\ref{d-ff1}), and (\ref{d-ff2}), we obtain the following estimate
on the number of cells in $\Psi $: $$
\begin{array}{rl}
Area (\Psi )& = N_{\mathcal R} (\Delta )+\sum\limits_{\Pi \in
S(\Delta)} Area (\Upsilon (\Pi )) \\ & \\& \le Area^{rel} (q) +
\sum\limits_{\Pi \in S(\Delta)} l(\partial \Pi ) \le (k+1)
Area^{rel}(q)+l(q).
\end{array} $$

We are going to define a $k$--filling of $q$ using $\Psi $.
Obviously every $2$--cell of $\Psi $ is an $l$--gone for some
$l\le k$. The only difficulty is that $\Psi $ may not be
homeomorphic to a disk. However in this case we can transform
$\Psi $ into a simply--connected diagram $\Psi ^\prime $ by using
the so-called $0$--bordering of the contour of $\Psi $, which is a
particular case of a $0$-refinement (we refer the reader to
\cite[Ch.4, Section 11.5]{Ols-book} for the definition). Applying
this process, we add new cells to $\Psi $, the so--called
$0$--faces, such that the contour of every $0$--face has length
$4$ and the image of every $0$--face under the canonical map into
the Cayley graph $\G $ is an edge (this is one of the reasons why
we need singular maps in the definition of the $k$--filling).
Moreover, the number of additional cells equals $l(q)$. The
diagram $\Psi ^\prime $ defines a $k$--partition of $D^2$ and the
natural map $Sk^{(1)}\Psi ^\prime \to \Gamma $ such that $\partial
\Psi ^\prime $ is mapped onto $q$ gives a $k$--filling of $q$.
Thus we have
$$
A^{(k)}(q)\le Area (\Psi ^\prime)= Area (\Psi )+l(q)\le(k+1)
Area^{rel} (q)+2l(q).
$$
\end{proof}

From the above lemma we obtain

\begin{thm} \label{d-f}
Suppose that a group $G$ is finitely presented with respect to a
collection of subgroups $\{ H_\lambda \} _{\lambda \in \Lambda }$
and the Dehn relative function $\delta ^{rel}_G$ of $G$ with
respect to $\{ H_\lambda \} _{\lambda \in \Lambda })$ is
well--defined. Then the Denh function $f_{\G } $ of the Cayley
graph $\Gamma (G, X\cup\mathcal H)$  is equivalent to $\delta
^{rel}_G$.
\end{thm}

\begin{proof}
Lemma \ref{Ak} gives us the estimate
$$
f_{\G }\preceq \delta ^{rel}_G.
$$
To prove the reverse inequality, we take $k$ such that any cycle
in $\G $ admits a $k$--filling. Let $q$ be such a cycle, $l(q)\le
n$. Suppose that $\Phi :Sk^{(1)} D^2\to \G $ is a $k$--filling
having the minimal number of faces among all $k$--fillings of $q$.
Thus we have
\begin{equation} \label{d-f1}
|\Phi |\le f_{\G } (n).
\end{equation}

Using $\Phi $ we can define labels and orientations on edges of
$D^2$ as follows. If $e$ is an edge of $D^2$ and $\Phi (e)$ is an
edge of $\G $, we endow $e$ with the induced label and orientation
in the obvious way; if $\Phi (e)$ is a vertex of $\G $, then we
put $\phi (e)=1$ (the orientation does not matter in this case).
Further, let $\Pi $ be a cell of $D^2$. Since $\Phi (\partial \Pi
)$ has length at most $k$, there exists a van Kampen diagram $\Xi
(\Pi )$ over (\ref{G1}) such that $\phi (\partial \Xi (\Pi
))\equiv \phi (\partial \Pi )$ and
\begin{equation} \label{d-f2}
N_\mathcal R(\Xi (\Pi ))\le \delta ^{rel}_G(k)<\infty .
\end{equation}
Finally we can obtain a ($0$--refinement of) van Kampen diagram
over (\ref{G1}) such that $\phi (\partial \Delta )\equiv \phi (q)$
by replacing all cells $\Pi $ of $D^2$ with the corresponding
diagrams $\Xi (\Pi )$.  Evidently, from (\ref{d-f1}) and
(\ref{d-f2}) we obtain
$$
Area ^{rel}(q)\le N_\mathcal R (\Delta )\le | \Phi | \delta ^{rel}
_G(k)\le f_{\G }(n)\delta ^{rel} _G(k).
$$
This leads to the inequality $\delta ^{rel} _G\preceq f_{\G }$.
\end{proof}

\begin{cor}\label{Gammahyp}
Suppose that a group $G$ is finitely presented with respect to a
collection of subgroups $\{ H_\lambda \} _{\lambda \in \Lambda }$
and the relative Dehn function $\delta ^{rel}_G$ of $G$ with
respect to $\{ H_\lambda \} _{\lambda \in \Lambda })$ is
well--defined. Then the following conditions are equivalent.

1) $G$ is hyperbolic relative to $\{ H_\lambda \} _{\lambda \in
\Lambda }$.

2) The Cayley graph $\G $ is a hyperbolic metric space (we refer
to Section 3.1 for the definition).
\end{cor}

\begin{proof}
By Theorem \ref{d-f} the relative Dehn function of $G$ with
respect to $\{ H_\lambda \} _{\lambda \in \Lambda }$ is linear if
and only if the Dehn function of $\G $ is linear. As is
well--known, the las condition is equivalent to the hyperbolicity
of $\G $. Up to notation, the proof can be found in \cite[Ch.
III.H, Sec. 2]{BriH}; see also \cite{Lys}, \cite{Ol91},
\cite{Short}.
\end{proof}

\begin{ex}
The requirement $\delta ^{rel}_G (n) < \infty $ for any $n$ in the
statement of Corollary \ref{Gammahyp} is essential. Indeed,
consider the group $G\cong \mathbb Z$ and a subgroup $H$ of finite
index in $G$. Then obviously $\G $ has finite diameter. In
particular, $\G $ is hyperbolic. However, $\delta ^{rel} _G$ is
not well--defined by Proposition \ref{malnorm}.
\end{ex}

\begin{defn}
A function $f:\mathbb N\to \mathbb N$ is said to be {\it
subquadratic } if $f=o(n^2)$ as $n\to \infty $.
\end{defn}

Recall that every finitely presented group with subquadratic Dehn
function is hyperbolic. This insight is due to M. Gromov
\cite{Gr1}, and was clarified by Ol'shanskii \cite{Ol91} and
others (see \cite{Bow95}, \cite{Paps95} and references therein).
Our next goal is to show that if $G$ has a subquadratic relative
Dehn function $\delta ^{rel}_G$ with respect to a collection $\{
H_\lambda \} _{\lambda \in \Lambda }$, then , in fact, $\delta
^{rel}_G$ is linear and $\G $ is a hyperbolic metric space.

\begin{cor} \label{Ghyp}
Suppose that a group $G$ is finitely presented with respect to a
collection of subgroups $\{ H_\lambda \} _{\lambda \in \Lambda }$
and the relative Denh function $\delta ^{rel}_G$ of $G$ with
respect to $\{ H_\lambda \} _{\lambda \in \Lambda })$, is
well--defined. Then the following conditions are equivalent.

1) $\delta ^{rel}_G $ is subquadratic.

2) $\delta ^{rel}_G$ is linear.

\end{cor}

\begin{proof}
If $\delta ^{rel}_G $ is subquadratic, then the Dehn function of
$\G $ is equivalent to $\delta ^{rel}_G $ by Theorem \ref{d-f}. It
is well--known that if a graph possess a subquadratic Dehn
function, it possess a linear one (see \cite[Ch. III.H, Sec.
2]{BriH}, \cite{Bow95}). Applying Theorem \ref{d-f} again, we
obtain that $\delta ^{rel}_G$ is linear.
\end{proof}


\chapter{Geometry of finitely generated relatively hyperbolic groups}


\section{Conventions and notation}

In this chapter we study the geometry of the Caley graph of a
finitely generated in the usual non--relative sense group $G$
which is hyperbolic relative to a collection $H_1, \ldots , H_m$
of subgroups. (Recall that by Corollary \ref{Lambda}, if a
finitely generated group is hyperbolic relative to a collection of
subgroups, then the collection is finite). Throughout the next
three sections we accept the following technical agreements about
$G$ and $H_1, \ldots , H_m$.

(i) {\it $G$ is represented by a finite relative presentation
\begin{equation}
G=\langle X, H_1, \ldots , H_m\; | \; R=1, R\in \mathcal R\rangle
\label{Gfg}
\end{equation}
with respect to $\{ H_1, \ldots , H_m\} $ and the relative Dehn
function of (\ref{Gfg}) satisfies
\begin{equation} \label{iidelta}
\delta^{rel}_G(n)\le Ln
\end{equation}
for some constant $L$; also, as in the previous chapter, we set
$$M=\max\limits_{R\in \mathcal R} ||R||.$$
By technical reasons it is convenient to increase $L$ in order to
satisfy the inequality $$ML>1.$$}

(ii) {\it  $G$ is generated by the set $X$ in the ordinary sense.}

(iii) {\it The set $X$ is chosen in such a way that the following
lemma holds}.

\begin{lem}\label{31}
Let $G$ be a finitely generated group hyperbolic relative to a
collection of subgroups $\{ H_1, \ldots , H_m\} $. Then there
exists a finite generating set $X$ of $G$ satisfying the following
condition. Let $q$ be a cycle in $\G $, $p_1, \ldots , p_k$ a
certain set of isolated $H_i $--components of $q$. Then
\begin{equation}\label{311}
\sum\limits_{i=1}^k \dx ((p_i)_-, (p_i)_+)\le MLl(q).
\end{equation}
\end{lem}

Let us show that we can always ensure the fulfillment of
(i)--(iii). We start with an arbitrary finite reduced relative
presentation
\begin{equation}\label{Gfgstart}
G=\langle X^\prime , H_1, \ldots , H_m\; | \; R=1, R\in \mathcal
R^\prime \rangle .
\end{equation}
Let $\Omega $ be the set given by Definition \ref{defofOmega}.
Note that $\Omega $ is finite in our case. We set $$X=X^\prime
\cup \Omega $$ and $$\mathcal R=\mathcal R^\prime \cup \{ \omega
W_\omega^{-1}, \; \omega \in \Omega \} ,$$ where $W_\omega $ is a
fixed word over $X^\prime $ representing the element $\omega $ in
$G$. It is clear that for $X$ and $\mathcal R$ chosen in this way,
the presentation (\ref{Gfg}) can be obtained from (\ref{Gfgstart})
by a finite number of Tietze transformations. It remains to prove
Lemma \ref{31}.

Let $u_i$ (respectively $v_i$) denote $(p_i)_-$ (respectively
$(p_i)_+$) regarded as an element of $G$. By Lemma \ref{Omega}, we
have
$$\sum\limits_{i=1}^k |u_i^{-1}v_i|_\Omega \le M^\prime L^\prime
l(q),$$ where
$$M=\max\limits_{R\in \mathcal R^\prime} ||R||\le M,$$
and the relative Dehn function of the presentation
(\ref{Gfgstart}) does not exceed $L^\prime n$. Increasing the
constant $L$ in (\ref{iidelta}) if necessary, we may assume that
$L^\prime \le L$. Thus we obtain
$$
\sum\limits_{i=1}^k \dx ((p_i)_-, (p_i)_+)\le\sum\limits_{i=1}^k
|u_i^{-1}v_i|_X\le \sum\limits_{i=1}^k |u_i^{-1}v_i|_\Omega \le
M^\prime L^\prime l(q)\le MLl(q).
$$

\begin{defn}\label{hypms}
Recall that a metric space $Y$ is called {\it $\delta
$--hyperbolic} (or simply {\it hyperbolic}) if it satisfies the
following {\it Rips condition}. For any geodesic triangle, each
side of the triangle belongs to the union of the closed $\delta
$--neighborhood of the other two sides.
\end{defn}

By Corollary \ref{Gammahyp}, the Cayley graph $\G $ is hyperbolic.
In what follows we denote by $\delta $ the hyperbolicity constant
of $\G $.


\section{Properties of quasi--geodesics}


First of all we recall some facts concerning quasi--geodesics in
hyperbolic metric spaces. Below we assume all paths under
consideration to be rectifiable (i.e., to have finite length).

\begin{defn}
A path $q$ in a metric space $Y$ is said to be $(\lambda,
c)$--{\it quasi--geodesic} for some $\lambda \ge 1$, $c\ge 0$, if
for every subpath $q$ of $p$ the inequality $$ l(q)\le \lambda\,
dist(q_-,q_+) +c $$ holds.
\end{defn}

The following useful lemma is quite obvious.

\begin{lem}\label{plus1}
Let $p$ be a $(\lambda , c)$--quasi--geodesic in a metric space
$Y$, $e$ a path of length $k$ such that $e_-=p_+$. Then the path
$pe$ is $(\lambda , c+(\lambda +1)k)$--quasi--geodesic.
\end{lem}

\begin{proof}
Let $q$ be a subpath of $pe$, $q_0$ the maximal common subpath of
$q$ and $p$. Then
$$
\begin{array}{rl}
l(q) & \le l(q_0)+k \le \lambda \dxh ((q_0)_-, (q_0)_+) +c+k\\ &
\\ & \le \lambda (\dxh ((q)_-, (q)_+) +k)+c+k.
\end{array}
$$
\end{proof}

Next result is well known and can be found in \cite{Gr1} or
\cite{GhH}.

\begin{lem}\label{qg}
For any $\delta \ge 0$, $\lambda \ge 1$, $c\ge 0$, there exists a
constant $H=H(\delta , \lambda , c)$ with the following property.
If $Y$ is a $\delta $--hyperbolic space and $p, q$ are $(\lambda ,
c)$--quasi--geodesic paths in $Y$ with same endpoints, then $p$
and $q$ belong to the closed $H$--neighborhoods of each other.
\end{lem}

From the definition of a hyperbolic space, one obtains

\begin{lem}
Let $Y$ be a $\delta $--hyperbolic metric space, $Q=p_1p_2p_3p_4$
a geodesic quadrangle in $Y$. Then each side of $Q$ belongs to the
closed $2\delta $--neighborhood of the other three sides.
\end{lem}

\begin{cor}\label{quad-rel}
Let $Y$ be a $\delta $--hyperbolic space and $p, q$ are geodesic
paths in $Y$ such that $dist (p_-, q_-)\le k $ and $dist (p_+,
q_+)\le k$, then $p$ and $q$ belong to the closed $(k+2\delta
)$--neighborhood of each other.
\end{cor}

From the above corollary and Lemma \ref{qg}, we can easily derive

\begin{lem} \label{K}
For any $\delta \ge 0 $, $\lambda \ge 1$, $c\ge 0$, and $k\ge 0$,
there exists a constant $K=K(\delta , \lambda , c,k)$ such that
the following condition holds. Suppose that $Y$ is a $\delta
$--hyperbolic space and $p, q$ are $(\lambda ,
c)$--quasi--geodesic paths in $Y$ such that $dist (p_-, q_-)\le k
$ and $dist (p_+, q_+)\le k$; then $p$ and $q$ belong to the
closed $K$--neighborhoods of each other.
\end{lem}

Given the data described in Section 3.1, there are two Cayley
graphs, namely $\Gamma (G, X)$, the Cayley graph of $G$ with
respect to the generating set $X$, and $\Gamma (G, X\cup\mathcal
H)$, the Cayley graph of $G$ with respect to the generating set
$X\cup \mathcal H.$ Obviously we have a natural embedding $$\iota
: \Gamma (G, X)\to \G ,$$ which is bijective on the set of
vertices. For simplicity we will identify $\Gamma (G, X)$ with its
image in $\Gamma (G, X\cup\mathcal H)$ under the embedding $\iota
$.

Assuming the length of each edge of $\G $ to be equal to $1$, we
obtain combinatorial metrics $dist _X$ and $\dxh $ in $\Gamma
(G,X)$ and $\G $ respectively. When speaking about geodesics
(quasi--geodesics) in $\Gamma (G,X)$ or in $\G $ we always mean
geodesics (quasi--geodesics) with respect to corresponding
combinatorial metric. Note that the restrictions of $\dx $ and
$\dxh $ on the vertex sets $V(\Gamma (G,X)=V(\G )=G$ coincide with
the word metric on $G$ with respect to the generating sets $X$ and
$X\cup \mathcal H$.

\begin{defn} \label{backtrac}
A path $p$ in $\G $ is called a {\it path without backtracking} if
for any $i=1, \ldots , m$, every $H_i$--component of $p$ is
isolated.
\end{defn}

\begin{defn}
Let $p$ be a path in $\G $, $v$ a vertex of a component $s$ of
$p$. If $v\ne s_-$ and $v\ne s_+$, we say that $v$ is an inner
vertex of $s$. A vertex $u$ of $p$ is called {\it non--phase}, if
$u$ is a inner vertex of some component of $p$. All other vertices
of $p$ are called {\it phase}.
\end{defn}

\begin{defn}
A path $p$ in $\G $ is said to be {\it locally minimal} if for any
$i=1, \ldots ,m$, every $H_i$--component of $p$ has length $1$ or,
equivalently, every vertex of $p$ is phase.
\end{defn}

We need the following simple observation.

\begin{lem}\label{LM}
Let $p$ be a path in $\G $. Then there exists a locally minimal
path $\hat p$ having the same set of phase vertices as $p$.
Moreover, if $p$ is $(\lambda, c)$--quasi--geodesic, then $\hat p$
is $(\lambda , c)$--quasi--geodesic; if $p$ is a path without
backtracking, then so is $\hat p$.
\end{lem}

\begin{proof}
Given an $H_i$--component $s$ of $p$, we can replace $s$ with the
single edge of $\G $ having the label $\phi (e)=_{H_i}\phi (s)$.
Doing this for all components of $p$ we obtain the desirable path
$\hat p$. The verification of the additional properties is
straightforward.
\end{proof}

\begin{defn}
Let $p$, $q$  be two paths in $\G $. We say that $p$ and $q$ are
{\it $k$--similar} if
\begin{equation}\label{QG1}
\max \{ \dx (p_-, q_-),\;  \dx (p_+, q_+)\} \le k.
\end{equation}
\end{defn}

The next proposition is an improved version of Lemma \ref{K} for
relatively hyperbolic groups.

\begin{defn}
We say that two paths $p$ and $q$ in $\G $ are $k$--similar if
$$
\max \{ \dx (p_-, q_-),\;  \dx (p_+, q_+)\} \le k.
$$
\end{defn}

\begin{prop} \label{QG}
For any $\lambda \ge 1 $, $c\ge 0$, $k\ge 0$ there exists a
constant $\varepsilon = \varepsilon (\lambda , c, k)
>0$ having the following property. Let $p$ and $q$ be two
$k$--similar $(\lambda , c)$--quasi--geodesic paths in $\G $ (with
respect to the relative metric $dist _{X\cup \mathcal H}$) such
that $p$ is a path without backtracking. Then for any phase vertex
$u$ of $p$ there exists a phase vertex $v$ of $q$ such that
\begin{equation}\label{QGmain}
dist_X(u,v)\le \varepsilon .
\end{equation}
\end{prop}

\begin{proof}
Before proving the theorem, we note that (\ref{QGmain}) is much
stronger that the inequality $\dxh (u,v)\le \varepsilon $, which
follows from Lemma \ref{K}. In view of Lemma \ref{LM}, it is
sufficient to prove the proposition for locally minimal
quasi--geodesic paths. Thus we assume that every vertex of $p$ and
$q$ is phase.

Let $u$ be a (phase) vertex of $p$. Recall that $\G $ is a $\delta
$--hyperbolic space. We set $K_0=K(\delta , \lambda , c, k)$ and
$K=K(\delta , \lambda , c,K_0) +1/2$, where  $K(\delta , \lambda ,
c, k)$ and $K(\delta , \lambda , c,K_0) $ are constants provided
by Lemma \ref{K}. Without loss of generality we may assume that
$K\in \mathbb N$ and
\begin{equation} \label{KKk}
K\ge K_0\ge k.
\end{equation}

Let us choose vertices $u_1$, $u_2$ on $p$ as follows. If $\dxh
(p_-, u)\le 2K$ (respectively $\dxh (u, p_+)\le 2K$), we put
$u_1=p_- $ (respectively $u_2=p_+$). If $\dxh (p_-, u)> 2K$
(respectively $\dxh (u, p_+)> 2K$), we take $u_1$ on the segment
$[p_-, u]$ of the path $p$ (respectively $u_2$ on the segment $[u,
p_+]$ of the path $p$) such that
\begin{equation} \label{Rips0}
\dxh (u, u_i)=2K
\end{equation}
for $i=1$ (respectively $i=2$).

Further, by Lemma \ref{K}, there exist two points $v_1, v_2$ of
$q$ such that
\begin{equation} \label{Rips1}
\dxh  (v_i, u_i)\le K_0, \;\;\; i=1,2.
\end{equation}
Without loss of generality we may assume that $v_1, v_2$ are
vertices of $q$. Moreover, we assume that in case $u_1=p_-$
(respectively $u_2=p_+$) the vertex $v_1$ coincides with $q_-$
(respectively the vertex $v_2$ coincides with $q_+$).

We denote by $p_1, p_2$ (respectively by $q_0$) the segments of
$p$ (respectively the segment of $q$) such that $(p_1)_-=u_1$,
$(p_1)_+=u$, $(p_2)_-=u$, $(p_2)_+=u_2$ (respectively
$(q_0)_-=v_1$, $(q_0)_+=v_2$). We also denote by $o_1, o_2$ the
paths in $\G $ such that $(o_i)_-=u_i$, $(o_i)_+=v_i$, $i=1,2$
(see Fig. \ref{QH}) and $o_1$, $o_2$ are chosen according to the
following agreement. If $u_1=p_-$ (respectively $u_2=p_+$), then
$o_1$ (respectively $o_2$) is a geodesic path in $\Gamma (G,X)$.
If $u_1\neq p_1$ (respectively $u_2\neq p_+$), then $o_1$
(respectively $o_2$) is geodesic in $\G $.  It follows from
(\ref{QG1}) and (\ref{KKk}) in the first case and from
(\ref{Rips1}) in the second case that
\begin{equation}
\label{QG2} l(o_i)\le K_0.
\end{equation}

Let $V$ denote the set of all vertices $z$ on $q_0$ that are
closest to $u$, i.e., satisfy the condition
$$
dist _{X\cup \mathcal H} (u,z)= \min\limits_{v\in q_0}\, \dxh
(u,v),
$$
where $v$ ranges among all vertices of $q_0$.  Taking into account
(\ref{Rips1}) and  Lemma \ref{K} we obtain
\begin{equation} \label{dzu}
\dxh (u,z)\le K(\delta , \lambda , c, K_0)+1/2=K
\end{equation}
for any $z\in V$. To each $z \in V$, we associate the set $O(z)$
of all a geodesic paths $o$ such that $o_-=u$, $o_+=z$. Each of
the paths $o\in O(z)$ cuts the cycle $p_1p_2o_2q_0^{-1}o_1^{-1}$
into two parts denoted $c_1$ and $c_2$. More precisely, let
$$
c_1=oq_1^{-1}o_1^{-1}p_1, \;\;\;\;\; {\rm and} \;\;\;\;\;
c_2=oq_2o_2^{-1}p_2^{-1},
$$
where $q_1=[v_1, z]$, $q_2=[z, v_2]$ are segments of $q_0$. To
prove the proposition we need a few auxiliary lemmas.

\begin{figure}
\begin{picture}(120,50)(0,10)
\qbezier(0,40)(60,70)(120,40) \qbezier(0,30)(60,0)(120,30)
\multiput(0,30)(120,0){2}{\circle*{1}}
\multiput(0,40)(120,0){2}{\circle*{1}}

\put(3,31){$p_-$} \put(114,31){$p_+$} \put(3,39){$q_-$}
\put(114,39){$q_+$}

\put(60,15){\circle*{1}} \put(60,55){\circle*{1}}

\put(25,20){\circle*{1}} \put(25,50){\circle*{1}}
\put(95,20){\circle*{1}} \put(95,50){\circle*{1}}

\put(23,17){$u_1$} \put(23,52){$v_1$} \put(96,17){$u_2$}
\put(96,52){$v_2$} \put(59,12){$u$} \put(59,57){$z$}

\qbezier(25,20)(35,35)(25,50) \qbezier(95,20)(85,35)(95,50)

\put(60,15){\line(0,1){40}} \put(0,30){\line(0,1){10}}
\put(120,30){\line(0,1){10}}

\put(33,18){\vector(4,-1){1}} \put(33,52){\vector(4,1){1}}
\put(87,18){\vector(4,1){1}} \put(87,52){\vector(4,-1){1}}

\put(35,20){$p_1$} \put(35,49){$q_1$} \put(83,20){$p_2$}
\put(83,49){$q_2$}

\put(30.1,34.6){\vector(0,1){1}} \put(90.05,34.6){\vector(0,1){1}}
\put(60,35){\vector(0,1){1}}

\put(32,35){$o_1$} \put(86,35){$o_2$} \put(62,35){$o$}

\put(47,11){$p$} \put(47,58){$q$}

\put(47,35){$c_1$} \put(72,35){$c_2$}
\end{picture}
\caption{} \label{QH}
\end{figure}

\begin{lem}\label{Rips39}
Let $z\in V$, $o\in O(z)$, and let $s$ be an $H_i$--component of
the path $o$ for a certain $i$. Then for $j=1,2,$ there exist no
$H_i$--components of $o_j$ connected to $s$.
\end{lem}

\begin{proof}
First assume that $s$ is connected to an $H_i$--component $t$ of
$o_1$ (see Fig. \ref{QH-1}). This means that
$$
\dxh (s_-, t_-)\le 1.
$$
Using (\ref{QG2}), (\ref{dzu}), and (\ref{KKk}), we obtain
$$
\begin{array}{rl}
\dxh (u, u_1)& \le \dxh (u,s_-)+ \dxh (s_-, t_-)+ \dxh (t_-, u_1)\\
& \\ & \le \left( \dxh (u, z)-1\right) + 1 + (l(o_1)-1)
\\ & \\ &  < K+K_0-1<2K.
\end{array}
$$

\begin{figure}
\begin{picture}(120,50)(0,10)
\qbezier(0,40)(60,70)(120,40) \qbezier(0,30)(60,0)(120,30)
\multiput(0,30)(120,0){2}{\circle*{1}}
\multiput(0,40)(120,0){2}{\circle*{1}}

\put(3,31){$p_-$} \put(114,31){$p_+$} \put(3,39){$q_-$}
\put(114,39){$q_+$}

\put(60,15){\circle*{1}} \put(60,55){\circle*{1}}

\put(25,20){\circle*{1}} \put(25,50){\circle*{1}}
\put(95,20){\circle*{1}} \put(95,50){\circle*{1}}

\put(23,17){$u_1$} \put(23,52){$v_1$} \put(96,17){$u_2$}
\put(96,52){$v_2$} \put(59,12){$u$} \put(59,57){$z$}

\qbezier(25,20)(35,35)(25,50) \qbezier(95,20)(85,35)(95,50)

\put(60,15){\line(0,1){40}} \put(0,30){\line(0,1){10}}
\put(120,30){\line(0,1){10}}

\put(33,18){\vector(4,-1){1}} \put(33,52){\vector(4,1){1}}
\put(87,18){\vector(4,1){1}} \put(87,52){\vector(4,-1){1}}

\put(35,20){$p_1$} \put(35,49){$q_1$} \put(83,20){$p_2$}
\put(83,49){$q_2$}

\put(90.05,34.6){\vector(0,1){1}}

\put(86,35){$o_2$} \put(47,11){$p$} \put(47,58){$q$}

\thicklines \put(60,35){\line(0,1){10}} \put(60,35){\circle*{1}}
\put(60,45){\circle*{1}} \put(60,40){\vector(0,1){1}}
\put(30.1,34.6){\vector(0,1){1}}

\qbezier(28.8,28)(31,35)(28.8,42) \put(28.9,42){\circle*{1}}
\put(28.9,28){\circle*{1}} \qbezier(28.8,28)(38,40)(60,35)
\put(62,35){$s_-$} \put(62,45){$s_+$} \put(24,28){$t_-$}
\put(24,42){$t_+$} \put(32,35){$t$} \put(57,40){$s$}
\end{picture}
\caption{} \label{QH-1}
\end{figure}

By our choice of $u_1$ and $v_1$, this inequality implies
$u_1=p_-$, $v_1=q_-$. Therefore, according to our choice of $o_i$,
$o_1$ is a path in $\Gamma (G,X)$. Thus $o_1$ contains no
$H_i$--components at all. A contradiction. The case $j=2$ is
completely analogous.
\end{proof}

For convenience, we give one more auxiliary definition. Given
$z\in V$ and $o\in O(z)$, an $H_i$--component $s$ of the path $o$
is called an {\it ending component} if $s$ contains $z$; for
otherwise $s$ is called a {\it non--ending} component. It can
happen that $o$ has no ending component, since the last edge of
$o$ can be labelled by a letter from $X$.

\begin{lem}\label{non-end}
Let $z\in V$, $o\in O(z)$, and let $s$ be a non--ending
$H_i$--component of the path $o$ for a certain $i$. Then for any
$j=1,2$ the following holds. If $s$ is not an isolated
$H_i$--component of $c_j$, then there exists an $H_i$--component
$t$ of $p_j$ that is connected to $s$.
\end{lem}

\begin{proof}
Again we consider the case $j=1$ only. According to Lemma
\ref{Rips39}, there are no $H_i$--component of $o_1$ connected to
$s$. Assume that there exists an $H_i$--component $t$ of $q_1$
connected to $s$ (see Fig. \ref{QH-2}). Since $s$ is a non--ending
component, we have $dist_{X\cup\mathcal H}(s_-, z)\ge 2.$ However,
$dist_{X\cup\mathcal H}(s_-, t_-)\le 1$. Hence
$$
\begin{array}{rl}
\dxh (u,t_-) & \le \dxh (u, s_-) + \dxh (s_-, t_-)\\
&\\ & < \dxh (u, s_-) +\dxh (s_-, z) \\ & \\= \dxh (u,z)
\end{array}
$$

\begin{figure}
\begin{picture}(120,50)(0,10)
\qbezier(0,40)(60,70)(120,40) \qbezier(0,30)(60,0)(120,30)
\multiput(0,30)(120,0){2}{\circle*{1}}
\multiput(0,40)(120,0){2}{\circle*{1}}

\put(3,31){$p_-$} \put(114,31){$p_+$} \put(3,39){$q_-$}
\put(114,39){$q_+$}

\put(60,15){\circle*{1}} \put(60,55){\circle*{1}}

\put(25,20){\circle*{1}} \put(25,50){\circle*{1}}
\put(95,20){\circle*{1}} \put(95,50){\circle*{1}}

\put(23,17){$u_1$} \put(23,52){$v_1$} \put(96,17){$u_2$}
\put(96,52){$v_2$} \put(59,12){$u$} \put(59,57){$z$}

\qbezier(25,20)(35,35)(25,50) \qbezier(95,20)(85,35)(95,50)

\put(60,15){\line(0,1){40}} \put(0,30){\line(0,1){10}}
\put(120,30){\line(0,1){10}}

\put(33,18){\vector(4,-1){1}} \put(87,18){\vector(4,1){1}}
\put(87,52){\vector(4,-1){1}}

\put(35,20){$p_1$}  \put(83,20){$p_2$} \put(83,49){$q_2$}

\put(30.1,34.6){\vector(0,1){1}} \put(90.05,34.6){\vector(0,1){1}}

\put(32,35){$o_1$} \put(86,35){$o_2$}

\put(47,11){$p$}

\thicklines \put(60,35){\line(0,1){10}} \put(60,35){\circle*{1}}
\put(60,45){\circle*{1}} \put(60,40){\vector(0,1){1}}

\put(62,35){$s_-$} \put(62,45){$s_+$} \put(42,50){$t_-$}
\put(54,51){$t_+$} \put(49,56){$t$} \put(57,40){$s$}
\put(45,54){\circle*{1}} \put(55,54.9){\circle*{1}}
\qbezier(45,54)(50,54.7)(55,54.9) \put (53,54.8){\vector(1,0){1}}
\qbezier(45,54)(50,43)(60,35)
\end{picture}
\caption{} \label{QH-2}
\end{figure}
This contradicts to the assumption that $z$ is a closest vertex to
$u$ on $q_0$.
\end{proof}

From Lemma \ref{non-end} we immediately obtain

\begin{cor} \label{non-end1}
Let $z\in V$, $o\in O(z)$. Then for any $i=1, \ldots , m$, every
non--ending $H_i$--component of the path $o$ is an isolated
$H_i$--component of at least one of the cycles $c_1, c_2$.
\end{cor}

\begin{proof}
Indeed suppose that some non--ending $H_i$--component $s$ of $o$
is not isolated in both $c_1$ and $c_2$. By Lemma \ref{non-end},
this means that there are $H_i$--components $t_1$ and $t_2$ of
$p_1$ and $p_2$ respectively that are connected to $s$. In
particular, $t_1$ is connected to $t_2$. However this contradicts
to the assumption that $p$ is a locally minimal path without
backtracking.
\end{proof}

\begin{lem} \label{end}
Let $j=1$ or $j=2$. Suppose that for any $z\in V$ and any $o\in
O(z)$, $o$ has the ending component, which is not isolated in
$c_j$; then for any $z\in V$ and any $o\in O(z)$, the ending
component of $o$ is connected to a component of $p_j$ for
corresponding $i$.
\end{lem}

\begin{proof}
For definiteness assume $j=1$. We proceed by induction on $\dxh
(z,v_1)$. In case $\dxh(z, v_1)=0$,  Lemma \ref{end} is obvious.
Indeed, assume that the ending component $s$ of $o\in O(z)$ is not
isolated in $c_1$. By Lemma \ref{Rips39}, there are no components
of $o_1$ connected to $s$. As $o$ is a geodesic path, $s$ is an
isolated component of $o$. Therefore, the only possible case is
that there exists a component of $p_1$ connected to $s$, since
$q_1$ is trivial in this case.

\begin{figure}
\begin{picture}(120,50)(0,10)
\qbezier(0,40)(60,70)(120,40) \qbezier(0,30)(60,0)(120,30)
\multiput(0,30)(120,0){2}{\circle*{1}}
\multiput(0,40)(120,0){2}{\circle*{1}}

\put(3,31){$p_-$} \put(114,31){$p_+$} \put(3,39){$q_-$}
\put(114,39){$q_+$}

\put(60,15){\circle*{1}} \put(60,55){\circle*{1}}

\put(25,20){\circle*{1}} \put(25,50){\circle*{1}}
\put(95,20){\circle*{1}} \put(95,50){\circle*{1}}

\put(23,17){$u_1$} \put(23,52){$v_1$} \put(96,17){$u_2$}
\put(96,52){$v_2$} \put(59,12){$u$}

\qbezier(25,20)(35,35)(25,50) \qbezier(95,20)(85,35)(95,50)

\put(60,15){\line(0,1){40}} \put(0,30){\line(0,1){10}}
\put(120,30){\line(0,1){10}}

\put(33,18){\vector(4,-1){1}} \put(87,18){\vector(4,1){1}}
\put(87,52){\vector(4,-1){1}}

\put(35,20){$p_1$}  \put(83,20){$p_2$} \put(83,49){$q_2$}

\put(30.1,34.6){\vector(0,1){1}} \put(90.05,34.6){\vector(0,1){1}}

\put(32,35){$o_1$} \put(86,35){$o_2$}

\put(47,11){$p$}

\thicklines \put(60,35){\line(0,1){20}} \put(60,35){\circle*{1}}
\put(60,55){\circle*{1}} \put(60,45){\vector(0,1){1}}

\put(62,35){$s_-$} \put(58,57){$s_+=z$} \put(42,50){$t_-$}
\put(54,51){$t_+$} \put(49,56){$t$} \put(57,45){$s$}
\put(45,54){\circle*{1}} \put(55,54.9){\circle*{1}}
\qbezier(45,54)(50,54.7)(55,54.9) \put (53,54.8){\vector(1,0){1}}
\qbezier(45,54)(50,43)(60,35) \put (51.2,44){\vector(-2,3){1}}
\end{picture}
\caption{} \label{QH-3}
\end{figure}

Now assume that $\dxh (z,v_1)\ge 1$. Arguing as above, we can
easily show that it is sufficient to consider the case when $s$ is
connected to a component $t$ of $q_1$. Let $c$ be the edge of $\G
$ labelled by a letter from $\widetilde H_i\setminus \{ 1\} $ for
corresponding $i$  such that $c_-=s_-$ and $c_+=t_-$ (see Fig.
\ref{QH-3}). Since $\dxh (u,t_-)=\dxh (u,z)$, we have $t_-\in V$.
Consider the new path $o^\prime\in O(t_-)$ defined by $o^\prime
=[u,s_-]c$, where $[u,s_-]$ is the segment of the path $o$.
Obviously the distance between $t_-$ and $v_1$ is smaller that
$\dxh (v_1, z)$ and $c$ is the ending component of $o^\prime $.
Thus, by the inductive assumption, there is an $H_i$--component
$r$ of $p_1$ connected to $c$. As $c$ is connected to $s$, we get
what we need. The lemma is proved.
\end{proof}

Arguing as in the proof of Corollary \ref{non-end1}, we
immediately obtain

\begin{cor} \label{end1}
There exists a vertex $z\in V$ and a path $o\in O(z)$ such that
either the last edge of $o$ is labelled by a letter from $X$ or
the ending component of $o$ is isolated in at least one of the
cycles $c_1, c_2$.
\end{cor}

Now let us return to the proof of Proposition \ref{QG}. By
Corollaries \ref{non-end1} and \ref{end1}, there exists a vertex
$z$ on $q_0$ that satisfies (\ref{dzu}) and a (geodesic) path $o$
connecting $u$ and $z$ such that for any $i$, any $H_i$--component
of $o$ is isolated in $c_1$ or in $c_2$. Applying Lemma \ref{31},
we obtain
\begin{equation}\label{dcomp}
\dx (s_-, s_+)\le MLl(c_j)
\end{equation}
for each $H_i$--component $s$ of $o$. Let us estimate the length
of $c_j$. Since $p$ and $q$ are $(\lambda ,c)$--quasi--geodesics,
we have the following bounds on the lengths of the paths $p_j$:
\begin{equation} \label{Rips01}
l(p_j)\le \lambda \dxh (u_1, u)+c\le 2\lambda K+c.
\end{equation}
Further,
\begin{equation} \label{Rips02}
\begin{array}{rl}
l(q_j) & \le l(q_0)\le \lambda \dxh (v_1, v_2)+c\\ & \\
& \le \lambda ( \dxh (v_1, u_1) +\dxh (u_1, u_2) +\dxh (u_2,
v_2))+c\\ & \\
& \le 6\lambda K +c.
\end{array}
\end{equation}
Finally, by combining (\ref{KKk}), (\ref{QG2}), (\ref{dzu}),
(\ref{Rips01}), and (\ref{Rips02}), we obtain
\begin{equation} \label{Rips03}
l(c_j)\le l(p_j)+l(o)+l(q_j)+l(o_j)\le 8\lambda K +2c +2K.
\end{equation}
Inequalities (\ref{dcomp}) and (\ref{Rips03}) imply
\begin{equation}\label{fin}
\dx (u,z) \le \dxh (u,z) ML\max\limits_{j=1,2} l(c_j) \le
KLM(8\lambda K +2c +2K).
\end{equation}
It remains to assume $\varepsilon $ to be equal to the right hand
side of (\ref{fin}).
\end{proof}

Proposition \ref{QG} allows one to show that if $p$ and $q$ are
quasi--geodesic paths without backtracking in $\G $ with 'close'
endpoints, then for any 'long' component $s$ of $p$ there exists a
component $t$ of $q$ connected to $s$ and, moreover, the
$X$--distances between corresponding endpoints of $s$ and $t$ is
'small'.

\begin{lem} \label{lin-BCP}
For any $\lambda \ge 1$, $c\ge 0$, $k\ge 0$, there are $C=C(
\lambda, c, k)$ and $D=D(\lambda , c, k)$ satisfying the following
conditions. Let $p$ and $q$ be a pair of $k$--similar $(\lambda ,
c)$--quasi--geodesics in $\G $ such that $p$ is a path without
backtracking. Then for any $i=1, \ldots ,m$ and any
$H_i$--component $s$ of $p$ satisfying the condition $\dx (s_-,
s_+) >C$, there exists an $H_i$--component $t$ of $q$ such that
$t$ is connected to $s$;
\end{lem}

\begin{proof}
We set
\begin{equation}\label{lbsp1}
C=LM(1+\lambda (2\varepsilon +1)+c+2\varepsilon ),
\end{equation}
where $\varepsilon = \varepsilon (\lambda , c,k)$ is provided by
Proposition \ref{QG}. Let $s$ be an $H_i$--component of $p$ such
that $\dx (s_-, s_+)>C$. By Proposition \ref{QG}, there exist
vertices $w_1, w_2$ on $q$ such that $ \dx (s_-, w_1)\le
\varepsilon $ and $\dx (s_+, w_2)\le \varepsilon $. Since $q$ is
$(\lambda , c)$--quasi--geodesic, the length of the segment
$q_0=[w_1, w_2]$ of the path $q$ satisfies the inequality
\begin{equation}\label{lbsp3}
\begin{array}{rl}
l(q_0) & \le \lambda \dxh (w_1, w_2)+c \\ & \\  & \le \lambda
(\dxh (w_1, s_-) +\dxh (s_-, s_+) + \dxh (s_+, w_2)) +c
\\ & \\ & \le
\lambda (2\varepsilon +1) +c.
\end{array}
\end{equation}

\begin{figure}
\begin{picture}(120,50)(0,10)
\qbezier(0,40)(60,70)(120,40) \qbezier(0,30)(60,0)(120,30)
\multiput(0,30)(120,0){2}{\circle*{1}}
\multiput(0,40)(120,0){2}{\circle*{1}} \put(3,31){$p_-$}
\put(114,31){$p_+$} \put(3,39){$q_-$} \put(114,39){$q_+$}

\put(25,20){\circle*{1}} \put(25,50){\circle*{1}}
\put(95,20){\circle*{1}} \put(95,50){\circle*{1}}

\put(59,57){$q_0$} \put(59,12){$s$} \put(60,55){\vector(1,0){1}}

\put(23,17){$s_-$} \put(23,52){$w_1$} \put(96,17){$s_+$}
\put(96,52){$w_2$}

\qbezier(25,20)(35,35)(25,50) \qbezier(95,20)(85,35)(95,50)

\put(30.1,34.6){\vector(0,1){1}} \put(90.05,34.6){\vector(0,1){1}}

\put(32,35){$r_1$} \put(86,35){$r_2$}

\put(10,20){$p$} \put(10,50){$q$}

\put(0,30){\line(0,1){10}} \put(120,30){\line(0,1){10}}

\thicklines \qbezier(25,20)(60,10)(95,20)
 \put(59,15){\vector(1,0){1}}

\end{picture}
\caption{} \label{BCP0}
\end{figure}

We fix some paths $r_1$ and $r_2$ such that $(r_1)_-=s_-$,
$(r_1)+=w_1$, $(r_2)_-=s_+$, $(r_2)_+=w_2$, and $r_1, r_2$ are
geodesic paths in $\Gamma (G, X)$ (not in $\G $). By Proposition
\ref{QG},
\begin{equation}\label{lri} l(r_i)\le \varepsilon ,
\;\;\; i=1,2.
\end{equation}
Notice that since the labels of $r_1$, $r_2$ consist of letters
from $X$, $r_1$ and $r_2$ contain no $H_i$--components. Suppose
that there exist no $H_i$--components of $q_0$ connected to $s$.
Then $s$ is an isolated $H_i$--component of the cycle
$$
b=sr_2q_0^{-1}r_1^{-1}.
$$
By Lemma \ref{31} and inequalities (\ref{lbsp3}), (\ref{lri}), we
have
$$
\begin{array}{rl}
\dx (s_-, s_+) \le & LM l(b) \le LM
(l(s)+l(q_0) +l(r_1) + l(r_2)) \\& \\
& \le LM (1+\lambda (2\varepsilon +1)+c +2\varepsilon )=C,
\end{array}
$$
that contradicts to the assumption $\dx (s_-, s_+) >C.$
\end{proof}

The next lemma shows that if there exist two connected components
$s$ and $t$ of two quasi--geodesics in $\G $ with 'close'
endpoints, then the endpoints of $s$ and $t$ are 'close' even in
case both components are 'short'.

\begin{lem}\label{lin-bcp2}
For any $\lambda \ge 1, c\ge 0$ there exists a constant $D
=D(\lambda , c,k)$ such that the following condition holds. Let
$p$, $q$ be a pair of $k$--similar $(\lambda ,
c)$--quasi--geodesics without backtracking in $\G $. Suppose that
$s$ and $t$ are connected $H_i$--components of $p$ and $q$
respectively. Then
$$\max \{ \dx (s_-, t_-),\; \dx (s_+, t_+)\} \le
D.$$
\end{lem}
\begin{proof}
Let $p_1=[p_-, s_-]$ and $q_1=[q_-, t_-]$ denote the corresponding
segments of the paths $p$ and $q$ respectively, $c^\prime $ and
$c^{\prime\prime }$ denote the connectors of length at most $1$
such that $c^\prime _-=s_-$, $c^\prime _+=t_-$, $c^{\prime\prime
}_-=s_+$, and $c^{\prime\prime }_+=t_+$ (see Fig. \ref{BCP01}).
Let $C(\lambda , c,k)$ be the constant provided by Lemma
\ref{lin-BCP}. We set
$$
D= C(\lambda , c+\lambda +1, k).
$$
Let us estimate the $X$--length of $c^\prime $.

\begin{figure}
\begin{picture}(120,50)(0,10)
\qbezier(0,40)(60,70)(120,40) \qbezier(0,30)(60,0)(120,30)
\multiput(0,30)(120,0){2}{\circle*{1}}
\multiput(0,40)(120,0){2}{\circle*{1}}

\put(3,31){$p_-$} \put(114,31){$p_+$} \put(3,39){$q_-$}
\put(114,39){$q_+$}

\put(25,20){\circle*{1}} \put(95,20){\circle*{1}}
\put(45,54){\circle*{1}} \put(75,54){\circle*{1}}

\put(59,57){$t$} \put(59,12){$s$}

\put(23,17){$s_-$} \put(96,17){$s_+$}

\put(44,34){$c^\prime$} \put(73,34){$c^{\prime\prime }$}

\put(15,20){$p_1$} \put(15,49){$q_1$}
\put(15,23.4){\vector(3,-1){1}} \put(15,46.55){\vector(3,1){1}}

\put(0,30){\line(0,1){10}} \put(120,30){\line(0,1){10}}

\thicklines \qbezier(25,20)(60,10)(95,20)
\qbezier(45,54)(60,56)(75,54) \qbezier(25,20)(50,35)(45,54)
\qbezier(75,54)(70,35)(95,20) \put(42.1,35){\vector(1,2){1}}
\put(78.2,35){\vector(-1,2){1}} \put(59,15){\vector(1,0){1}}
\put(60,55){\vector(1,0){1}}

\end{picture}
\caption{} \label{BCP01}
\end{figure}

Assume that the path $p_1c^\prime $ has a backtracking. Then there
exists an $H_i$--component $b$ of $p_1$ that is connected to
$c^\prime $. Therefore $b$ is connected to $s$ contradicting to
the assumption that $p$ is a path without backtracking. The same
arguments show that no component of $q_1$ is connected to
$c^\prime $. Further, by Lemma \ref{plus1} the path $p_1c^\prime $
is $(\lambda , c+\lambda+1)$--quasi--geodesic. Thus we can apply
Lemma \ref{lin-BCP} to $q_1$, $p_1c^\prime $, and the
$H_i$--component $c^\prime $ of $p_1c^\prime $. We obtain
$$
\dx (s_-, t_-)\le C(\lambda , c+\lambda +1, k)
$$
The bound on $\dx (s_+, t+)$ can be obtained in the same way.
\end{proof}

Taking together Proposition \ref{QG} and Lemmas \ref{lin-BCP},
\ref{lin-bcp2}, we obtain the following theorem, which is closely
related to Farb's Bounded Coset Penetration property (see the
appendix). For simplicity we change the notation and denote by
$\varepsilon (\lambda , c, k)$ the maximal constant among
$\varepsilon (\lambda ,c,k)$ from Proposition \ref{QG},  $C
(\lambda ,c,k)$ from Lemma \ref{lin-BCP},  $D (\lambda ,c,k)$ from
Lemma \ref{lin-bcp2}.

\begin{thm}\label{TBCP}
For any $\lambda \ge 1$, $c\ge 0$, $k\ge 0$, there exists a
constant $\varepsilon = \varepsilon (\lambda , c, k)$ such that
for any two $k$--similar $(\lambda , c)$--quasi--geodesics without
backtracking $p$ and $q$ in $\G $, the following conditions hold.

1) The sets of phase vertices of $p$ and $q$ are contained in the
closed $\varepsilon $--neighborhoods (with respect to the metric $
\dx $) of each other.

2) Suppose that $s$ is an $H_i$--component of $p$ such that $\dx
(s_-, s_+) > \varepsilon $; then there exists an $H_i$--component
$t$ of $q$ which is connected to $s$.

3) Suppose that $s$ and $t$ are connected $H_i$--components of $p$
and $q$ respectively. Then $$\max \{ \dx (s_-, t_-),\;  \dx (s_+,
t_+)\} \le \varepsilon .$$
\end{thm}


\section{Geodesic triangles in Cayley graphs}


Sometimes it is useful to think of $\delta $--hyperbolic spaces as
being fattened versions of trees (see \cite{Bow91}, \cite{CoDP},
and \cite{GhH}). More precisely, given any three positive numbers
$a,b,c$, we can consider the metric tree $T(a,b,c)$ that has three
vertices of valence one, one vertex of valence three, and edges of
length $a,b,$ and $c$. For convenience we extend the definition of
$T(a,b,c)$ in the obvious way to cover the cases where $a,b$ and
$c$ are allowed to be zero.

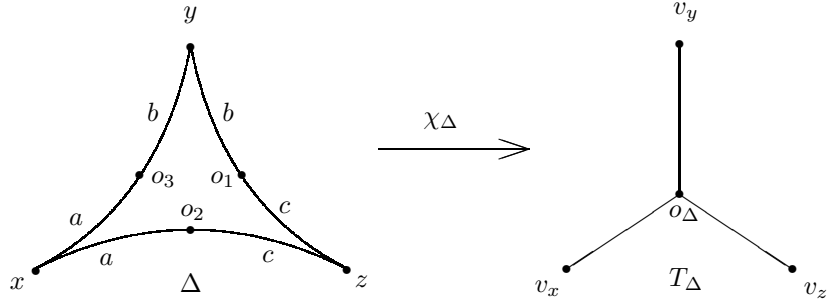
\begin{figure}
\unitlength 1mm 
\begin{picture}(108.36,42.25)(0,0)
\qbezier(28,37.65)(23.86,16.44)(7.07,7.95)
\qbezier(28,37.83)(32,16.62)(48.79,8.13)
\qbezier(7.07,7.95)(28,18.47)(48.61,8.13)

\put(53,24){\line(1,0){20}}

\put(73,24){\line(-3,-1){4}} \put(73,24){\line(-3,1){4}}

\put(93,18){\line(0,1){20}} \put(93,18){\line(3,-2){15}}
\put(93,18){\line(-3,-2){15}} \put(21.21,20.51){\circle*{1}}
\put(34.83,20.51){\circle*{1}} \put(48.8,7.96){\circle*{1}}
\put(7.43,7.78){\circle*{1}} \put(28,13.26){\circle*{1}}

\put(93,18){\circle*{1}} \put(93,38){\circle*{1}}
\put(78,8){\circle*{1}} \put(108,8){\circle*{1}}

\put(28,37.65){\circle*{1}} \put(5,6.19){\makebox(0,0)[cc]{$x$}}
\put(28,41.9){\makebox(0,0)[cc]{$y$}}
\put(50.73,6.36){\makebox(0,0)[cc]{$z$}}
\put(61.34,27.75){\makebox(0,0)[cc]{$\chi _\Delta $}}
\put(28,6.19){\makebox(0,0)[cc]{$\Delta $}}
\put(93.87,6.72){\makebox(0,0)[cc]{$T_\Delta $}}
\put(93.87,42.25){\makebox(0,0)[cc]{$v_y$}}
\put(75.48,5.48){\makebox(0,0)[cc]{$v_x$}}
\put(111.2,5.2){\makebox(0,0)[cc]{$v_z$}}
\put(23,29){\makebox(0,0)[cc]{$b$}}
\put(33,29){\makebox(0,0)[cc]{$b$}}
\put(12.8,14.67){\makebox(0,0)[cc]{$a$}}
\put(16.8,9.32){\makebox(0,0)[cc]{$a$}}
\put(38.3,9.72){\makebox(0,0)[cc]{$c$}}
\put(40.48,16.26){\makebox(0,0)[cc]{$c$}}
\put(28.3,15.56){\makebox(0,0)[cc]{$o_2$}}
\put(24.4,20){\makebox(0,0)[cc]{$o_3$}}
\put(32.3,20){\makebox(0,0)[cc]{$o_1$}}
\put(93.69,15.38){\makebox(0,0)[cc]{$o_\Delta $}}

\end{picture}

\caption{The map $\chi _\Delta $. The points $o_1,o_2,o_3$ are
mapped to $o_\Delta $}\label{Tabc}
\end{figure}

Given any three point $x,y,z$ in a metric space, the triangle
inequality tells us that there exists (unique) non--negative
numbers $a$, $b$, and $c$ such that $dist(x,y)=a+b$,
$dist(x,z)=a+c$, $dist(y,z)=b+c$. There is an isometry from $\{
x,y,z\} $ to a subset of vertices of $T(a,b,c)$ (the vertices of
valence one in the non--degenerate case); we label these vertices
$v_x, v_y, v_z$ in the obvious way (see Fig. \ref{Tabc}).

For a geodesic triangle $\Delta =\Delta (x,y,z)$ with vertices
$x,y,z$, we define $T_\Delta =T(a,b,c)$, where $a,b,c$ are chosen
as above. By $o_\Delta $ we denote the central vertex of $T_\Delta
$.  The above map $\{ x,y,z\} \to \{ v_x, v_y, v_z\} $ extends
uniquely to a map $\chi _\Delta :\Delta \to T_\Delta $ whose
restriction to each side of $\Delta $ is an isometry.

\begin{defn}\label{chi}
Let $\Delta $ be a geodesic triangle in a metric space $Y$.
Consider the map $\chi _\Delta : \Delta \to T_\Delta $ defined
above. We say that a point $u\in \Delta $ is {\it conjugate} to a
point $v\in \Delta $ if $\chi _\Delta  (u)=\chi _\Delta (v) $. The
triangle $\Delta $ is said to be {\it $\xi $--thin} if  $dist
(u,v)\le \chi $ for any two conjugate points $u,v\in \Delta $.
\end{defn}

The next lemma is well--known. It provides an equivalent
definition of relative hyperbolicity.

\begin{lem}\label{thin}
A geodesic metric space $Y$ is hyperbolic if and only if there
exists $\xi \ge 0$ such that every geodesic triangle in $Y$ is
$\xi $--thin.
\end{lem}

The main result of this section is the relative analogue of the
Rips condition for hyperbolic spaces. Namely we show that geodesic
triangles in $\Gamma (G, X\cup\mathcal H)$ are thin with respect
to the metric $dist_X$ in the following sense.

\begin{thm} \label{r0}
There exists a constant $\nu >0$ having the following property.
Let $\Delta = pqr$ be a triangle whose sides $p,q, r$ are
geodesics in $\Gamma (G, X\cup \mathcal H)$. Then for any vertex
$v$ on $p$, there exists a vertex $u$ on the union $q\cup r$ such
that
$$\dx (u,v)\le \nu .$$
\end{thm}

\begin{proof}

\begin{figure}
\begin{picture}(120,90)(-15,-20)
\put(0,0){\line(3,4){45}} \put(0,0){\line(1,0){90}}
\put(90,0){\line(-3,4){45}} \put(35,0){\line(-1,2){14}}

\multiput(9,12)(12,16){3}{\circle*{1}}
\multiput(15,0)(20,0){3}{\circle*{1}}

\multiput(18,24)(12,16){2}{\vector(-3,-4){3}}
\multiput(25,0)(20,0){2}{\vector(1,0){0}}

\put(0,0){\circle*{1}} \put(45,60){\circle*{1}}
\put(90,0){\circle*{1}}

\put(4,12){$u_1$} \put(17,28){$z$} \put(28,44){$u_2$}
\put(13,-3){$v_1$} \put(34,-3){$v$} \put(53,-3){$v_2$}
\put(-14,0){$p_-=r_+$} \put(92,0){$p_+=q_-$}
\put(39,62){$r_-=q_+$} \put(25,2){$p_1$} \put(45,2){$p_2$}
\put(18,20){$r_1$} \put(30,36){$r_2$} \put(31,10){$o$}

\put(12,33){$r$} \put(77,33){$q$} \put(44,-10){$p$}

\put(16,5){\vector(0,1){0}} \put(30,10){\vector(-1,2){2}}
\put(45,17){\vector(-1,2){2}}

\put(17,5){$o_1$}
\put(45,19){$o_2$}

\qbezier(15,0)(18.4,10)(9,12) \qbezier(55,0)(39.7,24)(33,44)
\end{picture}

\begin{picture}(120,90)(-15,-20)
\put(0,0){\line(3,4){45}} \put(0,0){\line(1,0){90}}
\put(90,0){\line(-3,4){45}} \put(35,0){\line(-1,2){14}}

\multiput(9,12)(12,16){3}{\circle*{1}}
\multiput(15,0)(20,0){3}{\circle*{1}}
\multiput(72,24)(-12,16){2}{\circle*{1}}

\multiput(18,24)(12,16){2}{\vector(-3,-4){3}}
\multiput(25,0)(20,0){2}{\vector(1,0){0}}

\put(0,0){\circle*{1}} \put(45,60){\circle*{1}}
\put(90,0){\circle*{1}}

\put(4,12){$u_1$} \put(17,28){$z$} \put(28,44){$w_1$}
\put(61,40){$w_2$} \put(73,24){$u_2$} \put(13,-3){$v_1$}
\put(34,-3){$v$} \put(53,-3){$v_2$} \put(-14,0){$p_-=r_+$}
\put(92,0){$p_+=q_-$} \put(39,62){$r_-=q_+$} \put(25,2){$p_1$}
\put(45,2){$p_2$} \put(18,20){$r_1$} \put(30,36){$r_2$}
\put(30,12){$o$}

\put(12,33){$r$} \put(77,33){$q$} \put(44,-10){$p$}

\put(16,5){\vector(0,1){0}} \put(29,12){\vector(-1,2){2}}
\put(54,12){\vector(0,1){0}} \put(48,36){\vector(1,0){0}}

\put(17,5){$o_1$}
\put(47,33){$o_3$} \put(50,13){$o_2$}

\qbezier(15,0)(18.4,10)(9,12) \qbezier(55,0)(49.5,30)(72,24)
\qbezier(33,44)(45,30.4)(60,40)

\put(69,28){\vector(-3,4){3}} \put(62,31){$q_0$}
\end{picture}
\caption{The quadrilateral and hexagonal cases in the proof of
Lemma \ref{r0}.} \label{QH}
\end{figure}

The logical scheme of the proof is similar to that of Proposition
\ref{QG}. First we need an auxiliary construction. Let $T_\Delta $
be the tree related to $\Delta $ and $\chi _\Delta : \Delta \to
T_\Delta $ be the corresponding map. Recall that $p_-=r_+,
p_+=q_-$. Let us consider vertices $v_1, v_2$ of $p$ chosen as
follows. Since $\G $ is a $\delta $--hyperbolic space for some
$\delta >0$, according to Lemma \ref{thin}, there exists a
constant $\xi $ such that $\Delta $ is $\xi $--thin. Increasing
$\xi $ if necessary we can assume $\xi \in \mathbb N$.

We chose vertices $v_1$ and $v_2$ on $p$ in the following way. If
$dist_{X\cup \mathcal H} (p_-, v)\le 6\xi $ (respectively
$dist_{X\cup \mathcal H} (v, p_+)\le 6\xi $), we put $v_1=p_- $
(respectively $v_2=p_+$). If $dist_{X\cup \mathcal H} (p_-, v)>
6\xi  $ (respectively $dist_{X\cup \mathcal H} (v, p_+)> 6\xi $),
we take $v_1$ on the segment $[p_-, v]$ of the path $p$
(respectively $v_2$ on the segment $[v, p_+]$ of the path $p$)
such that
\begin{equation} \label{rips0}
dist_{X\cup \mathcal H} (v, v_i)=6\xi
\end{equation}
for $i=1$ (respectively $i=2$).

Denote by $u_1, u_2$ the vertices on $q\cup r$ such that $\chi
_\Delta (u_i)=\chi _\Delta (v_i)$, $i=1,2$. By Lemma \ref{thin},
we have
\begin{equation} \label{rips1}
dist_{X\cup \mathcal H} (v_i, u_i)\le \xi , \;\;\; i=1,2.
\end{equation}
Reversing the roles of $p$ and $r$ if necessarily, we can assume
that $u_1\in r$. Thus there are only two cases to consider (see
Fig. \ref{QH}).

{\it Case 1} (Quadrilateral).  The vertex $u_2$ lays on $r$. In
this case we denote by $p_1, p_2$ (respectively by $r_0$) the
segments of $p$ (respectively the segment of $r$) such that
$(p_1)_-=v_1$, $(p_1)_+=v$, $(p_2)_-=v$, $(p_2)_+=v_2$
(respectively $(r_0)_-=u_2$, $(r_0)_+=u_1$). We also denote by
$o_1, o_2$ the geodesic paths in $\G $ such that $(o_i)_-=v_i$,
$(o_i)_+=u_i$, $i=1,2$.

{\it Case 2} (Hexagonal). No vertex on $r$ is conjugate to $v_2$
(therefore, $u_2$ lays on $q$, see Fig. \ref{QH}). We take a
vertex $w_1$ on the segment $[r_-, u_1]$ of the path $r$ such that
$w_1=r_-$ if $dist_{X\cup \mathcal H}(r_-, u_1)\le 12\xi $ and
\begin{equation} \label{rips11}
dist_{X\cup \mathcal H}(w_1, u_1)=12\xi
\end{equation}
if $dist_{X\cup \mathcal H}(r_-, u_1)>12\xi $. Furthermore, if
$w_1=r_-$, we set $w_2=w_1$, and if $w_1\ne r_-$, let  $w_2$ be
the vertex on $p\cup q$ such that $\chi _\Delta (w_2)=\chi _\Delta
(w_1)$. Thus we have
\begin{equation} \label{rips12}
dist_{X\cup \mathcal H}(w_1, w_2)\le \xi .
\end{equation}
It is easy to see that $w_2\in q$. Indeed, if $w_2\in p$, then
$$\dxh (v_1, w_2)=\dxh (u_1, w_1)=12\xi \ge \dxh (v_1,v_2),$$ i.e.,
the point $v_2$ belongs to the segment $[v_1, w_2]$ of $p$. Since
$r$ contains conjugate points for both $v_1$ and $w_2$, there is a
conjugate point for $v_2$ on $r$ that contradicts to our
assumption. Similar arguments show that $w_2$ belongs to the
segment $[u_2,r_-]$ of the path $q$.

We denote by $p_1, p_2$ (respectively by $r_0$) the segments of
$p$ (respectively the segment of $r$) such that $(p_1)_-=v_1$,
$(p_1)_+=v$, $(p_2)_-=v$, $(p_2)_+=v_2$ (respectively
$(r_0)_-=w_1$, $(r_0)_+=u_1$), and denote by $q_0$ the segment of
$q$ such that $(q_0)_-=u_2,$ $(q_0)_+=w_2$. Also let $o_1, o_2,
o_3$ denote the geodesic paths in $\G $ such that $(o_i)_-=v_i$,
$(o_i)_+=u_i$, $i=1,2$, $(o_3)_-=w_1$, $(o_3)_+=w_2$.

We deal with the second case in details, the first one is only
slightly different. Let $V$ denote the set of all vertices $z$ on
$r_0\cup q_0$ such that
\begin{equation}
dist _{X\cup \mathcal H} (v,z)=\min\limits_{z_0\in r_0\cup q_0}
dist _{X\cup \mathcal H} (v,z_0),
\end{equation}
where the minimum is taken among all vertices of $r_0\cup q_0$. In
particular, we have
\begin{equation}
dist _{X\cup \mathcal H} (v,z)\le \xi
\end{equation}
since the conjugate point for $z$ is a vertex on $r_0\cup q_0$.

To every $z \in V$, we associate the set $O(z)$ of all a geodesic
paths $o$ such that $o_-=v, o_+=z$. These paths cut the cycle
$p_1p_2o_2q_0o_3^{-1}r_0(o_1)^{-1}$ into two parts denoted by
$c_1$ and $c_2$. More precisely, if $z\in r_0$, we set
$$ c_1=or_1(o_1)^{-1}p_1,
$$ and $$c_2=o(r_2)^{-1}o_3(q_0)^{-1}(o_2)^{-1}(p_2)^{-1}, $$
where $r_1$ (respectively $r_2$) is the segment $[z, u_1]$
(respectively $[w_1, z]$) of the path $r_0$. If $z\in q_0$, we set
$$ c_1=oq_2o_3^{-1}r_0o_1^{-1}p_1
$$ and $$c_2=p_2o_2q_1o^{-1}, $$
where $q_1$ (respectively $q_2$) is the segment $[u_2,z]$
(respectively $[z, w_2]$) of the path $q_0$. We emphasize that
$c_1$, $c_2$ depend on the choice of $z\in V$ and $o\in O(z)$.

\begin{lem}\label{rips39}
Let $z\in V$, $o\in O(z)$, and let $s$ be an $H_i$--component of
the path $o$ for a certain $i$. Then for any $j=1,2,3$ there exist
no $H_i$--components of $o_j$ connected to $s$.
\end{lem}

\begin{figure}
\begin{picture}(120,90)(-15,-20)
\put(0,0){\line(3,4){45}} \put(0,0){\line(1,0){90}}
\put(90,0){\line(-3,4){45}} \put(35,0){\line(-1,2){14}}

\multiput(9,12)(12,16){3}{\circle*{1}}
\multiput(15,0)(20,0){3}{\circle*{1}}
\multiput(72,24)(-12,16){2}{\circle*{1}}

\multiput(18,24)(12,16){2}{\vector(-3,-4){3}}
\multiput(25,0)(20,0){2}{\vector(1,0){0}}

\put(0,0){\circle*{1}} \put(45,60){\circle*{1}}
\put(90,0){\circle*{1}}

\put(4,12){$u_1$} \put(17,28){$z$} \put(28,44){$w_1$}
\put(61,40){$w_2$} \put(73,24){$u_2$} \put(13,-3){$v_1$}
\put(34,-3){$v$} \put(53,-3){$v_2$} \put(-14,0){$p_-=r_+$}
\put(92,0){$p_+=q_-$} \put(39,62){$r_-=q_+$} \put(25,2){$p_1$}
\put(45,2){$p_2$} \put(11,20){$r_1$} \put(30,36){$r_2$}

\put(12,33){$r$} \put(77,33){$q$} \put(44,-10){$p$}

\put(16,5){\vector(0,1){0}} \put(54,12){\vector(0,1){0}}
\put(48,36){\vector(1,0){0}}

\put(17,5){$o_1$}

\put(47,33){$o_3$}

\qbezier(15,0)(18.4,10)(9,12) \qbezier(55,0)(49.5,30)(72,24)
\qbezier(33,44)(45,30.4)(60,40)

\put(69,28){\vector(-3,4){3}} \put(62,31){$q_0$}

\thicklines \put(29,12){\vector(-1,2){2}}
\put(30,10){\line(-1,2){5}} \put(30,10){\circle*{1}}
\put(25,20){\circle*{1}} \qbezier(54.2,5)(53.3,12)(54.7,17)
\put(54.2,5){\circle*{1}} \put(54.7,17){\circle*{1}}
\qbezier(54.2,5)(42,12)(30,10) \put(54,12){\vector(0,1){0}}
\put(29.3,15){$s$} \put(25.5,10){$s_-$} \put(20.5,20){$s_+$}
\put(55.5,5){$t_-$} \put(56,16){$t_+$} \put(51.5,11.5){$t$}
\end{picture}

\begin{picture}(120,90)(-15,-20)
\put(0,0){\line(3,4){45}} \put(0,0){\line(1,0){90}}
\put(90,0){\line(-3,4){45}} \put(35,0){\line(-1,2){14}}

\multiput(9,12)(12,16){3}{\circle*{1}}
\multiput(15,0)(20,0){3}{\circle*{1}}
\multiput(72,24)(-12,16){2}{\circle*{1}}

\multiput(18,24)(12,16){2}{\vector(-3,-4){3}}
\multiput(25,0)(20,0){2}{\vector(1,0){0}}

\put(0,0){\circle*{1}} \put(45,60){\circle*{1}}
\put(90,0){\circle*{1}}

\put(4,12){$u_1$} \put(17,28){$z$} \put(28,44){$w_1$}
\put(61,40){$w_2$} \put(73,24){$u_2$} \put(13,-3){$v_1$}
\put(34,-3){$v$} \put(53,-3){$v_2$} \put(-14,0){$p_-=r_+$}
\put(92,0){$p_+=q_-$} \put(39,62){$r_-=q_+$} \put(25,2){$p_1$}
\put(45,2){$p_2$} \put(11,20){$r_1$} \put(30,36){$r_2$}

\put(12,33){$r$} \put(77,33){$q$} \put(44,-10){$p$}

\put(16,5){\vector(0,1){0}} \put(54,12){\vector(0,1){0}}
\put(54,12){\vector(0,1){0}}

\put(17,5){$o_1$}

\qbezier(15,0)(18.4,10)(9,12) \qbezier(55,0)(49.5,30)(72,24)
\qbezier(33,44)(45,30.4)(60,40)

\put(69,28){\vector(-3,4){3}} \put(62,31){$q_0$}

\thicklines \put(29,12){\vector(-1,2){2}}
\put(30,10){\line(-1,2){5}} \put(30,10){\circle*{1}}
\put(25,20){\circle*{1}}

\qbezier(40,38.1)(46,35)(53,37)

\put(40,38.1){\circle*{1}} \put(53,37){\circle*{1}}
\qbezier(40,38.1)(40,32)(25,20) \put(29.3,15){$s$}

\put(49.5,36.2){\vector(1,0){0}} \put(50,13){$o_2$}

\put(25.5,10){$s_-$} \put(20.5,20){$s_+$}

\put(40,39){$t_-$} \put(53,39){$t_+$} \put(49,32.5){$t$}
\end{picture}

\caption{Two cases in the proof of Lemma \ref{rips39}.}
\label{j12}
\end{figure}

\begin{proof}
Without loss of generality we may assume that $z\in r_0$; the case
$z\in q_0$ can be treated in the similar way and we leave it to
the reader.

The proof of the lemma in the cases $j=1$ and $j=2$ almost
coincide with the proof of the Lemma \ref{Rips39}. Indeed assume
that $s$ is connected to an $H_i$--component $t$ of $o_j$ for
$j=1$ or $j=2$. This means that
$$
\dxh (s_-, t_-)\le 1.
$$
Thus we have
$$
\begin{array}{rl}
dist _{X\cup \mathcal H}(v, v_j)& \le dist _{X\cup \mathcal
H}(v,s_-)+dist _{X\cup \mathcal H}(s_-, t_-) +dist _{X\cup
\mathcal H}(t_-, v_j)\\ &
\\& \le (\xi -1)+1+(\xi -1)<2\xi .
\end{array}
$$
By our choice of $v_j$, this implies $v_j=p_-$ in case $j=1$ and
$v_j=p_+$ in case $j=2.$ Therefore, $v_j$ coincides with $u_j$ as
$$
dist _{X\cup \mathcal H} (v_1, p_-)=dist _{X\cup \mathcal H}(u_1,
p_-)
$$
and
$$
dist _{X\cup \mathcal H} (v_2, p_+)=dist _{X\cup \mathcal H}(u_2,
p+).
$$
Thus $o_j$ is trivial. A contradiction.

Suppose that $s$ is connected to an $H_i$--component $t$ of $o_3$.
In particular this means that $o_3$ is non--empty. By the choice
of $w_1$ and $w_2$ this implies the equality $\dxh (u_1,
w_1)=12\xi $. Note that
$$
dist _{X\cup \mathcal H}(u_1, z) \le dist _{X\cup \mathcal
H}(u_1,v_1)+dist _{X\cup \mathcal H}(v_1, v) +dist _{X\cup
\mathcal H}(v,z) \le 8\xi .
$$
On one hand, we obtain
\begin{equation} \label{r39}
dist _{X\cup \mathcal H}(z, w_1) \ge dist _{X\cup \mathcal H}(u_1,
w_1)- dist_{X\cup \mathcal H}(u_1, z) \ge 12\xi -8\xi =4\xi .
\end{equation}
On the other hand, we have
$$
\begin{array}{rl}
dist _{X\cup \mathcal H}(z, w_1)& \le dist _{X\cup \mathcal H}(z,
s_+)+ dist_{X\cup \mathcal H}(s_+, t_-)+ dist_{X\cup\mathcal
H}(t_-, w_1)\\ &\\ & \le (\xi -1) + 1+ (\xi -1) <2\xi
\end{array}
$$
that contradicts to (\ref{r39}). The lemma is proved.
\end{proof}

As in the previous section, given $z\in V$ and $o\in O(z)$, an
$H_i$--component $s$ of the path $o$ is called an {\it ending
component} if $s$ contains $z$; for otherwise $s$ is called a {\it
non--ending} component. The proof of the next lemma is completely
analogous to the proof of Lemma \ref{non-end}. We leave details to
the reader.

\begin{lem}\label{1non-end}
Let $z\in V$, $o\in O(z)$, and $s$ be a non--ending
$H_i$--component of the path $o$ for a certain $i$. For any
$j=1,2$, if $s$ is not an isolated $H_i$--component of $c_j$, then
there exists an $H_i$--component $t$ of $p_j$ that is connected to
$s$.
\end{lem}

Arguing as in the proof of Corollary \ref{non-end1}, we
immediately obtain the following.

\begin{cor}\label{1non-end1}
Let $z\in V$, $o\in O(z)$. Then for any $i=1, \ldots , m$, every
non--ending  $H_i$--component of the path $o$ is an isolated
$H_i$--component of at least one of the cycles $c_1, c_2$.
\end{cor}

The next lemma is the analogue of Lemma \ref{end}.

\begin{lem}\label{1end}
For any $j=1,2$, the following assertion is true. Suppose that for
any $z\in V$ and any $o\in O(z)$, $o$ has the ending component,
which is not isolated in $c_j$; then for any $z\in V$ and any
$o\in O(z)$, the ending component of $o$ is connected to an
$H_i$--component of $p_j$ for corresponding $i$.
\end{lem}

\begin{proof}
For definiteness assume $j=1$. For any $z\in V$, we define a
non-negative integer number $\pi (z)$ as follows
$$
\pi (z)=\left\{
\begin{array}{l}
\dxh (u_1, z), \; {\rm if }\; z\in r_0, \\ \\ \dxh (u_1, w_1)
+\dxh (w_2, z), \; {\rm if }\; z\in q_0.
\end{array}
\right.
$$
The reader will have no difficulties in proving this lemma in the
same way as Lemma \ref{end}. The only difference is that we have
to proceed by induction on $\pi (z)$.
\end{proof}

Obviously Lemma \ref{1end} yields

\begin{cor}\label{1end1}
There exists a vertex $z\in V$ and a path $o\in O(z)$ such that
either the last edge of $o$ is labelled by a letter from $X$ or
the ending component of $o$ is isolated in at least one of the
cycles $c_1, c_2$.
\end{cor}

Let us return to the proof of Theorem \ref{r0}.

By Corollaries \ref{1non-end1} and \ref{1end1}, there exists a
vertex $z\in V$ and a path $o\in O(z)$ such that every component
of $o$ is isolated in at least one of the cycles $c_1, c_2$.
Therefore, we have
\begin{equation} \label{r1}
\dx (u,z)\le l(o) ML\max \{ l(c_1), \; l(c_2)\} .
\end{equation}
Since
\begin{equation} \label{r2}
l(o)\le \xi
\end{equation}
by our choice of $z$, it remains to estimate the lengths of $c_1$
and $c_2$. Evidently we have
$$
l(q_0)\le l(o_3) + l(r_0) + l(o_1) + l(p_1) +l(p_2) +l(o_2) \le
27\xi .
$$
Therefore,
\begin{equation}\label{r3}
l(c_i)< l(o) +l(p_1)+l(p_2)+l(o_2)+l(q_0)+l(o_3) +l(r_0) +l(o_1)
\le 55\xi .
\end{equation}
In view of (\ref{r1}), (\ref{r2}), and (\ref{r3}), to complete the
proof it suffices to set $\nu = 55ML \xi ^2. $
\end{proof}

\begin{rem}
For metric spaces, the Rips condition can be regarded as the
definition of hyperbolicity. We note that the fulfilment of
Theorem \ref{r0} for a given group and a collection of subgroups
does not imply the relative hyperbolicity. Indeed for the pair
$G\cong \mathbb Z$ and $H\cong 2\mathbb Z$ with the natural
embedding $H\to G$, the statement of Theorem \ref{rips0} obviously
holds. However, $G$ is not hyperbolic relative to $H$. Moreover,
the corresponding relative Dehn function is not well--defined as
follows from Proposition \ref{malnorm}.
\end{rem}

By drawing the diagonal, we obtain the following corollary of
Theorem \ref{rips0}. It will be used in the next chapter to study
the root problem for relatively hyperbolic groups.

\begin{cor}\label{quad}
Let $p_1p_2p_3p_4$ be a geodesic quadrangle in $\G $. Then for any
vertex $u\in p_1$ there is a vertex $v\in p_2\cup p_3\cup p_4$
such that $$\dx (u,v)\le 2\nu .$$
\end{cor}


\section{Symmetric geodesics}\label{SecSymGeod}


The proofs of a number of theorems about ordinary hyperbolic
groups (in particular, the solution of the conjugacy problem) are
based on the following well--known property of hyperbolic metric
spaces \cite{BriH,Ep-etal}. Let $p, q$ be two geodesics such that
the distances between $p_-, q_-$ and between $p_+, q_+$ are
'small', say less than or equal to $1$, and $u,v$ are two points
on $p$ and $q$ respectively such that $dist (p_-, u)=dist (q_-,
v)$. Then the distance between $u$ and $v$ is not greater than
$k$, where $k$ is independent of $p$, $q$, $u$, and $v$. The
straightforward relative analogue of this property in the spirit
of Proposition \ref{QG} and Theorem \ref{r0} can be stated as
follows.

\begin{conj} \label{false}
For any $k\ge 0$, there exists a constant $\varkappa =\varkappa
(k)$ such that the following condition holds. Let $p,q$ be a pair
of $k$--similar geodesics in $\G $. Then for any two vertices
$u\in p$ and $v\in q$ such that $dist (p_-, u)=dist (q_-, v)$, we
have
\begin{equation} \label{false0}
\dx(u,v)\le \varkappa .
\end{equation}
\end{conj}

Unfortunately in general Conjecture \ref{false} is false. Indeed
consider the free group $G=\langle x,y\rangle $ and the subgroup
$H=\langle x\rangle $, which is a free factor of $G$. Evidently
$G$ is hyperbolic relative to $H$. Let us also consider two
geodesics, denoted by $p$ and $q$, in the corresponding relative
Cayley graph such that $$\phi (p)\equiv yx^n,\;\;\;\;\; \phi
(q)\equiv x^n,$$ and $$p_-=y^{-1}, \;\;\;\;\; q_-=1.$$ It is clear
that $p$ and $q$ are $1$--similar. Take the vertices $u\in p$ and
$v\in q$ such that $\dxh (p_-, u)=\dxh (q_-, v)=1$. Thus $u=1$,
$v=x^n$, and the length $u^{-1}v$ with respect to the generating
set $\{ x,y\} $ of $G$ equals $n$. As $n$ can be taken arbitrary
large, this obviously violates (\ref{false0}).

However a certain refined version of the above mentioned
conjecture can be proved for relatively hyperbolic groups. The
results of this section will be used in the next chapter in order
to study cyclic subgroups of relatively hyperbolic groups. We
begin with definitions.

\begin{defn}
Let $p,q$ be two paths in $\G $. We say that the pair $(p,q)$ is
{\it symmetric} if the labels of $p$ and $q$ coincide, i.e., $\phi
(p)\equiv \phi (q)$. To each such a pair we associate the pair of
elements $g_1=(p_-)^{-1}q_-$, and $g_2=(p_+)^{-1}q_+$ of $G$,
called the {\it characteristic elements} of $(p,q)$.
\end{defn}

It is easy to see that two elements $g_1, g_2$ are conjugate by
some element $t\in G$, i.e.,
$$
g_2=t^{-1}g_1t
$$
if and only if there exists a symmetric pair of geodesics $(p,q)$
in $\G $ such that $(g_1, g_2)$ is the characteristic pair of
$(p,q)$ and $\overline {\phi (p)}=\overline {\phi (p)}=t$.

\begin{defn}
Let $(p,q)$ be a pair of symmetric paths in $\G $. We say that two
vertices $u\in p$, and $v\in q$ are {\it synchronous}, if $l([p_-,
u])=l([q_-, v])$, where $[p_-, u]$ and $[q_-, v]$ are segments of
$p$ and $q$ respectively. Similarly, if $a$ and $b $ are
$H_i$--components of $p$ and $q$ respectively for a certain $i$,
we say that $a$ and $b$ are {\it synchronous components} when the
vertices $a_-$ and $b_-$ are synchronous.
\end{defn}

Now we are ready to state the main result of this section.

\begin{thm}\label{sc}
For any $k\ge 0 N$, there exists a constant $\varkappa =\varkappa
(k)$ with the following property. Let $(p,q)$ be a symmetric pair
of $k$--similar geodesics in $\G $, $u,v$ a pair of synchronous
vertices on $p$ and $q$ respectively. Then
$$
\dx (u,v)\le \varkappa .
$$
\end{thm}

The proof provided below can be slightly simplified by using
references to some technical lemmas from \cite{Bum}. However for
convenience of the reader we give a complete proof in the spirit
of our paper.

As usual, we divide the proof into a sequence of lemmas. The
following is a particular case of a well--known property of
hyperbolic spaces (see for example, \cite{BriH}).

\begin{lem}\label{sk1}
For any $k\ge 0$, there exists a constant $E=E(k)$ such that the
following condition holds. Let $p$, $q$ be a pair of $k$--similar
geodesics in $\G $, and $u$, $v$ synchronous vertices on $p$ and
$q$ respectively. Then $\dxh (u,v)\le E. $
\end{lem}

We note that the constant $E$ can be chosen effectively for a
given $k$. In fact, one can take $E(k)=E^\prime k$ for an
appropriate constant $E^\prime $, which depends only on the
hyperbolicity constant of $\G $, not on $k$.

\begin{defn}\label{minpair}
We say that a symmetric pair of geodesics $(p,q)$ in $\G $ is {\it
minimal}, if for any other symmetric pair of geodesics $(p^\prime
, q^\prime )$ in $\G $ having the same characteristic elements,
the inequality $l(p)\le l(p^\prime )$ holds.
\end{defn}

\begin{lem}\label{CP1}
Let $(p,q)$ be a minimal pair of symmetric geodesics in $\G $.

1) Suppose that, for some $i$, two $H_i$--components $a$ and $b$
of $p$ and $q$ respectively are connected. Then $a$ and $b$ are
synchronous.

2) Let $u_1,v_1$ and $u_2,v_2$ be two pairs of synchronous
vertices of $p$ and $q$ respectively. Then $(u_1)^{-1}v_1\ne
(u_2)^{-1}v_2$.
\end{lem}

\begin{proof}

\begin{figure}
\begin{picture}(120,122)(0,-7)

\put(60,30){\oval(100,60)} \put(70,30){\oval(20,20)}
\put(10,30){\line(1,0){50}}

\put(15,30){\vector(1,0){1}} \put(25,30){\vector(1,0){1}}
\put(45,30){\vector(1,0){1}} \put(55,30){\vector(1,0){1}}
\put(15,32){$q_1$} \put(25,32){$a$} \put(45,26){$b$}
\put(55,27){$p_2$}

\put(100,30){\vector(0,1){1}} \put(97,30){$e$}
\put(80,30){\vector(0,1){1}} \put(82,30){$r_2$}
\put(110,30){\vector(0,1){1}} \put(112,30){$r_1$} \put(30,43){$\Xi
$}

\put(60,30){\oval(80,40)[b]} \put(75,30){\oval(50,40)[t]}

\multiput(10,30)(10,0){6}{\circle*{1}}

\put(20,80){\line(1,0){80}} \put(20,110){\line(1,0){80}}
\qbezier(20,110)(25,95)(20,80) \qbezier(100,110)(95,95)(100,80)
\multiput(70,110)(10,0){2}{\circle*{1}}
\multiput(40,80)(10,0){2}{\circle*{1}}

\put(40,80){\line(4,3){40}} \put(60,95){\vector(4,3){1}}
\put(57,95){$e$}

\put(75,110){\vector(1,0){1}} \put(90,110){\vector(1,0){1}}
\put(45,110){\vector(1,0){1}}

\put(75,111){$a$} \put(90,107){$p_2$} \put(45,107){$p_1$}

\put(45,80){\vector(1,0){1}} \put(30,80){\vector(1,0){1}}
\put(70,80){\vector(1,0){1}}

\put(45,76){$b$} \put(30,77){$q_1$} \put(70,77){$q_2$}

\put(22.6,95){\vector(0,-1){1}} \put(97.6,95){\vector(0,-1){1}}
\put(24,95){$r_1$} \put(93,95){$r_2$}

\end{picture}
\caption{Gluing and cutting diagrams in the proof of Lemma
\ref{CP1}.} \label{ann}
\end{figure}

1) Suppose that $a$ and $b$ are not synchronous. Let
$$ p=p_1ap_2, $$
$$ q=q_1bq_2. $$
Assume for definiteness  that
\begin{equation}\label{cp3}
l(q_1)<l(p_1).
\end{equation}
Since $a$ and $b$ are connected, there is an edge $e$ in $\G $
such that $\phi (e) \in \widetilde H_i$ and $$ e_-=b_-, e_+=a_+.$$
Also denote by $r_1$, $r_2$ some paths in $\G $ such that
$(r_1)_-=p_-$, $(r_1)_+=q_-$, $(r_2)_-=p_+$, $(r_2)_+=q_+$ (see
Fig. \ref{ann}). Consider the cycles
$$
c_1=e^{-1}q_1^{-1}r_1^{-1}p_1a
$$
and
$$
c_2= ep_2r_2q_2^{-1}b^{-1}.
$$

\begin{figure}
\unitlength 1.11mm 

\begin{picture}(91.54,40)(-8,-3)
\qbezier(3.01,2.83)(45.25,11.31)(91.04,2.83)
\qbezier(2.9,32.4)(45.15,23.91)(90.93,32.4)
\put(3.01,2.83){\circle*{.9}} \put(2.9,32.4){\circle*{.9}}
\put(28.97,6.51){\circle*{.9}} \put(28.86,28.72){\circle*{.9}}
\put(60.61,6.72){\circle*{.9}} \put(60.5,28.51){\circle*{.9}}
\put(90.98,2.93){\circle*{.9}} \put(90.88,32.29){\circle*{.9}}
\put(2.89,18){\vector(0,-1){.07}}\put(2.84,32.39){\line(0,-1){29.446}}
\put(28.91,18){\vector(0,-1){.07}}\put(28.91,28.71){\line(0,-1){22.298}}
\put(60.44,18){\vector(0,-1){.07}}\put(60.44,28.39){\line(0,-1){21.667}}
\put(90.92,18){\vector(0,-1){.07}}\put(90.92,32.39){\line(0,-1){29.341}}
\put(13.77,30.5){\vector(4,-1){.07}}
\put(13.03,4.6){\vector(4,1){.07}}
\put(44.09,7.06){\vector(1,0){.07}}
\put(43.52,28.15){\vector(1,0){.07}}
\put(76.94,30.24){\vector(4,1){.07}}
\put(76.94,5){\vector(4,-1){.07}}
\put(13,33.5){\makebox(0,0)[cc]{$t_1\, (f_1)$}}
\put(42.57,30.5){\makebox(0,0)[cc]{$t_2$}}
\put(76,32.7){\makebox(0,0)[cc]{$t_3\, (f_3)$}}
\put(49.4,33.02){\makebox(0,0)[cc]{$p$}}
\put(48.88,3.36){\makebox(0,0)[cc]{$q$}}
\put(28.59,31.02){\makebox(0,0)[cc]{$u_1$}}
\put(60.44,31){\makebox(0,0)[cc]{$u_2$}}
\put(28.8,4.2){\makebox(0,0)[cc]{$v_1$}}
\put(60.54,4.2){\makebox(0,0)[cc]{$v_2$}}
\put(7,18.45){\makebox(0,0)[cc]{$(g_1)$}}
\put(11.98,1.5){\makebox(0,0)[cc]{$(f_1)$}}
\put(76,1.5){\makebox(0,0)[cc]{$(f_3)$}}
\put(25.6,19){\makebox(0,0)[cc]{$(w)$}}
\put(64.12,19){\makebox(0,0)[cc]{$(w)$}}
\put(87,18){\makebox(0,0)[cc]{$(g_2)$}}
\end{picture}
\caption{The elements represented by the labels of the
corresponding paths are written in brackets.} \label{ann-1}
\end{figure}
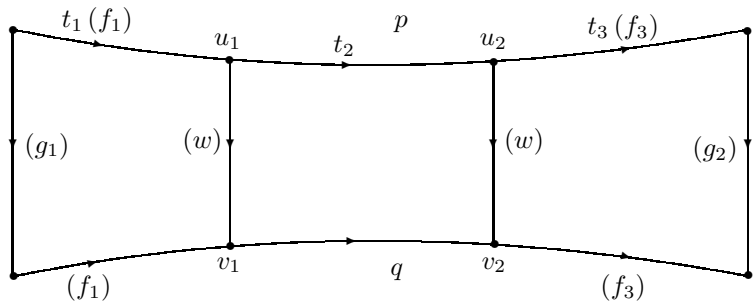

By $\Delta _1$, $\Delta _2$ we denote certain van Kampen diagrams
over (\ref{Gfg}) whose boundaries have the same labels as $c_1$
and $c_2$ respectively. For simplicity, we will identify $\partial
\Delta _j$ with $c_j$. Gluing $\Delta _1$ and $\Delta _2$ together
along subpaths $e^{-1}$ and $e$ of their boundaries, we obtain a
new diagram $\Delta $ over (\ref{Gfg}) with the boundary
$r_1^{-1}pr_2q^{-1}$. Further by gluing together the subpaths $p$
and $q$ (which have the same labels), we transform $\Delta $ into
an annular diagram $\Xi $ (Fig. \ref{ann}). Finally, we cut $\Xi $
along the image of the path $q_1ep_2$ in $\Xi $ and obtain a new
diagram $\Sigma $ over (\ref{Gfg}) with the boundary label
$$
\phi (\Sigma )\equiv \phi (r_1)^{-1} U\phi (r_2)U^{-1},
$$
where
$$
U\equiv \phi (q_1ep_2).
$$
Thus
\begin{equation}\label{cp4}
g_2=\overline{U} ^{-1}g_1\overline{U}
\end{equation}
in the group $G$, where $g_1=\overline{\phi (r_1)}$,
$g_2=\overline{\phi (r_2)}$ are the characteristic elements of
$(p,q)$. The equality (\ref{cp4}) leads to the symmetric pair
$(p^\prime, q^\prime )$ of geodesics in $\G $ such that $(g_1,
g_2)$ is the characteristic pair of $(p^\prime , q^\prime )$ and
$\phi (p^\prime )\equiv \phi (q^\prime )\equiv U$. According to
(\ref{cp3}), we have $$ l(p^\prime )= \| U\| =l(q_1) +1
+l(p_2)<l(p_1) +1 +l(p_2)=l(p)$$ that contradicts to the
assumption that the pair $(p,q)$ is minimal.

2) Let $(u_1)^{-1}v_1= (u_2)^{-1}v_2=w$ and $p=t_1t_2t_3$, where
$t_1=[p_-,u_1]$, $t_2= [u_1, u_2]$, $t_3=[u_2, p_+]$. Denote by
$f_1, f_3$ the elements represented by labels of $t_1, t_3$
respectively (see Fig. \ref{ann-1}). Then we have the following
equalities in the group $G$: $g_2=f_3^{-1}wf_3$,
$w=f_1^{-1}g_1f_1$, where $g_1,g_2$ are the characteristic
elements of $(p,q)$. Thus $g_1, g_2$ are conjugate by the element
$f_1f_3$ whose length is smaller than $l(p)$. We get a
contradiction with the minimality of $(p,q)$ again.
\end{proof}

It is worth to note that the condition (\ref{sk2-0}) in the
following lemma is weaker than $k$--similarity of $p$ and $q$.
(This is important for our goals.)

\begin{lem}\label{sk2}
Let $(p,q)$ be a minimal pair of symmetric geodesics in $\G $ such
that
\begin{equation} \label{sk2-0}
\max \{ \dxh (p_-,q_-),\; \dxh (p_+, q_+)\} \le k
\end{equation}
and let $u$ and $v$ be synchronous vertices on $p$ and $q$
respectively such that
$$ \min \{ \dxh (p_-, u), \; \dxh (u, p_+)\} \ge 2E, $$ where $E=E(k)$
is the constant provided by Lemma \ref{sk1}. Then $$ \dx (u,v)\le
6MLE^2k.$$
\end{lem}

\begin{proof}
We repeat the trick used in the proof of Proposition \ref{QG}. By
the conditions of our lemma, there exist vertices $u_1, u_2\in p$
and $v_1,v_2\in q$ such that $u_1\ne u_2$, $v_1\ne v_2$ and
\begin{equation}\label{scE1}
\dxh (u_i, u)=\dxh (v_i,v)=2E.
\end{equation}
For definiteness we assume that going along $p$ (respectively
along $q$) we first meet $u_1$ (respectively $v_1$) and then $u_2$
(respectively $v_2$).

Let $o_i$ be a geodesic path in $\G $ such that $(o_i)_-=u_i$,
$(o_i)_+=v_i$ for $i=1,2$, and $o$ be a geodesic path in $\G $
such that $o_-=u$, $o_+=v$. We consider the cycles
$$
c_1=[v_1,v]o^{-1}[u_1,u]^{-1}o_1
$$
and
$$
c_2= [v,v_2]o_2^{-1}[u,u_2]^{-1}o,
$$
where $[u_1, u]$ and $[u, u_2]$ (respectively $[v_1,v]$ and $[v,
v_2]$) are the segments of $p$ (respectively $q$).

Note that every component of $o$ is an isolated component in at
least one of the cycles $c_1,c_2$. Indeed, the same arguments as
in the proof of Lemma \ref{Rips39} together with the equalities
(\ref{scE1}) show that no component of $o_i$ is connected to a
component of $o$ for $i=1,2$. Suppose that a component $s$ of the
path $o$ is connected to a component $r_1$ in $c_1$ and $r_2$ in
$c_2$. Since $p$ and $q$ are geodesics, $r_1$ and $r_2$ can not
belong simultaneously to $p$ or $q$. Hence we can assume that
$r_1\in p$ and $r_2\in q$ (see Fig. \ref{sk2-fig}). By Lemma
\ref{CP1}, $r_1$ and $r_2$ are synchronous components. However
this can not happen if $r_1\in c_1$ and $r_2\in c_2$.

\begin{figure}
\unitlength 1mm 
\begin{picture}(103.06,37.31)(0,0)
\qbezier(6.01,6.89)(59.31,12.9)(103.06,6.89)
\qbezier(6.01,36.94)(59.31,30.93)(103.06,36.94)
\put(38.79,9.46){\circle*{1}} \put(38.79,34.37){\circle*{1}}
\put(58.13,9.99){\circle*{1}} \put(58.13,33.85){\circle*{1}}
\put(76,9.57){\circle*{1}} \put(76,34.27){\circle*{1}}
\put(29.54,8.9){\vector(1,0){.07}}\put(29.43,8.93){\line(1,0){.105}}
\put(29.01,34.9){\vector(1,0){.07}}\put(28.17,34.9){\line(1,0){.841}}
\put(58.02,33.95){\line(0,-1){23.86}}
\qbezier(38.58,34.48)(41.26,21.55)(38.68,9.46)
\qbezier(76,34.27)(73.74,24.7)(75.89,9.67) \thicklines
\qbezier(44.04,34.06)(47.35,33.95)(52.98,33.85)
\put(43.83,34.16){\circle*{1}} \put(57.92,25.96){\circle*{1}}
\put(58.02,16.08){\circle*{1}} \put(63.91,9.78){\circle*{1}}
\put(70.74,9.57){\circle*{1}} \put(52.98,34.06){\circle*{1}}
\qbezier(63.8,9.78)(67,9.77)(70.74,9.57)
\put(57.92,26.07){\line(0,-1){9.881}}
\qbezier(52.98,34.06)(52.61,28.33)(57.92,26.17)
\qbezier(57.92,16.19)(62.7,14.98)(63.91,9.78)
\put(48.56,33.95){\vector(1,0){.07}}
\put(68,9.78){\vector(1,0){.07}}
\put(57.92,21.23){\vector(0,1){.07}} \thinlines
\put(40,21.02){\vector(0,1){.07}}\put(39.94,20.29){\line(0,1){.736}}
\put(74.95,21.97){\vector(0,1){.07}}\put(74.84,21.23){\line(0,1){.736}}
\put(28.28,6){\makebox(0,0)[cc]{$p$}}
\put(27.64,37.31){\makebox(0,0)[cc]{$q$}}
\put(38.26,36.89){\makebox(0,0)[cc]{$v_1$}}
\put(47.72,36.16){\makebox(0,0)[cc]{$r_1$}}
\put(57.92,36.68){\makebox(0,0)[cc]{$v$}}
\put(75.89,36.89){\makebox(0,0)[cc]{$v_2$}}
\put(77.7,21.02){\makebox(0,0)[cc]{$o_2$}}
\put(37.7,20.18){\makebox(0,0)[cc]{$o_1$}}
\put(55.81,20.18){\makebox(0,0)[cc]{$s$}}
\put(38.37,6.5){\makebox(0,0)[cc]{$u_1$}}
\put(58.02,7.38){\makebox(0,0)[cc]{$u$}}
\put(76,6.65){\makebox(0,0)[cc]{$u_2$}}
\put(66.12,7){\makebox(0,0)[cc]{$r_2$}}
\end{picture}

  \caption{}\label{sk2-fig}
\end{figure}

Thus every component of $o$ is an isolated component in at least
one of the cycles $c_1,c_2$ and has $X$--length at most
$ML\max\limits_{i=1,2} l(c_i)$ by Lemma \ref{31}. By Lemma
\ref{sk1}, $l(o_i)\le E$ for $i=1,2$ and $l(o)\le E$. Therefore,
$l(c_i)\le 6E$ for $i=1,2$. This implies
$$
\dx (u,v)\le l(o)ML\max\limits_{i=1,2} l(c_i)\le 6MLE^2.
$$
\end{proof}

\begin{cor}\label{sk3}
For any $k\ge 0$ there exists a constant $\rho =\rho (k)$ such
that for any two conjugate elements $f,g$ of $G$ of relative
lengths $$\max\{ |f|_{X\cup \mathcal H}, \; |g|_{X\cup \mathcal H}
\} \le k,$$ there exists an element $t\in G$ such that $f^t=g$ and
$$|t|_{X\cup \mathcal H}\le \rho .$$
\end{cor}

\begin{proof}
Let $$\rho =(card\, X)^{6MLE^2k}+1+4E, $$ where $E=E(k)$ is the
constant from Lemma \ref{sk1}. Let $(p,q)$ be the minimal pair of
symmetric geodesics in $\G $ with the characteristic elements
$f,g$. Set $t=\overline{\phi (p)}$. If $l(p)>4E$, let $p_0$ be the
segment of $p$ such that $$ \dxh ((p_0)_-, p_-)=\dxh ((p_0)_+,
p_+)=2E.$$ By the second statement of Lemma \ref{CP1} and Lemma
\ref{sk2}, the length of $p_0$ does not exceed the number of
elements of $G$ having length at most $6MLE^2k$ with respect to
the generating set $X$. Thus we have $$l(p_0)\le (card\,
X)^{6MLE^2k}+1$$ and $$|t|_{X\cup\mathcal H}=l(p)\le l(p_0)+4E\le
\rho .$$
\end{proof}

\begin{lem}\label{sk4}
For any $k\ge 0$, $\lambda \ge 1$, $c\ge 0$, there exists a
constant $\eta =\eta (\lambda, c,k)$ such that the following
condition holds. Let $(p,q)$ be an arbitrary symmetric pair of
$k$--similar $(\lambda , c)$--quasi--geodesics without
backtracking in $\G $ such that no synchronous components of $p$
and $q$ are connected. Then for any $i=1,\ldots , m$ and any
$H_i$--component $e$ of $p$, we have
\begin{equation}\label{en}
\dx (e_-, e_+)\le \eta l([p_-, e_+]),
\end{equation}
where $[p_-,e_+]$ is the segment of $p$.
\end{lem}

\begin{proof}
Let $\varepsilon _0= \varepsilon (\lambda , c, k)$ be the constant
given by Theorem \ref{TBCP}. Let also $\varepsilon _1=\varepsilon
(\lambda , c, \max \{ k, \varepsilon _0\} )$. Set $$ \eta = \max
\{ 2\varepsilon _0, \varepsilon _1, 1\} .$$

We proceed by induction on $n=l([p_-, e_+])$. The proof in the
case $n=1$ is given below together with the proof in the general
case.

First suppose that no component of $q$ is connected to $e$. Then
the inequality (\ref{en}) follows from Theorem \ref{TBCP}.
Further, assume that there exists an $H_i$--component $e^\prime $
of $q$ connected to $e$. Let $m$ be the length of the segment
$[q_-, e^\prime _+]$ of $q$. By our assumption, $m\ne n$. Thus
there are two possibilities (see Fig. \ref{sk4-fig}).

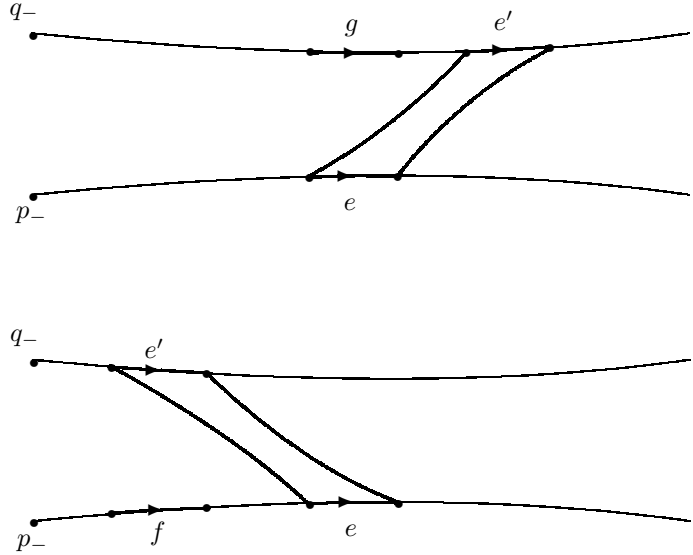
\begin{figure}
\unitlength .9mm 
\linethickness{0.4pt}
\ifx\plotpoint\undefined\newsavebox{\plotpoint}\fi 
\begin{picture}(104.12,92)(-12,-6)
\qbezier(7.07,5.83)(64.26,11.49)(104.12,5.83)
\qbezier(6.95,54.08)(64.14,59.74)(103.99,54.08)
\qbezier(7.07,29.7)(64.26,24.04)(104.12,29.7)
\qbezier(6.95,77.95)(64.14,72.29)(103.99,77.95)
\put(18.74,6.89){\circle*{1}} \put(7.24,5.52){\circle*{1}}
\put(7.12,53.77){\circle*{1}} \put(7.24,29.27){\circle*{1}}
\put(7.12,77.52){\circle*{1}} \put(32.74,7.77){\circle*{1}}
\put(47.99,8.27){\circle*{1}} \put(47.87,56.52){\circle*{1}}
\put(61.11,8.39){\circle*{1}} \put(60.98,56.64){\circle*{1}}
\put(71.09,75.07){\circle*{1}} \put(83.43,75.82){\circle*{1}}
\put(61.13,74.92){\circle*{1}} \put(48.05,75.22){\circle*{1}}
\put(32.74,27.52){\circle*{1}} \put(18.74,28.52){\circle*{1}}
\thicklines
\put(26,28.2){\vector(1,0){.07}}\qbezier(18.75,28.63)(26.25,28.06)(32.75,27.75)
\put(26.38,7.5){\vector(1,0){.07}}\qbezier(18.75,7)(26.94,7.63)(32.88,7.75)
\put(54.5,8.63){\vector(1,0){.07}}\qbezier(47.88,8.38)(54.56,8.81)(61,8.5)
\put(54.38,56.88){\vector(1,0){.07}}\qbezier(47.75,56.63)(54.44,57.06)(60.88,56.75)
\qbezier(18.5,28.75)(37.13,18.44)(47.75,8.38)
\qbezier(32.63,27.75)(47.19,13.81)(61,8.63)
\qbezier(47.72,56.78)(60.13,63.4)(71.06,75.07)
\qbezier(60.95,56.93)(70.24,68.68)(83.39,75.66)
\put(7,2.5){\makebox(0,0)[cc]{$p_-$}}
\put(6.88,51){\makebox(0,0)[cc]{$p_-$}}
\put(5.88,32.88){\makebox(0,0)[cc]{$q_-$}}
\put(5.75,81.13){\makebox(0,0)[cc]{$q_-$}}
\put(25.5,4){\makebox(0,0)[cc]{$f$}}
\put(54,4.5){\makebox(0,0)[cc]{$e$}}
\put(53.88,52.75){\makebox(0,0)[cc]{$e$}}
\put(25,31.38){\makebox(0,0)[cc]{$e^\prime $}}
\put(55.3,75.07){\vector(1,0){.07}}\qbezier(47.87,75.07)(56.19,75.14)(60.95,74.92)
\put(77.15,75.51){\vector(1,0){.07}}\qbezier(70.91,75.22)(77.22,75.44)(83.24,75.96)
\thinlines \put(54.11,78.78){\makebox(0,0)[cc]{$g$}}
\put(76.56,79.83){\makebox(0,0)[cc]{$e^\prime $}}
\end{picture}
\caption{The cases $m>n$ and $m<n$ in the proof of Lemma
\ref{sk4}.} \label{sk4-fig}
\end{figure}

{\it Case 1.} $m<n$ (this case is impossible if $n=1$). Note that
\begin{equation}\label{en1}
\max\{ \dx (e_- , e^\prime _- ), \, \dx (e_+ , e^\prime _+ )\} \le
\varepsilon _0
\end{equation}
by Theorem \ref{TBCP}. We denote by $f$ the $H_i$--component of
$p$ which is synchronous to $e^\prime $. Since $p$ and $q$ are
symmetric, we have
\begin{equation} \label{en2}
\dx (e^\prime _-, e^\prime _+)=\dx (f_-, f_+) \le \eta m
\end{equation}
by the the inductive assumption. Taking (\ref{en1}) and
(\ref{en2}) together, we obtain
$$
\dx (e_-, e_+) \le  2\varepsilon _0 +\eta m \le \eta (m+1)\le \eta
n
$$

{\it Case 2.} $m>n$. Denote by $g$ the $H_i$--component of $q$
which is synchronous to $e$. If no $H_i$--component of the segment
$[p_-,e_-]$ of $p$ is connected to $g$ (in particular, this is so
if $n=1$ since the segment $[p_-,e_-]$ is trivial in this case),
then we obtain
$$
\dx (e-, e_+) = \dx (g_-, g_+) \le \varepsilon _1
$$
by applying Theorem \ref{TBCP} for the segments $[p_-, e_-]$ and
$[q_-, e^\prime _-]$ of $p$ and $q$. The case when $g$ is
connected to a component of the segment $[p_-,e_-]$ can be reduced
to {\it Case 1} by reversing the roles of $p$ and $q$.
\end{proof}

\begin{lem}\label{sk5}
For any $k\in \mathbb N$, there exists a constant $\varkappa
_0=\varkappa _0(k)$ such that the following condition holds. Let
$(p,q)$ be a symmetric pair of $k$--similar geodesics in $\G $
such that no synchronous components of $p$ and $q$ are connected.
Then for every pair of synchronous vertices $u\in p$, $v\in q$, we
have $$ \dx (u,v)\le \varkappa _0.$$
\end{lem}

\begin{figure}

\unitlength 1.00mm 
\linethickness{0.4pt}
\begin{picture}(120,60)(15,20)
\qbezier(16.19,61.81)(82.09,90.29)(136.65,78.83)
\qbezier(15.35,34.06)(81.25,62.54)(135.8,51.08)
\put(83.25,80.94){\circle*{1}} \put(82.41,53.19){\circle*{1}}
\put(72.53,37.84){\circle*{1}} \put(64.75,32.16){\circle*{1}}
\put(16.4,61.81){\circle*{1}} \put(56.24,54.82){\circle*{1}}
\put(55.35,27.47){\circle*{1}} \put(15.14,33.64){\circle*{1}}
\put(135.8,51){\circle*{1}} \put(136.65,78.83){\circle*{1}}

\put(102.8,54.03){\circle*{1}} \put(119.83,53.82){\circle*{1}}

\qbezier(83.25,80.73)(77.47,64.22)(56.13,54.87)
\qbezier(82.41,52.98)(76.63,36.47)(55.29,27.12)
\put(36.37,69.79){\vector(2,1){0}}
\put(35.53,42.04){\vector(2,1){0}}
\put(79.67,73.58){\vector(-3,-4){0}}
\put(78.83,45.83){\vector(-3,-4){0}}
\multiput(56.13,54.66)(-.0336,-1.0848){25}{\line(0,-1){1.0848}}
\multiput(83.04,80.52)(-.033158,-1.427368){19}{\line(0,-1){1.427368}}
\multiput(16.19,61.6)(-.032813,-.860625){32}{\line(0,-1){.860625}}
\multiput(136.86,78.62)(-.032692,-1.051154){26}{\line(0,-1){1.051154}}

\put(59.91,29.28){\vector(-3,-2){0}}

\put(61.24,31.96){\makebox(0,0)[cc]{$r_2$}}
\put(68.68,36.72){\makebox(0,0)[cc]{$a$}}
\put(77.4,47){\makebox(0,0)[cc]{$r_1$}}
\put(78,73){\makebox(0,0)[cc]{$s$}}
\put(111.19,56.49){\makebox(0,0)[cc]{$b$}}
\put(83.39,83.99){\makebox(0,0)[cc]{$v$}}
\put(84.28,55.89){\makebox(0,0)[cc]{$u$}}
\put(29.14,37.46){\makebox(0,0)[cc]{$p$}}
\put(25.72,68.68){\makebox(0,0)[cc]{$q$}}
\put(13.97,63.33){\makebox(0,0)[cc]{$q_-$}}
\put(12.64,33){\makebox(0,0)[cc]{$p_-$}}
\put(137.8,81){\makebox(0,0)[cc]{$q_+$}}
\put(136,48){\makebox(0,0)[cc]{$p_+$}}
\put(94.84,40.43){\makebox(0,0)[cc]{$c$}}

\thicklines \qbezier(102.72,54.11)(90.08,33.59)(64.96,32.11)
\qbezier(102.8,54.3)(111.3,54.3)(119.83,53.6)
\qbezier(72.53,37.84)(69,35)(64.75,32.16)
\put(69.27,35.23){\vector(-4,-3){0}}
\put(92.91,42.51){\vector(-4,-3){0}}
\put(109.85,54.2){\vector(1,0){0}}
\end{picture}
\caption{}\label{kappaprove}
\end{figure}

\begin{proof}
Without loss of generality, we may assume that $q_-=1$. Since the
labels of the segments $[p_-,u]$ and $[1,v]$ of $p$ and $q$
respectively coincide, the element $p_-$ is conjugate to the
element $w=v^{-1}u$ in $G$.

Note that by Lemma \ref{sk1}, the relative length of $w$ satisfies
$ |w|_{X\cup \mathcal H}\le E$. We consider the shortest (with
respect to the relative metric) element $t\in G$ such that
$p_-=t^{-1}wt$. By Lemma \ref{sk3}, the we have
$$
|t|_{X\cup \mathcal H}\le \rho,
$$
where $\rho =\rho (\max \{ E, k\} )$. Let $(r,s)$ be the symmetric
pair of geodesics in $\G $ such that $r_-=u$, $s_-=v$, and the
labels $\phi (r)\equiv \phi (s)$ represent $t$ in $G$. Thus the
label of any path from $s_+$ to $r_+$ in $\G $ represents the
element $p_-$ in $G$.

Suppose that there exist components of $r$ connected to some
components of $p$. We consider the last such a component $a$ of
$r$. Thus $r=r_1ar_2$, there exists a component $b$ of $p$
connected to $a$, and no component of $r_2$ is connected to a
component of $p$. To be definite we assume that $b$ belongs to
$[u, p_+]$. (The case when $b$ belongs to $[p_-,u]$ can be treated
in the same way and is left to the reader.) Let $c$ be a path in
$\G $ of length at most $1$ such that $c_-=b_-$, $c_+=a_+$ (see
Fig. \ref{kappaprove}). Notice that
$$
l([u, b_-])=l([u, b_+])-1\le \dxh (u, a_-)+\dxh (a_-,
b_+)-1=l(r_1)
$$
as $p$ and $r$ are geodesic. Hence
$$
l([u,b_-]cr_2)\le l(r_1)+1+l(r_2)=l(r).
$$
Obviously the path $[p_-,b_-]cr_2 $ is $(1,\rho
)$--quasi--geodesic without backtracking as follows from the
choice of $a$. Thus replacing $r$ with $r^\prime =[u,b_-]cr_2$ and
$s$ with the geodesic symmetric to $r^\prime $, we may assume that
$[p_-,u]r$ has no backtracking.

Let $\varepsilon =\varepsilon (1,\rho , k)$ and $\eta =\eta
(1,\rho, \max \{ \varepsilon , k\} )$ be constants provided be
Theorem \ref{TBCP} and Lemma \ref{sk4} respectively. We consider
two cases.

{\it Case 1}. Suppose that there are no synchronous connected
components of $r$ and $s$. Applying Lemma \ref{sk4} for the
symmetric $(1, \rho )$--quasi--geodesic paths $r^{-1}[p_-,u]^{-1}$
and $s^{-1}[q_-,v]^{-1}$, we obtain that the $X$--length of any
$H_i$--component $e$ of $r$ satisfies the inequality
$$
\dx (e_-, e_+)\le \eta l(r)\le \eta\rho .
$$
Hence,
$$
\dx (u, r_+)\le l(r)\eta\rho \le \eta\rho ^2.
$$
Finally we have
$$
\dx (u,v)\le  \dx (u,r_+)+\dx (r_+,s_+)+\dx (s_+,v)\le 2\eta\rho
^2 +k.
$$

{\it Case 2.} Now assume that there is at least one pair of
synchronous connected components in $r$ and $s$. Let $i$, $j$ be
the connected synchronous components of $r$ and $s$ respectively
such that there are no connected synchronous components of the
segments $[u,i_-]$ and $[v,j_-]$ of $r$ and $s$. By Theorem
\ref{TBCP}, we have $\dx (i_-, j_-)\le \varepsilon $. Therefore,
the paths $[p_-,u][u,i_-]$ and $[1,v][v,j_-]$ form a symmetric
$\max \{ \varepsilon , k\} $--similar pair. Arguing as in the Case
1, we obtain
$$
\dx (u,v) \le 2\eta \rho ^2 + \varepsilon .
$$
In both cases it suffices to set $\varkappa _0=2\eta \rho ^2 +
\max \{ \varepsilon , k\} $.
\end{proof}

Now we are ready to prove the main result of this section.

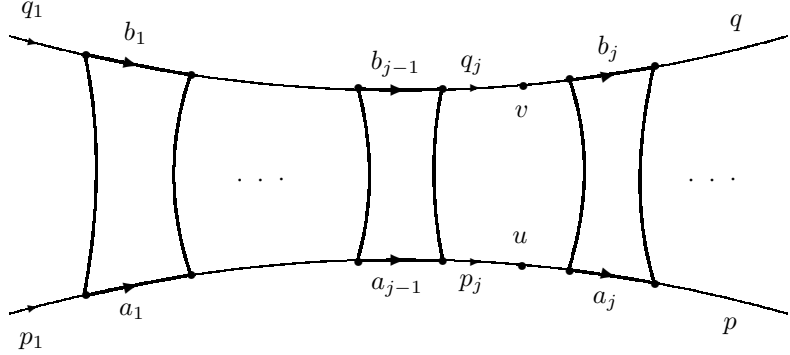
\begin{figure}
\unitlength 1.1mm 
\linethickness{0.4pt}
\ifx\plotpoint\undefined\newsavebox{\plotpoint}\fi 
\begin{picture}(100.06,42.57)(0,0)
\qbezier(6.19,5.83)(53.65,18.92)(100.06,5.83)
\qbezier(6.19,39.36)(53.65,26.28)(100.06,39.36)
\put(15.38,8.13){\circle*{.9}} \put(15.38,37.07){\circle*{.9}}
\put(28.1,10.55){\circle*{.9}} \put(28.1,34.65){\circle*{.9}}
\put(48.28,12.13){\circle*{.9}} \put(48.28,33.07){\circle*{.9}}
\put(68.25,33.39){\circle*{.9}} \put(68.04,11.63){\circle*{.9}}
\put(58.48,12.23){\circle*{.9}} \put(58.48,32.97){\circle*{.9}}
\put(73.82,10.97){\circle*{.9}} \put(73.82,34.23){\circle*{.9}}
\put(84.12,9.39){\circle*{.9}} \put(84.12,35.81){\circle*{.9}}
\put(9.57,6.73){\vector(4,1){.07}}
\put(9.57,38.47){\vector(4,-1){.07}}

\put(63,12.1){\vector(1,0){.07}} \put(63,33.1){\vector(1,0){.07}}

\thicklines
\put(21.76,9.35){\vector(4,1){.07}}\qbezier(15.24,7.99)(21.86,9.46)(28.06,10.51)
\put(21.76,35.84){\vector(4,-1){.07}}\qbezier(15.24,37.21)(21.86,35.74)(28.06,34.69)
\put(54.03,12.3){\vector(1,0){.07}}\qbezier(48.14,12.19)(54.76,12.4)(58.44,12.19)
\put(54.03,32.9){\vector(1,0){.07}}\qbezier(48.14,33.01)(54.76,32.79)(58.44,33.01)
\put(79.46,10.3){\vector(4,-1){.07}}\qbezier(73.68,11.14)(80.15,10.25)(83.88,9.57)
\put(79.46,34.9){\vector(4,1){.07}}\qbezier(73.68,34.06)(80.15,34.95)(83.88,35.63)
\qbezier(15.24,37)(18.03,23.86)(15.14,8.2)
\qbezier(27.96,34.69)(23.97,22.76)(27.96,10.62)
\qbezier(48.25,33.01)(51.03,22.18)(48.14,12.19)
\qbezier(58.34,33.01)(56.39,22.91)(58.44,12.4)
\qbezier(73.79,34.16)(77.42,23.39)(73.68,11.14)
\qbezier(83.98,35.74)(80.57,20.5)(84.09,9.46)
\put(33.95,21.97){\circle*{.4}} \put(88.5,21.97){\circle*{.4}}
\put(36.47,22.07){\circle*{.4}} \put(91.03,22.07){\circle*{.4}}
\put(39,22.18){\circle*{.4}} \put(93.55,22.18){\circle*{.4}}
\put(8.93,42.57){\makebox(0,0)[cc]{$q_1$}}
\put(8.83,2.52){\makebox(0,0)[cc]{$p_1$}}
\put(20.92,6.1){\makebox(0,0)[cc]{$a_1$}}
\put(21.34,39.73){\makebox(0,0)[cc]{$b_1$}}
\put(52.77,36.05){\makebox(0,0)[cc]{$b_{j-1}$}}
\put(52.98,9.04){\makebox(0,0)[cc]{$a_{j-1}$}}
\put(78.52,37.84){\makebox(0,0)[cc]{$b_j$}}
\put(78.1,7.36){\makebox(0,0)[cc]{$a_j$}}
\put(93.23,4.31){\makebox(0,0)[cc]{$p$}}
\put(93.97,40.99){\makebox(0,0)[cc]{$q$}}
\put(62,35.84){\makebox(0,0)[cc]{$q_j$}}
\put(62,9.35){\makebox(0,0)[cc]{$p_j$}}
\put(68.01,30.27){\makebox(0,0)[cc]{$v$}}
\put(67.8,15.03){\makebox(0,0)[cc]{$u$}}
\end{picture}
\caption{The decomposition of $p$ and $q$ in the proof of Theorem
\ref{sc}.} \label{sc-fig}
\end{figure}

\begin{proof}[Proof of Theorem \ref{sc}]
Let $a_1, \ldots , a_l$ be the set of all components of $p$ that
are connected to the corresponding synchronous components of $q$.
We denote by $b_j$ the component of $q$ connected to $a_j$.
Changing if necessary the order of enumeration of $a_1, \ldots ,
a_l$, we may assume that
$$
p=p_1a_1\ldots p_la_lp_{l+1},
$$
$$
q=q_1b_1\ldots q_lb_lq_{l+1},
$$
where for any $j=1, \ldots , l+1$,  $(p_j, q_j)$ is a symmetric
pair of $\varepsilon $--similar geodesics for $\varepsilon
=\varepsilon (1,0,k)$ given by Theorem \ref{TBCP} (see Fig 3.15).
Notice that there exist no connected synchronous components of
$p_j$ and $q_j$ for any $j=1, \ldots , l+1$. As $u$ and $v$ are
synchronous vertices of $p_j$ and $q_j$ for a certain $j$, it
suffices to set $\varkappa =\varkappa _0(\max \{ \varepsilon , k\}
)$, where $\varkappa _0(\max \{ \varepsilon , k\} )$ is provided
by Lemma \ref{sk5}.
\end{proof}


\chapter{Algebraic properties}



\section{Elements of finite order}


Recall that the number of conjugacy classes of elements of finite
order in any hyperbolic group is finite. A generalization of this
result to the class of hyperbolic products of groups can be found
in \cite{PanPhD}. In this section we extend these results by
proving the following.

\begin{defn}
Let $G$ be a group hyperbolic with respect to a collection of
subgroups $\{ H_\lambda \} _{\lambda \in \Lambda }$. An element
$g\in G$ is called {\it parabolic} if it is conjugate to an
element of one of the subgroups $H_\lambda $. Otherwise $g$ is
said to be {\it hyperbolic}.
\end{defn}

\begin{thm}\label{FinOrd}
Suppose $G$ is a group hyperbolic with respect to a collection of
subgroups $\{ H_\lambda \} _{\lambda \in \Lambda }$. Then the
number of conjugacy classes of hyperbolic elements of finite order
in $G$ is finite.
\end{thm}

It is well known that if $G$ acts on a tree without inversions and
with a compact quotient, then every element of finite order has a
fixed point (see \cite{Trees} or \cite{DD}). In other terms, if
$G$ is a fundamental group of a a graph of groups, then each
element of finite order in $G$ is conjugate to an element of one
of the vertex groups. In combination with Splitting Theorem from
Section 2.4, this shows that it suffices to prove the theorem in
case $G$ is finitely generated (and therefore the collection of
subgroups is finite by Corollary \ref{Lambda}). In the rest of
this section we assume that $G$ is generated by a finite set $X$
and is hyperbolic relative to subgroups $\{ H_1, \ldots , H_m\} $.

Recall, that a path in a metric space is said to be $k$--local
geodesic if any its subpath of length at most $k$ is geodesic. The
following lemma is well known (see, for example, \cite[CH. III.H,
Theorem 1.13]{BriH}).

\begin{lem}\label{k-loc}
Let $r$ be a $k$--local geodesic in a $\delta $--hyperbolic metric
space for some $k>8\delta $. Then $r$ is $(\lambda ,
c)$--quasi--geodesic for $\lambda =\frac{k+4\delta }{k-4\delta
}\le 3$ and $c=2\delta $.
\end{lem}

\begin{cor}\label{FinOrd1}
Let $g$ be an element of finite order in $G$. Then $G$ is
conjugate to an element of relative length at most $8\delta +1$,
where $\delta $ is the hyperbolicity constant for $\G $.
\end{cor}

\begin{proof}
We assume that $g$ is a shortest element in the conjugacy class
$g^G$ and $|g|_{X\cup\mathcal H}>8\delta +1$. Let us take a
shortest word $U\in (X\cup\mathcal H)^\ast $ representing $g$ and
consider the path $p_n$ such that $(p_n)_-=1$, $\phi (p_n)\equiv
U^n$, $n\in \mathbb N$. Since $g$ is a shortest element in $g^G$,
$p_n$ is $(8\delta +1)$--local geodesic in $\G $. Therefore, by
Lemma \ref{k-loc}, we have
$$
|g^n|_{X\cup \mathcal H}=\dxh (1, (p_n)_+) \ge \frac13l(p_n)
-2\delta\ge \frac13|n|-2\delta .
$$
Note that $1/3 |n|-2\delta \ne 0$ for any $n$ big enough. Hence
the order of $g$ is infinite.
\end{proof}

In contrast to the case of ordinary hyperbolic groups, the above
corollary does note imply the desired result since balls in $\G $
are, in general, not finite. The next lemma is the crucial
ingredient of the proof of Theorem \ref{FinOrd}.

\begin{lem}\label{FinOrd2}
Suppose that a group $G$ is generated by a finite set $X$ and is
hyperbolic relative to subgroups $\{ H_1, \ldots , H_m\} $. Then
there exists a constant $B$ such that the following condition
holds. Let $f$ be a hyperbolic element of finite order in $G$ such
that $f$ has smallest relative length among all elements of the
conjugacy class $f^G$. Then
$$|f|_X\le B|f|_{X\cup \mathcal H} ^2.$$
\end{lem}

\begin{proof}

\begin{figure}

\unitlength 0.8mm 

\begin{picture}(90,45)(-5,23)
\put(62.88,49.75){\oval(72.25,40.5)[]}

\put(26.75,58.75){\vector(0,1){0}}
\put(26.75,41.63){\vector(0,1){0}}

\put(98.88,52.38){\vector(0,-1){0}}

\thicklines \put(62,70){\vector(1,0){0}}
\put(61.5,29.5){\vector(-1,0){0}}
\put(68.25,70){\line(-1,0){14.8}} \put(68.25,
29.5){\line(-1,0){14.8}} \put(68.13,29.5){\line(0,1){40.33}}
\put(68.13,51.88){\vector(0,-1){0}}
\put(26.63,48.88){\circle*{1.2}} \put(53.88,69.88){\circle*{1.2}}
\put(68.13,69.88){\circle*{1.2}} \put(68.13,29.5){\circle*{1.2}}
\put(53.88,29.5){\circle*{1.2}}

\put(24,58){\makebox(0,0)[cc]{$a$}}
\put(24,40){\makebox(0,0)[cc]{$c$}}
\put(106.5,53.13){\makebox(0,0)[cc]{$p_0=b$}}
\put(61,72.63){\makebox(0,0)[cc]{$s_1$}}
\put(61.13,26.75){\makebox(0,0)[cc]{$s_2$}}
\put(70.5,51.5){\makebox(0,0)[cc]{$e$}}
\end{picture}
\label{EFO} \caption{}
\end{figure}

Let $U$ be a shortest word in $(X\cup\mathcal H)^\ast $
representing $f$ in $G$. Let $n$ be the order of $f$, $l$ the
relative length of $f$. We consider an arbitrary cycle $p$ in $\G
$ with the label $U^n$ and take a subpath $p_0$ of $p$ obtained as
follows. If $p$ is a cycle without backtracking, we set $p_0=p$.
Further suppose that there exists two connected components $s_1$
and $s_2$ of $p$. Let
$$p=as_1bs_2c.$$ Passing to another pair of connected components if necessary,
we may assume that $b$ is a path without backtracking and no
component of $b$ is connected to the component $s_1$, $s_2$. In
this case we set $p_0=b$.

Since $\phi (p_0)$ is a (cyclic) subword of $U^n$, we have $$ \phi
(p_0)\equiv U_0^kV,$$ where $U_0$ is a cyclic shift of $U$ and $V$
is a cyclic subword of $U$ of length less than $l$. It is easy to
check that $k>0$. Indeed if $k=0$, then $f$ is not the shortest
element in the conjugacy class $f^G$.

Let $e$ be a path in $\G $ of length at most $1$ such that
$e_-=(p_0)_-$ and $e_+=(p_0)_+$ (see Fig. \ref{EFO}). We consider
the cycle $c=p_0e^{-1}$. According to our choice of $p_0$, $c$ is
a cycle without backtracking. We have
$$
l(c)=l(p_0)+1\le kl+\| V\| +1\le (k+1)l.
$$

Let $W$ be an $H_j$--syllable in $U$. Then $W$ is a label of at
least $k$ components of $c$. Since all components of $c$ are
isolated, applying Lemma \ref{31} we obtain
$$
k|\overline{W}|_X\le MLl(c)\le ML(k+1)l.
$$
Hence,
\begin{equation}\label{FO1}
|\overline{W}|_X \le ML\frac{k+1}{k}l\le 2MLl.
\end{equation}
Finally, since the inequality (\ref{FO1}) holds for any syllable
$W$ of $U$, we have $ |\overline{U}|_X\le 2MLl^2$.
\end{proof}

\begin{proof}[Proof of Theorem \ref{FinOrd}]
By combining Lemmas \ref{FinOrd1} and \ref{FinOrd2}, we obtain
that each hyperbolic element of finite order in $G$ is conjugate
to an element of $X$--length at most $B(8\delta +1)^2$. Since $G$
is locally finite with respect to the metric $\dx $, this proves
the theorem.
\end{proof}

As a corollary, we have

\begin{cor}\label{SetO}
The set of orders of hyperbolic elements in $G$ is finite.
\end{cor}

\begin{cor}
If $G$ is residually finite and all subgroups $H_\lambda $ are
torsion free, then $G$ is virtually torsion free, that is $G$
contains a torsion free subgroup of finite index.
\end{cor}

\begin{proof}
Let $g_1, g_2,\ldots , g_k$ be elements of $G$ such that each
hyperbolic element of finite order is conjugate to one of them.
Then there exists a normal subgroup $N$ of finite index in $G$
such that $g_i\notin N$ for all $i=1,\ldots ,k$. Thus $N$ contains
no hyperbolic elements of finite order. Hence $N$ is torsion free.
\end{proof}

Note that the requirement of residual finiteness is essential.

\begin{ex}
Let $H$ be a finitely generated torsion free simple group. Let $w$
be an arbitrary nontrivial element of $H$. Consider the group $G$
given by the relative presentation
$$
G=\langle x, H\; |\; x^2=1, \; x=[x,w][x,w^2]\ldots [x,w^n]
\rangle .
$$
It is easy to check that $G$ is a quotient of the free product
$\langle x\; | \; x^2=1\rangle  \ast H$ by the relation
$x=[x,w][x,w^2]\ldots [x,w^n]$ satisfying $C^\prime (\lambda )$
hypothesis (as a relation over a free product) with $\lambda \to
0$ as $n\to \infty $. Therefore, by the Greendlinger Lemma
\cite{LS}, the relative Dehn function of $G$ with respect to $H$
is linear for all $n$ big enough and thus $G$ is hyperbolic
relative to $H$. However, any subgroup $N$ of finite index in $G$
contains $H$ (otherwise $N\cap H=\{ 1\} $ and hence $G/N$ is
infinite). Now the relation $x=[x,w][x,w^2]\ldots [x,w^n]$ implies
that $x\in G$. Thus $N=G$, i.e., $G$ contains no proper subgroups
of finite index. In particular, $G$ is not virtually torsion free.
\end{ex}


\section{Relatively quasi--convex subgroups}


Our discussion in this section is stimulated by some ideas of
Gromov \cite{Gr1} which were elaborated by Gersten and Short
\cite{GS}, Alonso and Bridson \cite{AloB}, and others (see
\cite{Ep-etal} and references therein). Our main goal here is to
introduce the (geometric) notion of a quasi--convex subgroup of a
relatively hyperbolic group and to obtain some analogues of
well--known theorems about quasi--convex subgroups in hyperbolic
groups. In the next section, adopting an idea of Gersten and Short
\cite{GS} to the relative case, we apply our results to the study
of translation numbers. For the dynamical notion of
quasi--convexity for convergence groups we refer to \cite{B96}.

\begin{defn}\label{qc}
Let $G$ be a group generated by a finite set $X$, $\{ H_1, \ldots
, H_m\} $ a collection of subgroups of $G$. A subgroup $R$ of $G$
is called relatively quasi--convex with respect to $\{ H_1, \ldots
, H_m\} $ (or simply {\it relatively quasi--convex} when the
collection $\{ H_1, \ldots , H_m\} $ is fixed) if there exists a
constant $\sigma
>0$ such that the following condition holds. Let $f$, $g$ be two
elements of $R$, and $p$ an arbitrary geodesic path from $f$ to
$g$ in $\G $. Then for any vertex $v\in p$, there exists a vertex
$w\in R$ such that
$$\dx (u,w)\le \sigma .$$
 Note that, without loss of generality, we may assume one of
the elements $f,g$ to be equal to the identity since both the
metrics $\dx $ and $\dxh $ are invariant under the left action of
$G$ on itself.
\end{defn}

It is easy to see that, in general, this definition depends on
$X$. However in case of relatively hyperbolic groups we have

\begin{prop}\label{qcinv}
Let $G$ be a group hyperbolic with respect to a collection of
subgroups $\{ H_1, \ldots , H_m\} $ and $R$ a subgroup of $G$.
Suppose that $X_1, X_2$ are two finite generating sets of $G$.
Then $Q$ is relatively quasi--convex with respect to $X_1$ if and
only if it is relatively quasi--convex with respect to $X_2$.
\end{prop}

\begin{proof}
Let $\Gamma _1=\Gamma (G, X_1\cup \mathcal H)$ and $\Gamma
_2=\Gamma (G, X_2\cup \mathcal H)$. For every $x\in X_1$, we fix a
word $W_x$ over $X_2$ representing $x$ in $G$. To each path $p$ in
$\Gamma _1$, we assign a path in $\Gamma _2$ which starts and ends
at the same elements as $p$ and has label obtained from $\phi (p)$
by replacing $x$ with $W_x$ for every $x\in X_1$.

Suppose that $R$ is quasi--convex with respect to $X_2$. Let $r$
be an element of $R$, $p$ a geodesic path in $\Gamma _1$ such that
$p_-=1$ and $p_+=r$. We also take a vertex $u\in p$. Denote by $q$
the path in $\Gamma _2$ corresponding to $p$, and by $v$ the
vertex corresponding to $v$ (thus $u=v$ being considered as
elements of $G$). By Proposition \ref{Lip}, $q$ is $(\lambda ,
0)$--quasi--geodesic for some constant $\lambda $ which is
independent of $p$. Moreover, since $p$ is a path without
backtracking, then obviously so is $q$.

By Theorem \ref{TBCP}, $q$ lies in the closed $\varepsilon
=\varepsilon (\lambda , 0,0)$--neighborhood (with respect to the
metric $dist _{X_2}$) of the geodesic path $s$ in $\Gamma _2$ with
$s_-=1$, $s_+=r$. As $R$ is quasi--convex with respect to $X_2$,
$s$ belongs to the closed $\sigma $--neighborhood (with respect to
$dist _{X_2}$) of $R$, where $\sigma $ is the quasi--convexity
constant. Thus
$$
dist_{X_2}(v,R)\le \varepsilon +\sigma .
$$
Applying Proposition \ref{Lip} again, we obtain the upper bound on
$dist _{X_1}(u,R)$, which is independent of $u$ and $p$. Thus $R$
is quasi--convex with respect to $X_1$.
\end{proof}

\begin{defn}
Let $G$ be as in the Definition \ref{qc}. A relatively
quasi--convex subgroup $R$ of $G$ is called {\it strongly
relatively quasi--convex} if the intersection $R\cap H_i^g$ is
finite for any $g\in G$, $i=1, \ldots , m$.
\end{defn}

\begin{defn}
Recall that a map $\iota : (M_1, dist_1)\to (M_2, dist_2) $
between two metric spaces $M_1$ and $M_2$ with metrics $dist_1$
and $dist_2$ is called a {\it quasi--isometric embedding} if there
exist $c_1,c_2 >0$ such that for every two points $x,y\in M_1$ we
have
$$
\frac{1}{c_1}dist _1(x,y) -c_2\le dist _2(\iota (x), \iota (y))\le
c_1dist _1(x,y)+c_2.
$$
\end{defn}

\begin{thm}\label{qc1}
Suppose that the group $G$ is hyperbolic relative to the
collection of subgroups $\{ H_1, \ldots , H_m\} $. Let $R$ be a
subgroup of $G$. Then the following conditions are equivalent.
\begin{enumerate}

\item $R$ is strongly relatively quasi--convex.

\item $R$ is generated by a finite set $Y$ and the natural map
$(R, dist _Y) \to (G, \dxh )$ is a quasi--isometric embedding.
\end{enumerate}
\end{thm}

\begin{proof}
For every $x,y\in G$ and every $i=1, \ldots , m$ we consider the
set
$$
Z_{x,y,i}=\{ xhy\; | \; h\in H_i\} \cap R.
$$
Also set
$$
Z_0=\{ r\in R\; | \; |r|_X\le 2\sigma +1\} ,
$$
where $\sigma $ is the quasi--convexity constant for $R$.

To prove the theorem we need two auxiliary lemmas.

\begin{lem}\label{qc11}
Let
\begin{equation}\label{QC10}
B_\sigma ^X =\{ g\in G\; |\; |g|_X\le \sigma \} .
\end{equation}
Then the subgroup $R$ is generated by the set
$$
Z=Z_0\bigcup \left( \bigcup\limits_{x,y\in B_\sigma ^X}
\bigcup\limits_{i=1}^m Z_{x,y,i} \right) .
$$
\end{lem}

\begin{proof}
Let $r$ be an arbitrary element of $R$. We consider a geodesic $p$
in $\G $ such that $p_-=1$, $p_+=r$. Let $g_0=1, g_1, \ldots,
g_n=r$ be the consecutive vertices of $p$ (see Fig. . By the
definition of a relatively quasi--convex subgroup, for any $i=1,
\ldots , n-1$, there exists an element $r_i\in R$ such that
\begin{equation}\label{QC1}
\dx (r_i, g_i)\le \sigma .
\end{equation}
We also set $r_0=1$ and $r_n=r$. Denote by $x_i$ the element
$r_i^{-1}g_i $ and by $e_{i+1}$ the edge of $p$ going from $g_i$
to $g_{i+1}$.

Obviously we have
$$
r_{i+1}=r_ix_i\overline{\phi (e_{i+1})}x_{i+1}^{-1}.
$$
By (\ref{QC1}), the $X$--length of $x_i$ satisfies $|x_i|_X\le
\sigma $. Therefore the element
$$
s_i=x_i\overline{\phi (e_{i+1})}x_{i+1}^{-1}
$$
either belongs to $Z_0$ in case $\phi (e_{i+1})\in X$, or belongs
to $Z_{x_i, x_{i+1}, j} $ if $\phi (e_{i+1})\in \widetilde H_j$
for some $j=1, \ldots , m$. Thus in both cases we have $s_i\in Z$.

\begin{figure}
\unitlength 1.00mm 
\begin{picture}(129.08,40)(10,20)
\multiput(12.4,27.33)(.0337218045,.0442606516){798}{\line(0,1){.0442606516}}
\put(39.31,62.65){\line(1,0){89.77}}
\multiput(129.08,62.65)(-.0337172012,-.0517930029){686}{\line(0,-1){.0517930029}}
\multiput(105.95,27.12)(-13.39429,-.03){7}{\line(-1,0){13.39429}}
\multiput(31.66,33)(.0336129,.08058065){155}{\line(0,1){.08058065}}

\multiput(36.87,45.49)(.036330935,.033669065){278}{\line(1,0){.036330935}}
\multiput(46.97,54.85)(.07174419,.03372093){172}{\line(1,0){.07174419}}
\multiput(59.31,60.65)(.16277372,.03364964){137}{\line(1,0){.16277372}}
\multiput(97.51,65.26)(.0796129,-.0336129){155}{\line(1,0){.0796129}}
\multiput(109.85,60.05)(.03361809,-.03432161){199}{\line(0,-1){.03432161}}
\multiput(37.16,45.49)(-.03,-3.242){5}{\line(0,-1){3.242}}

\put(37.1,38){\vector(0,1){.1}} \put(47.12,45){\vector(0,1){.1}}
\put(38,38){$x_1$} \put(48,45){$x_2$}

\multiput(37.01,29.28)(-.0481982,.03351351){111}{\line(-1,0){.0481982}}
\put(47.12,54.7){\line(0,-1){19.32}}
\multiput(47.12,35.38)(-.05711864,-.03361582){177}{\line(-1,0){.05711864}}
\multiput(59.16,60.65)(-.03222,-3.07222){9}{\line(0,-1){3.07222}}
\multiput(58.87,33)(-.165493,.0335211){71}{\line(-1,0){.165493}}
\multiput(109.85,59.91)(.03333,-1.32222){9}{\line(0,-1){1.32222}}
\multiput(110.15,48.01)(.03804878,.03353659){164}{\line(1,0){.03804878}}
\put(86.37,66){\circle*{.4}} \put(81.76,41.03){\circle*{.4}}
\put(88.89,66){\circle*{.4}} \put(84.29,41.03){\circle*{.4}}
\put(91.42,66){\circle*{.4}} \put(86.81,41.03){\circle*{.4}}
\put(100.49,31.07){\makebox(0,0)[cc]{$R$}}
\put(31.96,32.85){\circle*{1}} \put(36.87,29.28){\circle*{1}}
\put(46.97,35.38){\circle*{1}} \put(59.01,33){\circle*{1}}
\put(110.3,48.01){\circle*{1}} \put(116.54,53.37){\circle*{1}}
\put(110,59.91){\circle*{1}} \put(59.31,60.8){\circle*{1}}
\put(47.12,54.85){\circle*{1}} \put(37.16,45.64){\circle*{1}}
\put(71.8,66.74){\makebox(0,0)[cc]{$p$}}
\put(30,33){\makebox(0,0)[cc]{$1$}}
\put(119.07,54){\makebox(0,0)[cc]{$r$}}
\put(41,28.84){\makebox(0,0)[cc]{$r_1$}}
\put(47.87,32){\makebox(0,0)[cc]{$r_2$}}
\put(35.2,47){\makebox(0,0)[cc]{$g_1$}}
\put(44.89,56.6){\makebox(0,0)[cc]{$g_2$}}
\put(33.59,37.61){\vector(1,2){.07}}
\put(41.62,49.95){\vector(1,1){.07}}
\put(31.66,38.95){\makebox(0,0)[cc]{$e_1$}}
\put(40.58,52.18){\makebox(0,0)[cc]{$e_2$}}
\put(115.5,58){\makebox(0,0)[cc]{$e_n$}}
\put(113.12,56.64){\vector(3,-4){.07}}

\end{picture}
\caption{}\label{QC-fig1}
\end{figure}

It is clear that
\begin{equation} \label{QC0}
r_n=r_{n-1}s_{n-1}=\ldots = r_0s_0\ldots s_{n-1}=s_0\ldots s_{n-1}
\end{equation}
as $r_0=1$. Therefore $r_n\in \langle Z\rangle $. Since $r$ was an
arbitrary element of $R$, the lemma is proved.
\end{proof}

\begin{lem}\label{qc12}
For any $x,y\in G$ and any $i=1, \ldots , m$ the set $Z_{x,y,i}$
is finite.
\end{lem}

\begin{proof}
Assume that $Z_{x,y,i}=\{ z_0,z_1, \ldots \} $ is infinite for
some $x,y\in G$, $i\in \{ 1, \ldots , m\} $. Suppose that for any
$j=0,1, \ldots $,  $z_j=xh_jy$, where $h_j\in H_i$. The elements
$t_j=z_0^{-1}z_j$ are different for different $j\ge 0$. However we
have
$$
t_j=y^{-1}h_0x^{-1}xh_jy=y^{-1}h_0h_jy\in H_i^y
$$
for every $j$. Recall that $t_j\in R$ for any $j$ by the
definition of $Z_{x,y,i}$. Therefore, the intersection $R\cap
H_i^y$ is infinite contrary to our assumption.
\end{proof}

Now we are able to show that the first condition in Theorem
\ref{qc1} implies the second one. Using Lemma \ref{qc12}, one can
easily see that the set $Z$ is finite since the ball $B_\sigma ^X$
is finite. Further by Lemma \ref{qc11}, $R$ is generated by $Z$.
Therefore $R$ is finitely generated. Moreover, it follows from the
proof of Lemma \ref{qc11} (see (\ref{QC0})) that for any element
$r\in R$ of relative length $|r|_{X\cup \mathcal H}=n$, we have
$r=s_0\ldots s_{n-1}$ for certain $s_0, \ldots , s_{n-1}\in Z$.
This means that
$$
|r|_Z\le |r|_{X\cup \mathcal H} ,
$$
i.e., the map $(R, dist_Z)\to (G, \dxh )$ is a quasi--isometric
embedding.

To prove the converse implication, we assume that $R$ is generated
by a finite set $Y=Y^{-1}$ and the natural map $(R, dist _Y)\to
(G, \dxh )$ is a quasi--isometric embedding. We denote by $c_1,
c_2$ the corresponding quasi--isometry constants.

Note that if $r\in R\cap H_i^g$ for some $i\in \{ 1, \ldots , m\}
$ and $g\in G$, then
$$
|r|_Y\le c_1|r|_{X\cup \mathcal H} +c_2\le c_1(2|g|_{X\cup
\mathcal H} +1)+c_2.
$$
Since $Y$ is finite, we obtain $card\, (R\cap H_i^g) <\infty $ for
$i=1, \ldots , m$ and any $g\in G$.

It remains to show that $R$ is relatively quasi--convex. For any
element $y\in Y$, we fix a word $W_y$ over the alphabet $X\cup
\mathcal H$ representing $y$ in $G$. Set
\begin{equation}\label{mu}
\mu =\max\limits_{y\in Y} \| W_y\| .
\end{equation}
Given an element $r\in R$, we consider the shortest word
$V=y_1\ldots y_n$ over $Y$ representing $r$. Let $$U=W_{y_1}
\ldots W_{y_n}$$ be the word over $X$ obtained from $V$ by
replacing each $y_i$ with the corresponding $W_{y_i}$.

Let $U_0$ be a subword of $U$. Then
$$
U_0\equiv AW_{y_j}\ldots W_{y_{j+k}}B,
$$
where $\| A\|< \mu $ and $\| B\| < \mu $. Since any subword of $V$
is geodesic with respect to the metric $dist _Y$ on $R$, we have
$$
\begin{array}{rl}
\| U_0\| < & 2\mu +\mu (k+1) = \\ & \\ & 2\mu +\mu |y_j\ldots
y_{j+k}|_Y \le
\\ & \\ & 2\mu +\mu \left( c_1|y_j\ldots y_{j+k} |_{X\cup \mathcal
H} +c_2\right) \le \\&\\ & 2\mu +\mu \left( c_1\left(
|\overline{U_0} |_{X\cup \mathcal H} +2\mu \right )+ c_2\right) .
\end{array}
$$

Thus the path $p$ in $\G $ with $p_-=1$ labelled $U$ is $(\mu c_1,
2\mu +2\mu ^2c_1 +\mu c_2)$--quasi--geodesic. Let $\varepsilon =
\varepsilon (\mu c_1, 2\mu +2\mu ^2c_1 +\mu c_2,0) $ be the
constant provided by Proposition \ref{QG}. Then for any geodesic
path $q$ in $\G $ such that $q_-=1$ and $q_+=r$, and any vertex
$v\in q$, there exists a vertex $u\in p$ such that
$$
\dx (u,v)\le \varepsilon .
$$
It is clear that any vertex of $p$ belongs to the closed $\mu
$--neighborhood of $R$ with respect to the metric $\dx $. Hence,
$$
\dx (v, R)\le \varepsilon +\mu.
$$
Since the right hand side of the above inequality is independent
of $r$, $R$ is relatively quasi--convex.
\end{proof}

As is well known, any quasi--convex subgroup of a hyperbolic group
is hyperbolic itself. The theorem below generalizes this result.

\begin{thm}\label{qc-hyp}
Let $R$ be a strongly relatively quasi--convex subgroup of $G$.
Then $R$ is a hyperbolic group.
\end{thm}

\begin{proof}
By Theorem \ref{qc1}, $R$ is generated by a finite set $Y$. As in
the proof of Theorem \ref{qc1}, for any $y\in Y$, let $W_y$ denote
a word over $X\cup \mathcal H$ representing $y$ in $G$ and let
$\mu $ be defined by (\ref{mu}). We denote by $\Gamma (R, Y)$ the
Cayley graph of $R$ with respect to $Y$ and define the map
$$\psi : \Gamma (R, Y)\to \G $$ as follows. First we require the
restriction of $\psi $ to the vertex set of $\Gamma (R,Y)$ to
coincide with the natural embedding $R\to G$. Secondary, given an
edge $e$ of $\Gamma (R,Y)$ with label $\phi (e)=y$, we assume
$\psi (e)$ to be the path in $\G $ labelled $W_y$. Obviously this
two conditions uniquely define $\psi $. Note that $\psi $ is a
quasi--isometric embedding by Theorem \ref{qc1}. We denote by
$c_1, c_2$ the corresponding quasi--isometry constants.

Let $\Delta =pqr$ be a geodesic triangle in $\Gamma (R, Y) $. Then
the image $\Delta ^\ast =p^\ast q^\ast r^\ast $ of $\Delta $ under
$\psi $ is a triangle in $\G $ whose sides are $(\lambda ,
c)$--quasi--geodesics, where $\lambda $ and $c$ depend only on
$c_1, c_2$, not on $\Delta $. Since $\G $ is a $\delta
$--hyperbolic space, by Lemma \ref{qg}, there is a constant
$H=H(\lambda , c,0)$ such that each side of $\Delta ^\ast $
belongs to the closed $H$--neighborhood (with respect to $\dxh $)
of the other two sides.

Let $u$ be a vertex on one of the sides of $\Delta $, say $u\in
p$, and $u^\ast =\psi (u)$. Let also $v^\ast $ be a vertex on
$q^\ast \cup r^\ast $ such that
\begin{equation}\label{QC21}
\dxh (u^\ast , v^\ast )\le H.
\end{equation}
Clearly there exists a vertex $w^\ast \in q^\ast\cup r^\ast $ such
that $w^\ast =\psi (w) $ for some $w\in q\cup r$ and
\begin{equation}\label{QC22}
\dxh (w^\ast, v^\ast )\le \frac{1}{2}\mu .
\end{equation}
Combining (\ref{QC21}) and (\ref{QC22}), we obtain
$$
dist _Y(u,w)\le c_1\dxh (u^\ast , w^\ast )+c_2\le c_1\left(
H+\frac{1}{2}\mu \right) +c_2.
$$
Thus $\Gamma (R,Y)$ is $\delta ^\prime $--hyperbolic for $\delta
^\prime =c_1\left( H+1/2\mu \right) +c_2$.
\end{proof}

\begin{rem} \label{qcrem1}
The above theorem does not hold without the assumption $card\,
(R\cap H_i^g)< \infty $. Indeed, let $G=H_1\ast H_2$ for some
finitely generated groups $H_1$, $H_2$. Let $K_1$ and $K_2$ be
subgroups of $H_1$ and $H_2$ respectively. Suppose, in addition,
that at least one of the subgroups $K_1,K_2$ is not finitely
generated. Then the subgroup $R$ generated by $K_1, K_2$ is
obviously relatively quasi--convex, but not finitely generated by
the Grushko--Neumann theorem since $R\cong K_1\ast K_2$.
\end{rem}

We conclude with a proposition describing the intersections of
relatively quasi--convex subgroups. The logical scheme of the
proof is due to Short \cite{Sh90}.

\begin{prop}
Let $P$ and $R$ be two strongly relatively quasi--convex subgroups of $G$.
Then $P\cap R$ is strongly relatively quasi--convex.
\end{prop}

\begin{figure}
\unitlength 0.9mm 
\begin{picture}(138,98)(5,-4)
\multiput(9.02,32)(.0603142536,-.0337373737){891}{\line(1,0){.0603142536}}
\multiput(15.03,13.79)(.23233333,.03366667){210}{\line(1,0){.23233333}}
\multiput(89.27,74.78)(.23233333,.03366667){210}{\line(1,0){.23233333}}
\multiput(36.06,16.62)(.0410115448,.03372182518){1819}{\line(1,0){.0410115448}}
\multiput(8.92,31.96)(.041128943,.0337299391){1807}{\line(1,0){.041128943}}
\multiput(62.73,1.93)(.04096079514,.03373826615){1811}{\line(1,0){.04096079514}}
\multiput(63.62,20.81)(.04096079514,.03373826615){1811}{\line(1,0){.04096079514}}
\multiput(89.27,74.85)(-.03968,-.032788){19}{\line(-1,0){.03968}}
\multiput(87.76,73.6)(-.03968,-.032788){19}{\line(-1,0){.03968}}
\multiput(86.25,72.36)(-.03968,-.032788){19}{\line(-1,0){.03968}}
\multiput(84.75,71.11)(-.03968,-.032788){19}{\line(-1,0){.03968}}
\multiput(83.24,69.87)(-.03968,-.032788){19}{\line(-1,0){.03968}}
\multiput(81.73,68.62)(-.03968,-.032788){19}{\line(-1,0){.03968}}
\multiput(80.22,67.37)(-.03968,-.032788){19}{\line(-1,0){.03968}}
\multiput(78.71,66.13)(-.03968,-.032788){19}{\line(-1,0){.03968}}
\multiput(77.21,64.88)(-.03968,-.032788){19}{\line(-1,0){.03968}}
\multiput(75.7,63.64)(-.03968,-.032788){19}{\line(-1,0){.03968}}
\multiput(74.19,62.39)(-.03968,-.032788){19}{\line(-1,0){.03968}}
\multiput(72.68,61.14)(-.03968,-.032788){19}{\line(-1,0){.03968}}
\multiput(71.18,59.9)(-.03968,-.032788){19}{\line(-1,0){.03968}}
\multiput(69.67,58.65)(-.03968,-.032788){19}{\line(-1,0){.03968}}
\multiput(68.16,57.41)(-.03968,-.032788){19}{\line(-1,0){.03968}}
\multiput(66.65,56.16)(-.03968,-.032788){19}{\line(-1,0){.03968}}
\multiput(65.14,54.91)(-.03968,-.032788){19}{\line(-1,0){.03968}}
\multiput(63.64,53.67)(-.03968,-.032788){19}{\line(-1,0){.03968}}
\multiput(62.13,52.42)(-.03968,-.032788){19}{\line(-1,0){.03968}}
\multiput(60.62,51.18)(-.03968,-.032788){19}{\line(-1,0){.03968}}
\multiput(59.11,49.93)(-.03968,-.032788){19}{\line(-1,0){.03968}}
\multiput(57.6,48.68)(-.03968,-.032788){19}{\line(-1,0){.03968}}
\multiput(56.1,47.44)(-.03968,-.032788){19}{\line(-1,0){.03968}}
\multiput(54.59,46.19)(-.03968,-.032788){19}{\line(-1,0){.03968}}
\multiput(53.08,44.95)(-.03968,-.032788){19}{\line(-1,0){.03968}}
\multiput(51.57,43.7)(-.03968,-.032788){19}{\line(-1,0){.03968}}
\multiput(50.07,42.45)(-.03968,-.032788){19}{\line(-1,0){.03968}}
\multiput(48.56,41.21)(-.03968,-.032788){19}{\line(-1,0){.03968}}
\multiput(47.05,39.96)(-.03968,-.032788){19}{\line(-1,0){.03968}}
\multiput(45.54,38.72)(-.03968,-.032788){19}{\line(-1,0){.03968}}
\multiput(44.03,37.47)(-.03968,-.032788){19}{\line(-1,0){.03968}}
\multiput(42.53,36.23)(-.03968,-.032788){19}{\line(-1,0){.03968}}
\multiput(41.02,34.98)(-.03968,-.032788){19}{\line(-1,0){.03968}}
\multiput(39.51,33.73)(-.03968,-.032788){19}{\line(-1,0){.03968}}
\multiput(38,32.49)(-.03968,-.032788){19}{\line(-1,0){.03968}}
\multiput(36.49,31.24)(-.03968,-.032788){19}{\line(-1,0){.03968}}
\multiput(34.99,30)(-.03968,-.032788){19}{\line(-1,0){.03968}}
\multiput(33.48,28.75)(-.03968,-.032788){19}{\line(-1,0){.03968}}
\multiput(31.97,27.5)(-.03968,-.032788){19}{\line(-1,0){.03968}}
\multiput(30.46,26.26)(-.03968,-.032788){19}{\line(-1,0){.03968}}
\multiput(28.96,25.01)(-.03968,-.032788){19}{\line(-1,0){.03968}}
\multiput(27.45,23.77)(-.03968,-.032788){19}{\line(-1,0){.03968}}
\multiput(110.23,77.67)(.05575,-.031583){15}{\line(1,0){.05575}}
\multiput(111.9,76.72)(.05575,-.031583){15}{\line(1,0){.05575}}
\multiput(113.57,75.77)(.05575,-.031583){15}{\line(1,0){.05575}}
\multiput(115.25,74.83)(.05575,-.031583){15}{\line(1,0){.05575}}
\multiput(116.92,73.88)(.05575,-.031583){15}{\line(1,0){.05575}}
\multiput(118.59,72.93)(.05575,-.031583){15}{\line(1,0){.05575}}
\multiput(120.26,71.98)(.05575,-.031583){15}{\line(1,0){.05575}}
\multiput(121.94,71.04)(.05575,-.031583){15}{\line(1,0){.05575}}
\multiput(83.1,92.76)(.060313901,-.03367713){446}{\line(1,0){.060313901}}
\multiput(136.76,62.88)(-.05986111,.0337037){216}{\line(-1,0){.05986111}}
\multiput(15.01,13.68)(.040415094,.033622642){265}{\line(1,0){.040415094}}
\qbezier[1000](44,22.89)(20.96,113.94)(104.95,73.28)
\qbezier(39.39,48.91)(70.24,84.43)(97.22,67.04)
\qbezier(49.05,46.97)(47.72,52.92)(49.95,59.46)
\qbezier(67.04,58.27)(65.7,64.22)(67.93,70.76)
\qbezier(49.95,59.46)(65.85,55.6)(72.84,35.08)
\qbezier(67.93,70.76)(83.84,66.89)(90.83,46.38)
\put(43.85,23.04){\circle*{1.1}} \put(49.2,46.68){\circle*{1.1}}
\put(49.79,59.46){\circle*{1.1}} \put(67.93,70.76){\circle*{1.1}}
\put(97.21,66.9){\circle*{1.1}} \put(104.94,73.29){\circle*{1.1}}
\put(72.84,35.08){\circle*{1.1}} \put(90.68,46.53){\circle*{1.1}}
\put(67.04,58.42){\circle*{1.1}} \put(39.24,48.91){\circle*{1.1}}
\put(44,53.96){\vector(1,1){.07}}\multiput(43.26,53.22)(.032315,.032315){23}{\line(0,1){.032315}}
\put(79.83,72.84){\vector(1,0){.07}}\multiput(77.89,72.69)(.38649,.02973){5}{\line(1,0){.38649}}
\put(48.76,81.91){\vector(1,1){.07}}\multiput(47.57,80.87)(.038361,.033566){31}{\line(1,0){.038361}}
\put(45.5,21){\makebox(0,0)[cc]{$1$}}
\put(106.14,71.2){\makebox(0,0)[cc]{$g$}}
\put(46.97,83.24){\makebox(0,0)[cc]{$p$}}
\put(49.05,61.24){\makebox(0,0)[cc]{$a$}}
\put(67.19,73.14){\makebox(0,0)[cc]{$b$}}
\put(64.07,58.12){\makebox(0,0)[cc]{$bx$}}
\put(70,33.89){\makebox(0,0)[cc]{$ay$}}
\put(88,45.64){\makebox(0,0)[cc]{$by$}}
\put(46.5,45.49){\makebox(0,0)[cc]{$ax$}}
\put(98.11,65){\makebox(0,0)[cc]{$v$}}
\put(36.87,48.46){\makebox(0,0)[cc]{$u$}}
\put(123.68,75.51){\makebox(0,0)[cc]{$R$}}
\put(126.35,60.8){\makebox(0,0)[cc]{$P$}}
\put(41.92,55.2){\makebox(0,0)[cc]{$s_1$}}
\put(79.08,74.77){\makebox(0,0)[cc]{$s_2$}}
\end{picture}
  \caption{}\label{PintQ}
\end{figure}

\begin{proof}

Let us take an element $g\in P\cap R$ and consider a geodesic $p$
in $\G $ such that $p_-=1$, $p_+=g$. We also take an arbitrary
vertex $u\in p$. Let $\sigma $ denote the quasi--convexity
constant for the both subgroups $P$ and $R$.

We have to show that there exists a vertex $v\in P\cap R$ within
an $X$--distance at most $\sigma ^\prime $ from $v$, where $\sigma
^\prime $ is the constant which is independent of $p$ and $u$. Let
us consider a vertex $v\in P\cap R$ satisfying the following two
conditions.

\begin{enumerate}
\item There exists a path $s$ from $u$ to $v$ such that for every
vertex $w\in s$ we have $$ \max \{ \dx (w,P), \; \dx (w, R)\} \le
\sigma .$$ (Note that the set of vertices $v$ satisfying this
condition is non--empty; in particular, it contains $g$.)

\item $v$ is the closest vertex to $u$ with respect to $\dxh $
satisfying the first condition.
\end{enumerate}
We are going to prove that
\begin{equation} \label{QC31}
\dxh (u,v)\le \left( card\; B_\sigma ^X\right) ^2,
\end{equation}
where $B_\sigma ^X$ is the ball defined by (\ref{QC10}).

Suppose that (\ref{QC31}) is false. Then $l(s)> \left( card\;
B_\sigma ^X\right) ^2$. According to the first condition there
exist two vertices, say $a$ and $b$, on $s$ and two elements from
$B_\sigma ^X$, say $x$ and $y$, such that $ax\in P$, $ay\in R$ and
$bx\in P$, $by\in R$ (see Fig. \ref{PintQ}). Let $s=s_1[a,b]s_2$.
Consider the path $t$ in $\G $ such that $t_-=u$ and $\phi
(t)\equiv \phi (s_1)\phi(s_2)$. We state that $t_+\in P\cap R$.
Indeed $bx\in P$ and $b\overline{\phi (s_2)}=v\in P$ yield
$x^{-1}\overline{\phi (s_2)}\in P$. Therefore,
$$t_+=\overline{\phi (s_1)}\overline{\phi (s_2)}=
(ax)(x^{-1}\overline{\phi (s_2)})\in P.$$ Similarly, $t_+\in R$.
Thus $t_+\in P\cap R$. Moreover, for every vertex $w\in t$ we have
\begin{equation}\label{PQsigma}
\max \{ \dx (w,P), \; \dx (w, R)\} \le \sigma .
\end{equation}
Indeed let $t=s_1t_2$, $\phi (t_2)\equiv \phi (s_2)$. If $w\in
s_1$, then the fulfilment of (\ref{PQsigma}) is obvious. Suppose
that $w\in t_2$. Let $z$ be the vertex of $s_2$ such that the
segment $[w,t_+]$ of $t_2$ has the same label as the segment
$[z,v]$ of $s_2$. Then there exists $z^\prime \in P$ such that
$\dx (z, z^\prime )\le \sigma $. Note that
$$
\begin{array}{rl}
w(z^{-1}z^\prime )= & w(z^{-1}v)(v^{-1}z^\prime )=w\overline{\phi
([z,v])}v^{-1}z^\prime =w\overline{\phi ([w,t_+])}v^{-1}z^\prime
=\\&\\& ww^{-1}t_+v^{-1}z^\prime =t_+v^{-1}z^\prime \in P
\end{array}
$$
as $t_+,v,z^\prime \in P$. Thus $$\dx (w,P)\le \dx (w,
w(z^{-1}z^\prime))=\dx (z,z^\prime )\le \sigma .$$ Similarly $\dx
(w, R)\le \sigma $. We have proved that $t_+$ satisfies the first
condition for $v$. Since $t$ is shorter than $s$, we arrive at a
contradiction.

Now we want to estimate the $X$--distance between $u$ and $v$. Let
$r$ be an element of $R$ such that
\begin{equation}\label{QC32}
\dx (u,r)\le \sigma .
\end{equation}
Note that
$$
\dxh (v,r)\le \dxh (u,v)+\dxh (u,r)\le \sigma + \left( card\;
B_\sigma ^X\right) ^2.
$$
Let $Y$ be a finite generating set for $R$. Then, according to
Theorem \ref{qc1}, we have
$$
dist_Y (v,r)\le c_1\left( \sigma + \left( card\; B_\sigma
^X\right) ^2\right) +c_2
$$
where $c_1,c_2$ depend on $R$ and $G$ only.  Therefore,
\begin{equation}\label{QC33}
\dx (v,r)\le dist_Y(v,r) \max\limits_{y\in Y} |y|_X \le \left(
c_1\left( \sigma + \left( card\; B_\sigma ^X\right) ^2\right)
+c_2\right) \max\limits_{y\in Y} |y|_X.
\end{equation}
Summing (\ref{QC32}) and (\ref{QC33}), we obtain
$$
\dx (u,v)\le \left( c_1\left( \sigma + \left( card\; B_\sigma
^X\right) ^2\right) +c_2\right) \max\limits_{y\in Y} |y|_X +\sigma
.
$$
Thus $P\cap R$ is relatively $\sigma ^\prime $--quasi--convex for
$$\sigma ^\prime =\left( c_1\left( \sigma + \left( card\; B_\sigma
^X\right) ^2\right) +c_2\right) \max\limits_{y\in Y} |y|_X +\sigma
.$$
The fact that $P\cap R$ is strongly relatively quasi-convex is clear.
\end{proof}


\section{Cyclic subgroups and translation numbers}


It is well known that if a group $G$ is hyperbolic, then any
cyclic subgroup of $G$ is quasi-isometrically embedded into $G$
\cite{Gr1,GhH}. This result is one of the corner stones of the
small cancellation theory over hyperbolic groups \cite{OlsIJAC}.

It seems to be easy to prove the following by using arguments
similar to the ordinary hyperbolic case. Let $G$ be a group, $\{
H_1, \ldots , H_m\} $ a collection of subgroups of $G$. Suppose
that the relative Cayley graph $\G $ is hyperbolic. Then for any
element $g$ that is not conjugate to an element of one of the
subgroups $H_1, \ldots , H_m$, the cyclic subgroup generated by
$g$ is quasi--isometrically embedded into $G$ (with respect to the
relative metric on $G$).

Unfortunately, in general, this is not true. Indeed, consider an
arbitrary groups $H$ and the direct product $G=H_1\times H_2$,
where $H_i\cong H$, $i=1,2$. Then the relative Cayley graph
$\Gamma (G, \mathcal H)$ of $G$ with respect to $\{ H_1, H_2\} $
has finite diameter and, in particular, it is hyperbolic. To each
element $h\in H$, one can assign the element $g=(h,h)\in H\times
H$. Obviously $g$ is not conjugate to an element of one of the
copies of $H$ in $G$ whenever $h\ne 1$. However, the relative
length of $g^n$ is at most $2$ for every $n$.

In this section we establish the quasi--convexity of cyclic
subgroups generated by hyperbolic elements in relatively
hyperbolic groups. As the previous example shows, to this end we
need some additional arguments apart from the hyperbolicity of $\G
$. Moreover, we do not restrict ourselves to the case of finitely
generated groups, as the general case is important for the
development of the small cancellation theory over relatively
hyperbolic groups \cite{Osin} and some of its applications.
However the Splitting Theorem allows to reduce the proof to the
finitely generated case. After such a reduction, we will follow
the logical scheme suggested in \cite{GS} (the underlying idea has
also been used in \cite{Ep-etal} and \cite{AloB} to obtain similar
results in wider contexts). Our main tools will be the theorems
about quasi--convex subgroups and results from Section 3.4.

\begin{thm}\label{cyc1}
Let $G$ be a finitely generated group hyperbolic relative to a
collection of subgroups $\{ H_1, \ldots , H_m\} $, $g$ a
hyperbolic element of $G$. Then the centralizer $C(g)$ of $g$ in
$G$ is a strongly relatively quasi--convex subgroup in $G$.
\end{thm}

\begin{proof}

\begin{figure}
\unitlength 0.6mm 

\begin{picture}(175.75,80)(5,40)
\qbezier(41.5,40)(147.25,114)(173,64)
\qbezier(41.5,71.25)(147.25,145.25)(173,95.25)
\qbezier(133.25,83.5)(144.25,71.75)(107.25,51)
\qbezier(133.25,114.75)(144.25,103)(107.25,82.25)
\put(133.25,84){\circle*{1.41}}
\put(133.25,115.25){\circle*{1.41}}
\put(107.25,51.25){\circle*{1.41}}
\put(107.25,82.5){\circle*{1.41}} \put(172.5,64.5){\circle*{1.41}}
\put(172.5,95.75){\circle*{1.41}} \put(41.5,40.25){\circle*{1.41}}
\put(41.5,71.5){\circle*{1.41}}

\put(41.5,71.5){\line(0,-1){31.5}}

\put(107,82.75){\line(0,-1){31.25}}

\multiput(132.75,115.25)(.03125,-1.4375){8}{\line(0,-1){1.4375}}

\put(133,99.5){\line(0,-1){15.25}}

\put(172.5,95.5){\line(0,-1){30.75}}

\put(76,62.25){\vector(4,3){.07}}

\put(65,86.75){\vector(4,3){.07}}

\put(125.5,94.25){\vector(-4,-3){.07}}

\put(127,64.25){\vector(-1,-1){.07}}

\put(39.25,38){\makebox(0,0)[cc]{$1$}}
\put(78,59){\makebox(0,0)[cc]{$p$}}
\put(67,91){\makebox(0,0)[cc]{$q$}}
\put(175.5,62){\makebox(0,0)[cc]{$a$}}
\put(136.25,88){\makebox(0,0)[cc]{$v$}}
\put(107.75,48.25){\makebox(0,0)[cc]{$u=vt$}}
\put(37.5,72.75){\makebox(0,0)[cc]{$g$}}
\put(185,96.5){\makebox(0,0)[cc]{$ag=ga$}}
\put(134.5,119.25){\makebox(0,0)[cc]{$w=gv$}}
\put(69,83){$vtg=wt=gvt$}
\end{picture}

\caption{}\label{cent-fig}
\end{figure}

Let $a$ be an element of  $C(g)$, $p$ a geodesic path in $\G $
such that $p_-=1,$ $p_+=a$. We have to show that for any vertex
$v$ on $p$ there exists a vertex $u\in C(g)$ such that $\dx
(u,v)\le \sigma $, where $\sigma =\sigma (g)$ is independent of
$v$ and $a$. For this purpose, we also consider a geodesic $q$ in
$\G $  such that $\phi (p)\equiv \phi (q)$ and $q_-=g$ (see Fig.
\ref{cent-fig}). Since $a\in C(g),$ we have
$$
q_+=g\overline{\phi (q)}=g\overline{\phi (p)}=ga=ag.
$$
By $k$ we denote the length $|g|_X$. Then, $$\dx (p_+, q_+)=\dx
(a, ag)=|g|_X.$$ Thus $(p,q)$ is a pair of $k$--similar symmetric
geodesics in $\G $.

Let $\varkappa =\varkappa (k)$ be the constant from Theorem
\ref{sc}, $w$ the vertex on $q$ synchronous to $v$. By Theorem
\ref{sc}, $\dx (v,w)\le \varkappa $. Let $Conj $ denote the set of
all pair of elements $(f,g)\in G\times G$ such that $f$ and $g$
are conjugate in $G$. To each such a pair we assign an element
$t_{f,g}\in G$ such that $f^{t_{f,g}}=g$ and set
$$
\sigma = \max \{ |t_{f,g}|_X \; : \; (f,g)\in Conj, \;
|f|_X+|g|_X\le k+\varkappa \} .
$$
Since $G$ is locally finite with respect to the metric $\dx $,
$\sigma $ is well--defined and depends on $k$ and $\varkappa $
only (the collection of elements $t_{f,g}$ is supposed to be
fixed).

By the choice of $\sigma $, there exists an element $t\in G$ of
$X$--length at most $\sigma $ such that $(v^{-1}w)^t=g$. Note also
that $w=g\overline{\phi ([g,w])}=gv$. This yields $$
vtg=vt(v^{-1}w)^t= v(v^{-1}w)t=wt=gvt.$$ Hence $vt\in C(g)$. It
remains to note that $\dx (v,vt)=|t|_X\le \sigma $.

Let us show that $C(g)\cap H_i^f$ is finite for any $f\in G$,
$i=1, \ldots m$. Evidently we have $C(g)\cap H_i^f\le H_i^{fg}$.
Therefore,
$$
C(g)\cap H_i^f \le H_i^{fg}\cap H_i^f= \left( H_i^{fgf^{-1}} \cap
H_i\right) ^f.
$$
Since $g$ is hyperbolic, $fgf^{-1}\notin H_i$ and thus the
intersection $H_i^{fgf^{-1}} \cap H_i$ is finite by Proposition
\ref{malnorm}. The lemma is proved.
\end{proof}

\begin{cor}\label{cycfg}
Let $G$ be a finitely generated group hyperbolic relative to a
collection of subgroups $\{ H_1, \ldots , H_m\} $, $g\in G$ a
hyperbolic element of infinite order. Then there exist $\lambda
>0 , c\ge 0$ such that
\begin{equation}\label{cycfg-0}
| g^n| _{X\mathcal \cup H} \ge \lambda |n| -c
\end{equation}
for any $n\in \mathbb Z$.
\end{cor}

\begin{proof}
By Lemma \ref{cyc1}, $C(g)$ is strongly relatively quasi--convex.
Further according to Theorems \ref{qc1} and \ref{qc-hyp}, $C(g)$
is generated by a finite set $Y$ and hyperbolic. The center $Z$ of
$C(g)$ is infinite, as it contains $\langle g\rangle $. As is
well--known, any hyperbolic group with infinite center is
virtually cyclic. Hence the index of $\langle g\rangle $ in $C(g)$
is finite. This obviously implies that the map $(\langle g\rangle
, dist)\to (C(g), dist_Y)$ is a quasi--isometric embedding (here
$dist $ denotes the natural metric on $\langle g\rangle $ with
respect to the generating set $\{ g\} $).

Thus we have a sequence of quasi--isometric embeddings $$(\langle
g\rangle , dist)\to (C(g), dist_Y)\to (G, \dxh ).$$ Since
composition is also a quasi--isometric embedding, we obtain
(\ref{cycfg-0}).
\end{proof}

\begin{cor}
Let $g$ be a hyperbolic element of infinite order in $G$. If
$(g^k)^t=g^l$ for some $k,l\in \mathbb Z$, $t\in G$, then $k=\pm
l$.
\end{cor}

\begin{proof}
The argument is standard. If $k\ne \pm l$, we can assume
$|k|<|l|$. Then $(g^{k^n})^{t^n}=g^{l^n}$. This yields
$$
\begin{array}{rl}
|g^{l^n}|_{X\cup \mathcal H}= & |(g^{k^n})^{t^n}|_{X\cup \mathcal
H} \le 2n|t|_{X\cup \mathcal H}+|k|^n|g|_{X\cup \mathcal H}=\\
& \\ & 2n|t|_{X\cup \mathcal H}+(l^n)^{ \log\limits_l k}
|g|_{X\cup \mathcal H}
\end{array}
$$
for any $n\in \mathbb N$ contradictory (\ref{cycfg-0}).
\end{proof}

The reader can also derive

\begin{cor}
Suppose that $B$ is a subgroup of $G$. If $B$ is isomorphic to a
Baumslag--Solitar group, that is, $$B\cong \langle a,b\; | \;
(a^k)^b=a^l\rangle ,$$ then $B$ is conjugate to a subgroup of
$H_\lambda $ for some $\lambda \in \Lambda $.
\end{cor}

Our next goal is to show that, in fact, the constant $\lambda $ in
Corollary \ref{cycfg} is independent of $g$. It is convenient to
express this property in terms of translation numbers.

\begin{defn}
Let $K$ be a group generated by a finite set $S$ relative to a
collection of subgroups $\{ I_\lambda \} _{\lambda \in \Lambda }$,
$x$ an element of $K$. The {\it relative translation number} of
$x$ is defined to be
\begin{equation}\label{trn1}
\tau ^{rel}(x)=\lim\limits_{n\to \infty }\frac{1}{n} |x^n|_{S\cup
\mathcal I},
\end{equation}
where $|\cdot |_{S\cup \mathcal I}$ denotes the relative length
function with respect to $S$ and $\{ I_\lambda \} _{\lambda \in
\Lambda }$.
\end{defn}

The following lemma is quite trivial (The proof of its
non--relative analog can be found in \cite{GS}; it works in the
relative case without any changes.)

\begin{lem}

\begin{enumerate}
\item The limit in (\ref{trn1}) always exists.

\item $\tau ^{rel}(x)=\inf\limits_n \frac{|x|_{S\cup \mathcal
I}}{n} $. In particular, $\tau^{rel}(x)\le |x|_{S\cup \mathcal
I}.$

\item $\tau ^{rel}(x)=\tau ^{rel}(x^t)$ for any $x,t\in K$; thus
$\tau ^rel(x) $ depends only on the conjugacy class of $x$.

\item $\tau ^{rel} (x^n)= |n|\tau ^{rel}(x) $ for any $x\in K$,
$n\in \mathbb Z$.
\end{enumerate}
\end{lem}

\begin{thm}\label{trd00}
Let $G$ be a group, $\{ H_\lambda \} _{\lambda \in \Lambda }$ a
collection of subgroups of $G$. Suppose that $G$ is hyperbolic
with respect to $\{ H_\lambda \} _{\lambda \in \Lambda }$. Then
there exists $d >0$ such that for every hyperbolic element of
infinite order $g\in G$, the relative translation number of $g$
with respect to $\{ H_\lambda \} _{\lambda\in \Lambda }$ satisfies
the inequality $\tau ^{rel}(g)>d $.
\end{thm}

\begin{proof}
The proof consists of several lemmas. First of all we are going to
reduce the proof to the case when $G$ is finitely generated.

\begin{lem}\label{tgtn}
Let $K$ be the fundamental of a tree of groups with vertex groups
$I_1, \ldots , I_m$. Suppose that $x$ is an element of $K$. Then
either $x$ is conjugate to an element of $I_i$ for some $i$ or the
translation number of $g$ relative to $I_1, \ldots , I_m$ is at
least $2$.
\end{lem}

\begin{proof}
Obviously it suffices to prove the lemma in case $m=2$. (Then we
can apply inductive arguments). If $m=2$, $K$ is a free product of
$I_1,I_2$ with amalgamated subgroups, say $A$ and $B$. Below we
use terminology and some well--known results about amalgamated
products, which can be found, for example, in \cite[Sec.
4.2]{MKS}. If $x$ is not conjugate to an element of $I_1$ or
$I_2$, then $x$ is conjugate to a cyclically reduced (in the sense
of amalgamated products) element $y=ay_1\ldots y_d$, where $d\ge
2$, $a\in A$, $y_1,\ldots , y_d$ are coset representatives of
$I_1$ or $I_2$ with respect to $A$, and $y_i,y_{i+1}$ (indices are
$mod\; d$) are not in the same factor. Then the number of factors
in the reduced form of $y^n$ is at least $2dn\ge 2n$. Since the
reduced form is unique and the reduction process does not increase
the number of factors, we have $|y^n|_\mathcal I\ge 2|n|$. This
yields the assertion of the lemma as translation numbers depend
only on conjugacy classes.
\end{proof}

Recall that by the Splitting Theorem, $G=G_0\ast\left(\ast
_{\lambda \in \Lambda\setminus \Lambda _0} H_\lambda \right) $,
where $G_0$ is a tree of groups with vertex groups $H_\lambda,
\lambda \in \Lambda _0=\{ \lambda _1, \ldots , \lambda _m\} $, and
$Q$. If $g\in G$ is hyperbolic, then either $g$ has length at
least $2$ (as the element of a free product), or $g$ is conjugate
to an element of $G_0$. In the first case $\tau^{rel}(g)\ge 2$ by
the obvious reasons. In the second case Lemma \ref{tgtn} shows
that it suffices to prove the theorem for elements $g$ such that
$g$ is conjugate to an element of $Q$. Since translation numbers
depend only on the conjugacy classes, we may assume $g\in Q$.
Recall that $Q$ is finitely generated and hyperbolic relative to
some subgroups $L_i$, $i=1, \ldots , m$. By Proposition
\ref{QisomG}, it suffices to show that there exists $d>0$ such
that for any $g\in Q$ the relative translation number of $g$ in
$Q$ with respect to $\{ L_1, \ldots , L_m\} $ is greater than $d$.

Thus we can assume that the group $G$ is generated by a finite set
$X$ in the usual (non--relative) sense. To prove the theorem it
suffices to show that for some $a >0 $, the set
\begin{equation}\label{trda}
T_G(a)=\{ \tau ^{rel}(g)\; | \; g\in G, \; g \;{\rm is\;
hyperbolic}, \; \tau ^{rel} (g) < a \}
\end{equation}
is finite. To this end we use an auxiliary lemma below, which can
be regarded as a generalization of Lemma \ref{FinOrd2}. Throughout
the rest of this section, we use the notation $M$, $L$ (see the
beginning of the previous chapter) and denote by $\delta $ the
hyperbolicity constant of $\G $.

\begin{lem}\label{trd}
Suppose that $g$ is a hyperbolic element of $G$ satisfying the
following conditions.

\begin{enumerate}
\item $g$ has minimal relative length among all elements of the
conjugacy class $g^G$.

\item $|g|_{X\cup \mathcal H} \le 8\delta +1$.

\item $\tau ^{rel} (g)< \xi $ for $\xi =\frac{1}{8\delta +3}$.
\end{enumerate}

\noindent Then the $X$--length of $G$ satisfies
\begin{equation}\label{trd-main}
|g|_X \le (8\delta +1)(32\delta +6)ML.
\end{equation}
\end{lem}

\begin{proof}
Since $\tau ^{rel} (g)< \xi $, there exists $N\in \mathbb N$ such
that $|g^N|_{X\cup \mathcal H}< \xi N$. Let $U\in (X\cup \mathcal
H)^\ast $ be a shortest word representing $g$ in $G$ and $W\in
(X\cup \mathcal H)^\ast $ a shortest word representing $g^N$.
Obviously we have
\begin{equation}\label{trd1}
\| U\| \le 8\delta +1
\end{equation}
and
\begin{equation}\label{trd2}
\| W\| \le \xi N.
\end{equation}

We consider the cycle $pq^{-1} $ in $\G $ such that $p_-=q_-=1$,
$\phi (p)\equiv U^N$, and $\phi (q)\equiv W$. There are three
possibilities to consider.

{\it Case 1}. First suppose that any component of $pq^{-1}$ is
isolated. Given an $H_i$--syllable $V$ of $U$, we have at least
$N$ $H_i$--components of $p$ labelled $V$. By Lemma \ref{31}, we
have
$$
N|\overline{V}|_X\le MLl(pq^{-1})\le ML(N\| U\| +\| W\| ).
$$
Using (\ref{trd1}), (\ref{trd2}) and dividing both the sides of
the previous inequality by $N$, we obtain
\begin{equation}\label{trd0}
|\overline{V}|_X\le ML\frac{N(8\delta +1) +\xi N}{N}< ML(8\delta
+2)
\end{equation}
(note that $\xi <1$). Since (\ref{trd0}) is true for every
syllable of $U$, we obtain
$$
|g|_X=|\overline{U}|_X\le \| U\|ML(8\delta +2)\le (8\delta
+1)(8\delta +2)ML.
$$
Obviously this inequality is even stronger than (\ref{trd-main}).

{\it Case 2}. Assume that there are two connected components $s_1$
and $s_2$ of $p$. Then repeating the same arguments as in the
proof of Lemma \ref{FinOrd2}, we obtain
$$
|g|_X\le 2ML(8\delta +1)
$$
and the inequality (\ref{trd1}) obviously holds. We leave details
to the reader.

{\it Case 3}. Suppose that no different components of $p$ are
connected. As the path $q$ is geodesic, no different components of
$q$ are connected. Therefore, the only possibility for two
components $s$ and $t$ of $pq^{-1}$ to be connected is $s\in p$,
$t\in q$.

Let $s_1, \ldots , s_n$ be the set of all components of $p$ such
that for any $i=1, \ldots , n$, there exists a component of $q$,
denoted by $t_i$, that is connected to $s_i$. Let us denote by
$e_i$ and $f_i$ the paths of length at most $1$ in $\G $ such that
$(e_i)_-=(t_i)_-$, $(e_i)_+=(s_i)_-$, $(f_i)_-=(t_i)_+$,
$(f_i)_+=(s_i)_+$. We also set $e_0$, $f_0$ (respectively
$e_{n+1}$ and $f_{n+1}$) to be the paths consisting of just one
vertex $1$ (respectively $g^N$). Without loss of generality we may
assume that
$$
p=r_1s_1\ldots r_ns_nr_{n+1}.
$$

Let us choose a subsequence $s_{i_1}, \ldots , s_{i_l}$ of the
sequence $s_1, \ldots , s_n $ as follows. We set $s_{i_1}=s_1$ .
Further suppose we have already chosen $s_{i_k}$. Then
$s_{i_{k+1}}$ is defined to be the first component in the sequence
$s_{i_k+1}, s_{i_k+2},\ldots , s_n$ such that the corresponding
component $t_{i_{k+1}}$ belongs to the segment $[(t_{i_k})_+,
q_+]$ of $q$. Thus after completing this process, we will obtain a
sequence of components $s_{i_1}, \ldots , s_{i_l}$ of $p$ such
that (see Fig. \ref{trfig}):

(i) $ p=p_1s_{i_1}\ldots p_ls_{i_l}p_{l+1} $ for some $p_1, \ldots
p_{l+1}$;

(ii) $ q=q_1t_{i_1}\ldots q_lt_{i_l}q_{l+1}$ for some $q_1, \ldots
, q_{l+1}$;

(iii)for any $j=1, \ldots l+1$, every component of the cycle
$c_j=f_{i_{j-1}}p_je_{i_j}^{-1}q_{i_j}^{-1}$ is isolated in $c_j$
(Figure \ref{trfig}).

For simplicity, we change the notation and denote $s_{i_j}$,
$f_{i_j}$, and $e_{i_j}$ by $s_j$, $f_j$, and $e_j$ respectively.
Thus we have
$$
 p=p_1s_{1}\ldots p_ls_{l}p_{l+1},
$$
$$
q=q_1t_{1}\ldots q_lt_{l}q_{l+1},
$$
and
$$
c_j=f_{j-1}p_je_{j}^{-1}q_{j}^{-1}.
$$

\begin{figure}

\unitlength 0.8mm 

\begin{picture}(140,80)(18,50)

\put(29.01,64.96){\line(1,0){123.402}}

\put(29.01,65.17){\circle*{.94}} \put(29.22,85.14){\circle*{.94}}
\put(29.01,102.8){\circle*{.94}} \put(56.97,125.92){\circle*{.94}}
\put(81.57,125.71){\circle*{.94}}
\put(152.41,100.91){\circle*{.94}}
\put(152.56,86.82){\circle*{.94}}
\put(139.38,65.17){\circle*{.94}}
\put(126.13,65.17){\circle*{.94}}
\put(152.12,65.17){\circle*{.94}} \put(93.55,65.17){\circle*{.94}}
\put(73.16,65.17){\circle*{.94}} \put(55.71,65.17){\circle*{.94}}
\put(43.1,65.17){\circle*{.94}}

\put(29.05,74.95){\vector(0,1){.07}}
\put(29.05,112.43){\vector(0,1){.07}}
\put(152.5,76.19){\vector(0,-1){.07}}
\put(146.9,65){\vector(1,0){.07}}
\put(65.94,65){\vector(1,0){.07}}
\put(37.83,65){\vector(1,0){.07}}

\put(90.77,88.3){\oval(123.57,75.13)[t]}

\multiput(29.1,85.21)(-.02946,-3.32927){6}{\line(0,-1){3.32927}}
\multiput(152.59,100.94)(-.02946,-5.95143){6}{\line(0,-1){5.95143}}

\put(71.24,125.86){\vector(1,0){.07}}
\put(82.55,65){\vector(1,0){.07}}
\put(133.82,65){\vector(1,0){.07}}
\put(152.5,93.69){\vector(0,-1){.07}}
\put(50.03,65){\vector(1,0){.07}}
\put(29.17,94.93){\vector(0,1){.07}}

\put(104.12,94.93){\circle*{.5}} \put(109.6,94.93){\circle*{.5}}
\put(115.08,94.93){\circle*{.5}}

\thicklines

\put(93.55,65.1){\line(-1,0){20.2}}
\put(139.38,65.1){\line(-1,0){13.25}}
\put(152.41,100.91){\line(0,-1){14}}

\put(56.74,125.7){\line(1,0){24.749}}
\put(43.13,65.05){\line(1,0){12.551}}
\put(29.03,102.87){\line(0,-1){17.541}}

\put(71.24,125.86){\vector(1,0){.07}}
\put(82.55,65){\vector(1,0){.07}}
\put(133.82,65){\vector(1,0){.07}}
\put(152.5,93.69){\vector(0,-1){.07}}
\put(50.03,65){\vector(1,0){.07}}
\put(29.17,94.93){\vector(0,1){.07}}

\put(68.59,98.82){\vector(-1,4){.07}}
\put(89.98,99.17){\vector(-1,4){.07}}
\put(48.79,88.39){\vector(-3,4){.07}}
\put(40.13,78.84){\vector(-2,3){.07}}
\put(128.1,82.38){\vector(1,3){.07}}
\put(140.15,77.07){\vector(1,4){.07}}

\qbezier(29.17,85.21)(43.04,84.32)(43.13,65.05)
\qbezier(55.68,65.05)(52.77,96.7)(28.99,102.88)
\qbezier(56.92,126.04)(72.48,100.23)(73.18,65.23)
\qbezier(81.32,126.04)(90.77,113.14)(93.51,64.88)
\qbezier(139.12,65.05)(138.95,90.86)(152.2,86.97)
\qbezier(126.22,65.05)(125.16,100.14)(152.38,100.94)

\put(26.34,64.52){\makebox(0,0)[cc]{$1$}}
\put(26.7,76.37){\makebox(0,0)[cc]{$p_1$}}
\put(25.99,95.64){\makebox(0,0)[cc]{$s_{i_1}$}}
\put(26.52,113.14){\makebox(0,0)[cc]{$p_2$}}
\put(68.06,129.22){\makebox(0,0)[cc]{$s_{i_2}$}}
\put(156.5,96.7){\makebox(0,0)[cc]{$s_{i_l}$}}
\put(158,79){\makebox(0,0)[cc]{$p_{l+1}$}}
\put(157,67){\makebox(0,0)[cc]{$g^N$}}
\put(36.24,68){\makebox(0,0)[cc]{$q_1$}}
\put(48.79,68){\makebox(0,0)[cc]{$t_{i_1}$}}
\put(64.35,68){\makebox(0,0)[cc]{$q_2$}}
\put(81.32,68){\makebox(0,0)[cc]{$t_{i_2}$}}
\put(132.76,68){\makebox(0,0)[cc]{$t_{i_l}$}}
\put(144.96,68){\makebox(0,0)[cc]{$q_{l+1}$}}
\put(40.5,83){\makebox(0,0)[cc]{$e_{i_1}$}}
\put(51,93){\makebox(0,0)[cc]{$f_{i_1}$}}
\put(72,100.23){\makebox(0,0)[cc]{$e_{i_2}$}}
\put(94,103){\makebox(0,0)[cc]{$f_{i_2}$}}
\put(125.7,84.32){\makebox(0,0)[cc]{$e_{i_l}$}}
\put(138.3,81){\makebox(0,0)[cc]{$f_{i_l}$}}
\end{picture}
\caption{}\label{trfig}
\end{figure}

We will call $c_1, \ldots , c_{l+1}$  sections. We state that
there exists a section $c_j$ such that the following two
inequalities hold.
\begin{equation} \label{trd5}
\left\{
\begin{array}{l}
l(q_j)\le l(p_j),\\ \\
l(p_j)\ge 8\delta +1.
\end{array}
\right.
\end{equation}

Indeed let us denote by $S_1$ the set of all sections $c_j$ such
that $l(q_j)>l(p_j)$ and by $S_2$ the set of all sections $c_j$
such that $l(p_j)< 8\delta +1$. It suffices to show that there
exists a section $c_j\notin S_1\cup S_2$. Notice that the number
of all sections is at most $l(q)< \xi N$. According to the choice
of $U$, any component $s_j$ consists of a single edge. Let $\pi
(c_j)$ denote the number $l(p_j)+l(s_j)=l(p_j)+1$. On one hand, we
have
$$\sum\limits_{j=1}^{l+1} \pi (c_j)=l(p)+1.$$ On the other hand,
\begin{equation}\label{trd6}
\sum\limits_{c\in S_2} \pi (c)< (8\delta +2)card\, S_2< (8\delta
+2)\xi N
\end{equation}
and
\begin{equation}\label{trd7}
\sum\limits_{c\in S_1} \pi (c)< l(q)+1< \xi N +1
\end{equation}
Inequalities (\ref{trd6}) and (\ref{trd7}) yield
$$
\sum\limits_{c\in S_1\cup S_2} \pi (c)< (8\delta +3)\xi N+1= N+1
\le l(p)+1.
$$
Therefore there is a section which is not in $S_1\cup S_2$.

Let $c_j$ be a section satisfying (\ref{trd5}). We have
$$
\phi (p_j)\equiv U_0^kA,
$$
where $U_0$ is a cyclic shift of the word $U$ and $\| A\| < \| U\|
$. According to the second inequality of (\ref{trd5}) and the
second condition of the lemma, we have $k\ge 1$. Let $V$ be a
syllable of $U$. Then arguing as in Case 1, and taking into
account the inequality
$$
l(c_j)\le  l(p_j)+l(q_j)+2\le 2l(p_j)+2< 2(k+1)\| U\| +2\le
2(k+1)(8\delta +1)+2
$$
we obtain
$$
k|\overline{V}|_X\le MLl(c_j)\le  ML(2(k+1)(8\delta +1)+2 ).
$$
Hence
$$
|\overline{V}|_X\le ML(32\delta +6).
$$
Since this is true for every syllable of $U$, we obtain
$$
|g|_X=|\overline{U}|_X\le \| U\|ML(32\delta +6)\le (8\delta
+1)(32\delta +6)ML.
$$
\end{proof}

Let us return to the proof of Theorem \ref{trd00}.  Given an
arbitrary element $g$ in $G$ which has shortest relative length in
the conjugacy class $g^G$, there are two possibilities.

First suppose that $|g|_{X\cup \mathcal H}\ge 8\delta +1$. We take
a shortest word $U\in (X\cup\mathcal H)^\ast $ representing $g$
and consider the path $p_n$ such that $(p_n)_-=1$, $\phi
(p_n)\equiv U^n$, $n\in \mathbb N$. Since $g$ is a shortest
element in $g^G$, $p_n$ is $(8\delta +1)$--local geodesic in $\G
$. Therefore, by Lemma \ref{k-loc} we have
$$
|g^n|_{X\cup \mathcal H}=\dxh (1, (p_n)_+) \ge \frac{1}{\lambda }
(l(p_n) -c)= \frac{1}{\lambda } (n |g|_{X\cup \mathcal H} -c),
$$
for $\lambda \le 3$ and $c=2\delta $. Hence, $$\tau ^{rel} (g) \ge
\lim\limits_{n\to \infty } \frac{1/\lambda (n|g|_{X\cup \mathcal
H} -c)}{n} =\frac{1}{\lambda } |g|_{X\cup\mathcal H}\ge 1/3.$$

Now assume that $|g|_{X\cup \mathcal H}<8\delta +1 $. Then, by
Lemma \ref{trd}, either $\tau ^{rel}(g)\ge \xi =1/ (8\delta +3)$
or $g$ is conjugate to an element of the set $B=B_{(8\delta
+1)(32\delta +6) ML}^X$ of elements whose $X$--length is at most
$(8\delta +1)(32\delta +6) ML$. Since $B$ is finite and $\tau
^{rel} (g)$ depends only on the conjugacy class of $g$, $T_G(a)$
is finite for $a= 1/ (8\delta +3)$.
\end{proof}


\chapter{Algorithmic problems}



\section{The word and membership problems}


Recall that the word problem for a recursively presented group
generated by a recursive set $X$ is to decide, given a word $W$ in
the alphabet $X^{\pm1}$, whether $W$ represents $1$ in $G$. In
\cite{F}, Farb showed that the word problem is solvable for any
finitely generated group $G$ hyperbolic relative to subgroups $\{
H_1, \ldots , H_m\}$ provided it is solvable for each of the
subgroups $H_1, \ldots , H_m$. It is not hard to generalize this
result as follows.

\begin{thm}
Suppose that $G$ is a group given by a finite relative
presentation with respect to recursively presented subgroups $H_1,
\ldots , H_m$. Assume also that the corresponding relative Dehn
function $\delta ^{rel}(n)$ is bounded from above by some
recursive function and the word problem is solvable for all
subgroups $H_i$, $i=1, \ldots ,m$. Then the word problem is
solvable for $G$.
\end{thm}

\begin{proof}
Let $X$ be a finite generating set of $G$, $X=X^{-1}$, $\widetilde
H_i$ an isomorphic copy of $H_i$, $i=1, \ldots , m$,
$$F=F(X)\ast \widetilde H_1 \ast \ldots \ast \widetilde H_m,$$ and
let $N$ be the kernel of the naturally defined homomorphism $F\to
G$. Suppose that $N=\langle \mathcal R\rangle ^F$ for a certain
finite subset $\mathcal R\subset N$, and the corresponding
relative Dehn function $\delta ^{rel}(n)$ is bounded from above by
some recursive function. We also assume that the relative
presentation
\begin{equation}
\langle X, H_1, \ldots , H_m\; |\; R=1,\; R\in \mathcal R\rangle
\label{WordProb1}
\end{equation}
is reduced (see Definition \ref{redpres}). Let $\Omega _i$,
$\Omega $ denote the (finite) sets introduced in Definition
\ref{defofOmega}. By Proposition \ref{Hfg}, $\Omega _i$ generates
$H_i$. Thus the group $F$ is generated by the finite set $Z=X\cup
\widetilde\Omega $, where $\widetilde \Omega $ consists of the
preimages of elements of $\Omega $ under the canonical
homomorphisms $\widetilde H_i\to H_i$.

For any word $W$ over $X$, $W\in N$, of length $\|W\| \le n$, we
consider a van Kampen diagram $\Delta $ of minimal type over
(\ref{WordProb1}) with boundary label $W$. Note that any edge of
this diagram is labelled by an element of $Z$. Indeed this is so
for any external edge since $W$ is a word over $X$. If $e$ is an
internal edge, then $e$ belongs to the boundary of some $\mathcal
R$--cell by Lemma \ref{DiagMinType} and thus $\phi (e)\in Z$.
Fixing a basepoint in $\Delta $, we obtain a representation
\begin{equation}\label{WP2}
W=_F\prod\limits_{i=1}^k f_i^{-1}R_if_i,
\end{equation}
$f_i\in F$, $R_i\in \mathcal R$, $i=1, \ldots ,k$, where $k\le
\delta ^{rel}(n)$, and each element $f_i$ is a label of a path
without self--intersections in $\Delta $. In particular, $|f_i|_Z$
is not greater than the number of all edges in $\Delta $, i.e.,
$$|f_i|_Z\le M\delta ^{rel}(n) +\|W\| \le M\delta ^{rel} (n)+n$$
(recall that $M=\max\limits_{R\in \mathcal R} \| R\| $.) Thus the
number of factors in (\ref{WP2}) and the lengths of conjugating
elements $f_i$ with respect to the (finite) set $Z$ are bounded by
recursive functions of $n$. Since the word problem is solvable for
$F$, these bounds allow to derive that it is solvable for $G$.
\end{proof}

It is worth to notice that one can easily provide an example of a
finitely generated group $G$ and a finitely generated subgroup $H$
of $G$ such that the word problem is solvable for both $G$ and
$H$, but the corresponding relative Dehn function is not
well--defined. (For example, this is so if $G=H\times \mathbb Z$
and $H=\mathbb Z$.)

Recall that the {\it membership} problem for a subgroup $K$ of a
group $G$ is to decide for a given element $g\in G$ whether $g$
belongs to $K$.

\begin{thm}\label{GenWP}
Suppose that $G$ is a group given by a finite relative
presentation with respect to recursively presented subgroups $H_1,
\ldots , H_m$. Assume also that the corresponding relative Dehn
function $\delta ^{rel}(n)$ is bounded from above by some
recursive function and the word problem is solvable for all
subgroups $H_i$, $i=1, \ldots ,m$. Then the membership problem is
solvable for $H_i$, $i=1, \ldots , m$.
\end{thm}

Note that the theorem does not hold without the requirement of the
solvability of the word problem in $H_i$, $i=1, \ldots , m$. For
example, if $H$ is a finitely presented group with undecidable
word problem, then the group $G=H$ is hyperbolic relative to the
subgroups $H_1=H$ and $H_2=\{ 1\} $. Clearly the membership
problem for $H_2$ is unsolvable in this case. To prove the theorem
we need an auxiliary definition.

\begin{defn}
Let $H$ be a group generated by a finite set $A$ and $K$ a
subgroup of $H$ generated by a finite set $B$. The {\it distortion
function} of $K$ in $H$ with respect to the generating sets $A$
and $B$ is defined to be
$$
\Delta _K^H(n)=\max \{  |x|_B\; |\;  x\; {\rm satisfies}\;
|x|_A\le n\} .
$$
It is easy to see that $\Delta _K^H$ is independent (up to
equivalence) of the choice of finite generating sets in $H$ and
$K$. In case $\Delta _K^H(n)\sim n$, we say that $K$ is
undistorted in $H$.
\end{defn}

\begin{lem}\label{distortion}
Let $G$ be a finitely generated group, $\{ H_1, \ldots H_m\} $ a
collection of finitely generated subgroups of $G$. Suppose that
$G$ is finitely presented with respect to $\{ H_1, \ldots H_m\} $
and the corresponding relative Dehn function $\delta ^{rel}(n)$ is
well--defined. Then for any $i=1, \ldots , m$, the distortion of
the subgroup $H_i$ in $G$ satisfies
\begin{equation}\label{distr1}
\Delta _{H_i}^G (n)\preceq \delta ^{rel}(n).
\end{equation}
\end{lem}

\begin{proof}
Let $h$ be a non--trivial element of $H_i$, $W$ a shortest word
over $X$ representing $h$. Then there exists a cycle $c=pq^{-1}$
in $\G $ such that $\phi (p)\equiv W$ and $q$ is an edge in $\G $
labelled $h$. Note that $q$ is an isolated $H_i$--component of
$c$. Applying Lemma \ref{Omega} we obtain
$$
|h|_{\Omega _i}=|\overline{\phi (q)}|_{\Omega _i}\le M
Area^{rel}(c)\le M\delta ^{rel} (\| W\| +1)=M\delta ^{rel} (|h|_X
+1).
$$
This implies (\ref{distr1}).
\end{proof}

\begin{proof}[Proof of Theorem \ref{GenWP}]
As is well-known \cite{F0}, if the word problem in a finitely
generated group $H$ is solvable and a subgroup $K\le H$ is
finitely generated, then the membership problem for $K$ is
solvable if and only if $\Delta _K^H(n)$ is bounded from above by
a recursive function. Thus Theorem \ref{GenWP} follows from the
previous lemma.
\end{proof}

\begin{cor}\label{WP}
Let $G$ be a finitely generated group hyperbolic relative to
recursively presented subgroups $\{ H_1, \ldots , H_m\}$. Suppose
that the word problem is solvable in each of the subgroups $H_1,
\ldots , H_m$. Then:

1) (Farb, \cite{F}) The word problem is solvable in $G$.

2) For any $i=1, \ldots , m$, the membership problem is solvable
for $H_i$.
\end{cor}


\section{The parabolicity problems}


In the previous chapter we saw that some important properties
(such as the finiteness of the conjugacy classes of elements of
finite orders, strong quasi--convexity of cyclic subgroups, etc.)
hold for hyperbolic elements although they can be violated for
parabolic ones.  Thus given a finitely generated group $G$ which
is hyperbolic relative to a collection of subgroups $\{ H_1,
\ldots , H_m\} $, it is natural to consider the following two
algorithmic problems

1) ({\it The general parabolicity problem}) Given an element $g\in
G$, decide whether $g$ is parabolic or hyperbolic.

2) ({\it The special parabolicity problem}) Given an element $g\in
G$ and $i\in \{ 1, \ldots , m\} $, decide whether $g$ is conjugate
to an element of $H_i$.

In case $m=1$ these problems coincide. It is proved in \cite{Bum}
that they are solvable for any finitely generated group $G$
hyperbolic relative to a subgroup $H$ whenever the conjugacy
problem is decidable in $H$. Similar arguments allow to show that
the special parabolicity problem is solvable for any finitely
generated group $G$ hyperbolic relative to $\{ H_1, \ldots , H_m\}
$ whenever the conjugacy problem is solvable in all $H_1, \ldots ,
H_m $. We observe that the last requirement is essential for the
special parabolicity problem. On the other hand, the general
parabolicity problem is solvable in $G$ whenever the word problem
is solvable in $H_1, \ldots , H_m $. In this section we give the
proof of these results in the spirit of our paper.

\begin{thm}\label{ParProb}
Let $G$ be a group hyperbolic relative to a collection of
recursively presented subgroups $\{ H_1, \ldots , H_m\} $.

1) If the word problem is solvable for all $H_i$, $i=1, \ldots ,
m$, then the general parabolicity problem is solvable in $G$.
Moreover, there is an algorithm which allows, given a parabolic
element $g\in G$, to find some $t\in G$ and some $j\in \{ 1,
\ldots ,m\} $ such that $g^t\in H_j$.

2) (Bumagin, \cite{Bum}) If the conjugacy problem is solvable for
all $H_i$, $i=1, \ldots , m$, then the special parabolicity
problem is solvable in $G$. Moreover, there is an algorithm which
allows, given $i\in \{ 1, \ldots , m\} $ and $g\in G$ that is
conjugate to an element of $H_i$, to find an element $t\in G$ such
that $g^t\in H_i$.
\end{thm}

The proof of the theorem is based on the next two lemmas.

\begin{lem}\label{PP1}
Let $G$ be a finitely generated group hyperbolic relative to a
collection of subgroups $\{ H_1, \ldots , H_m\}$. There exists a
recursive function $\sigma (k)$ satisfying the following
condition. Let $g$ be a parabolic element of $G$ such that
$|g|_X\le k$. Then there exists $t\in G$ such that $g^t\in H_j$
for a certain $j\in \{ 1, \ldots , m\} $ and $|t|_X\le \sigma
(k)$.
\end{lem}

\begin{proof}
Let $\mathcal P$ be the set of all pairs of symmetric geodesics in
$\G $ with characteristic elements $g,h$, where $h\in H_j$ for
some $j=1, \ldots , n$ (see Section \ref{SecSymGeod} for
definitions). Let $(p,q)$ be a pair of geodesics of minimal
lengths in $\mathcal P$. By Corollary \ref{sk3}, the length of the
element $t=\overline{\phi (p)}=\overline{\phi (q)}$ is not greater
than $\rho (k)$, where the constant $\rho (k)$ can be effectively
calculated (one can notify that this is the common property of all
constants in our paper). We note that no synchronous components of
$p$ and $q$ are connected. Indeed if $p=p_1ap_2$, $q=q_1bq_2$,
where $a,b$ are connected synchronous components of $p$ and $q$
respectively, then $(p_1,q_1)\in \mathcal P$, which contradicts to
the minimality of length of $p$ and $q$. Therefore, by Lemma
\ref{minpair} there are no connected components of $p$ and $q$ at
all.

Further let $r$ be the edge in $\G $ labelled by an element $h\in
H_j$ for some $j$ such that $r_-=p_+$, $r_+=q_+$. Note that $r$
can not be connected to an $H_j$--component of $p$ or $q$. Indeed
if $p=s_1cs_2$, where $c$ is an $H_j$--component of $p$ connected
to $r$, then the label of $cs_2$ represents an element of $H_j$.
Thus, for $t_0=\overline{\phi (s_1)}$, we have
$$g^{t_0}=g^{t\overline{\phi (cs_2)}^{-1}}= h^{\overline{\phi
(cs_2)}^{-1}}\in H_i .$$ This contradicts to the choice of $(p,q)$
again. Hence no components of the paths $pr$ and $q$ are
connected. Obviously $pr$ and $q$ are $k$--connected
$(1,2)$--quasi--geodesics. By Theorem \ref{TBCP}, any component
$s$ of $p$ has X--length at most $|s|_X\le \varepsilon (1,2,k) $,
where $\varepsilon (1,2,k)$ can be effectively calculated for
given $k$. Thus we have $|t|_X\le \sigma (k)$ for $\sigma (k)=\rho
(k)\varepsilon (1,2,k)$.
\end{proof}

\begin{lem}\label{PP2}
Let $G$, $\{ H_1, \ldots , H_m\}$ be as in the previous lemma.
Suppose, in addition, that $H_1, \ldots , H_m$ are recursively
presented and have solvable conjugacy problem. Then there exists a
recursive function $\theta (k)$ satisfying the following
condition. Let $g$ be an element of $G$ conjugate to an element
$h\in H_i$ for some $i\in \{ 1, \ldots , m\} $ such that $|g|_X\le
k$. Then there exists $t\in G$ such that $g^t\in H_i$ and
$|t|_X\le \theta (k)$.
\end{lem}

\begin{proof}
Let $\mathcal Q$ be the set of all pairs of symmetric geodesics in
$\G $ with characteristic elements $g,h_0$, where $h_0\in H_i$.
Let $(p,q)$ be a pair of geodesics of minimal lengths in $\mathcal
Q$. As above, using Corollary \ref{sk3}, we obtain
$|t|_{X\cup\mathcal H}\le \rho (k)$, where $t=\overline{\phi
(p)}=\overline{\phi (q)}$. By Lemma \ref{minpair}, only
synchronous components of $p$ and $q$ can be connected.

By $r$ we denote the edge in $\G $ going from $p_+$ to $q_+$ and
labelled by the element $h_0\in H_i$, which is conjugate to $g$ by
$t$. The paths $pr$ and $q$ are $k$--similar
$(1,2)$--quasi--geodesics. As in the previous lemma, we can show
that $r$ can not be connected to a component of $p$ or $q$.

If there are no connected synchronous components of $p$ and $q$,
we can repeat the arguments from the proof of the previous lemma.
Further let $p=p_1a_1\ldots p_na_np_{n+1}$, $q=q_1b_1\ldots
q_nb_nq_{n+1}$, where for $j=1, \ldots , n$, $a_j, b_j$ are
connected synchronous components of $p$ and $q$ respectively. We
assume that $p_j$ and $q_j$ contain no connected components for
$j=1, \ldots , n+1$.

We denote by $v_j$ and $w_j$ the elements $((a_j)_-)^{-1}(b_j)_-$
and $((a_j)_+)^{-1}(b_j)_+$ of $G$ respectively.  By Theorem
\ref{TBCP}, for any $j=1, \ldots , n$, we have $\max \{ |v_j|_X,\,
|w_j|_X\} \le \varepsilon (1,2,k)$. Since the conjugacy  problem
is decidable in $H_1, \ldots , H_m$, there is a recursive function
$\tau :\mathbb N\to \mathbb N$ such that for every $i\in \{1,
\ldots , m\} $, and any two elements $v,w\in H_i$ that are
conjugate in $H_i$, there is an element $s\in H_i$ such that
$v^s=w$ and $|s|_X\le \tau (\max \{ |v|_X,\, |w|_X\} )$. Note that
elements $v_j$ and $w_j$ are conjugate in the subgroup $H_{i_j}$
for a suitable $i_j$. Let $s_j\in H_{i_j}$ be the corresponding
conjugating element such that $v_j^{s_j}=w_j$ and $|s_j|_X\le \tau
(\varepsilon (1,2,k))$. It is easy to check that for the element
$z=\overline{\phi (p_1)}s_1\ldots \overline{\phi
(p_n)}s_n\overline{\phi (p_{n+1})}$, we have $g^z=g^t=h_0$.

Let us estimate the length of $z$. By our assumption, no
components of subpaths $p_j$ and $q_j$, $j=1, \ldots , n$ are
connected. The same is true for the subpaths $p_{n+1}r$ and $q_n$.
Since for any $j=1, \ldots , n$, $p_j$ and $q_j$ are $l=\max\{
\varepsilon(1,2,k),\, k\} $--similar, as well as $p_{n+1}r$ and
$q_{n+1}$,  for any component $a$ of $p_j$, $j=1, \ldots, n+1$, we
have $|a|_X\le \varepsilon (1,2,l)$ by Theorem \ref{TBCP}. Hence,
$$
|t|_X=|\overline {\phi (p)}|_X\le \rho (k) \max\{ \varepsilon
(1,2,l),\, \tau (\varepsilon (1,2,k))\} .
$$
\end{proof}

\begin{proof}[Proof of Theorem \ref{ParProb}]
Let us prove the first statement of the theorem. By Corollary
\ref{WP}, the word problem is decidable in $G$. Moreover, the
membership problem is decidable for $H_i$, $i=1, \ldots , m$.
Given an element $g\in G$, we solve the membership problem for all
elements of type $g^t$, where $|t|_X\le \sigma (|g|_X)$, and $H_j$
for all $j=1, \ldots , m$. By Lemma \ref{PP1}, $g$ is parabolic if
and only if $g^t\in H_j$ for some $|t|_X\le \sigma (|g|_X)$ and
$j\in \{ 1, \ldots , m\} $.

The proof of the second statement follows from Lemma \ref{PP2} in
the same way.
\end{proof}

Finally we note that the solvability of the conjugacy problem in
each $H_i$, $i=1, \ldots , m$, is essential in the second
statement of Theorem \ref{ParProb}. That is, if we replace this
condition with the solvability of the word problem in each $H_i$,
$i=1, \ldots , m$, the statement would be false. Clearly the
minimal possible value of $m$ in any counterexample is $2$.

\begin{thm}\label{PP3}
There exists a group $G$ hyperbolic relative to finitely presented
subgroups $H_1, H_2$ such that the word problem is solvable in
$H_1$ and $H_2$ and the special parabolicity problem is unsolvable
in $G$.
\end{thm}

\begin{proof} Suppose that $H$ is a finitely generated group having
solvable word problem, $x\in H$ a fixed nontrivial element of
order $2$ such that there is no algorithm which allows to decide
whether $h\sim x$ for a given element $h\in H$. (We explain how to
construct such a group below.) Let $G=H_1\ast _{\langle x\rangle }
H_2$ be the amalgamated product of two copies of $H$, where the
amalgamated subgroups are generated by the elements corresponding
to $x$ in each copy of $H$. Note that $G$ is hyperbolic relative
to $\{ H_1, H_2\} $. Indeed the existence of the action of $G$ on
the Bass--Serre tree yields that $G$ is hyperbolic relative to $\{
H_1, H_2\} $ in the sense of Bowditch and thus in our sense (see
Appendix). Obviously an element $h\in H_1$ is conjugate to an
element of $H_2$ if and only if $h$ is conjugate to $x$ in $H_1$.
Thus the special parabolicity problem is unsolvable in $G$.

To construct the group $H$ with the desired properties, let us
take the abelian group $$A=\langle  x_i, i\in \mathbb N \; |\;
x_i^2=1,\, [x_i,x_j]=1,\, i,j\in \mathbb N\rangle .$$ Let $f:
\mathbb N\to \mathbb N $ be a recursive function such that the
range of $f$ is not recursive. Set $N=\{x_{f(i)},\, i\in \mathbb
N\} $. We fix an arbitrary element $x\in N$ and consider the
sequence of groups $Q(n)$, $n=0,1,\ldots $, such that $Q(0)=A$ and
$Q(n)$ is obtained from $Q(n-1)$ by adding one extra generator
$t_{n}$ subject to the relation
$$t_{n}^{-1}x_{f(n)}t_{n}=x.$$

\begin{lem}
For any non--negative integer $n$, the word problem and in $Q(n)$
is decidable.
\end{lem}

\begin{proof}
The case $n=0$ is obvious. Let $n>0$. The group $Q(n)$ is the
HNN--extension of $Q(n-1)$ with finite associated subgroups
$\langle x_{f(n)}\rangle $ and $\langle x \rangle $. Since the
subgroups $\langle x_{f(n)}\rangle $ and $\langle x \rangle $ are
finite and the word problem in $Q(n-1)$ is solvable, the
membership problem for $\langle x_{f(n)}\rangle $ and $\langle x
\rangle $ is solvable in $Q(n-1)$. Hence the word problem is
solvable in $Q(n)$ \cite[Corollary 2.2, Ch. IV]{LS}.
\end{proof}

We set $Q=\bigcup\limits_{i=1}^\infty Q(n)$. Clearly $Q$ is
recursively presented and the word problem is solvable in $Q$.
\begin{lem}
For an element $a\in A$, $a$ is conjugate to $x$ in $Q$ if and
only if $a\in N$.
\end{lem}

\begin{proof}
The 'if' part follows from our construction. We now suppose that
$a\in A$ and $a$ is conjugate to $x$ in $Q$. Then $a$ is conjugate
to $x$ in $Q(n)$ for some $n$. If $n=0$, then $a=x\in N$ since
$Q(0)=A$ is abelian. If $n>0$, without loss of generality we may
assume that $a$ is not conjugate to $x$ in $Q(n-1)$. Then
$a=x_{f(n)}$ by the Collins Lemma \cite[Theorem 2.5, Ch. IV]{LS}.
\end{proof}

Finally we embed $Q$ into a finitely presented group $H$ such that
the word problem is decidable in $H$ and two elements of $Q$ are
conjugate in $H$ if and only if they are conjugate in $Q$. (Such
an embedding exists by the Olshanskii--Sapir Theorem, see
\cite{OS}.)

It follows from the construction that an element $x_i$ is
conjugate to $x$ in $H$ if and only $i$ belongs to the range of
$f$. As the range of $f$ is not recursive, the problem of whether
a given element $h\in H$ is conjugate to $x$ in $H$ is unsolvable.
\end{proof}


\section{Algorithmic problems for hyperbolic elements}


In this section we assume that $G$ is a finitely generated group
hyperbolic relative to a collection $\{ H_1, \ldots , H_m \}$ of
subgroups. We also use some notation introduced at the beginning
of Chapter 3. For two elements $f,g\in G$ we write $f\sim g$ if
$f$ is conjugate to $g$.

Recall that the conjugacy problem for a group $G$ given by a
recursive presentation is to decide, for any two elements $g_1$
and $g_2$ of $G$, whether or not $g_1$ is conjugate to $g_2$ in
$G$. Recently Bumagin \cite{Bum} proved the following.

\begin{thm}[Bumagin, \cite{Bum}]
Suppose that a finitely generated group $G$ is hyperbolic relative
to recursively presented subgroups $\{ H_1, \ldots , H_m\}$ and
the conjugacy problem is solvable for each of the subgroups $H_1,
\ldots , H_m$. Then the conjugacy problem is solvable for $G$.
\end{thm}

By the {\it conjugacy problem for hyperbolic elements} we mean the
following: given two hyperbolic elements $f,g\in G$, decide
whether $f\sim g$. Observe that the next theorem is true without
any assumptions about the conjugacy problem in $H_1, \ldots ,
H_m$.

\begin{thm}\label{ConjProb}
Suppose that the word problem is solvable for all $H_i$, $i=1,
\ldots , m$. Then the conjugacy problem for hyperbolic elements is
solvable in $G$.
\end{thm}

\begin{proof}
Since the word problem is solvable in $G$ by Corollary \ref{WP} ,
it is enough to show that there is a recursive function $\alpha
:\mathbb N\to \mathbb N$ such that for any two conjugate
hyperbolic elements $f,g\in G$ satisfying the inequality $\max\{
|f|_X, \, |g|_X\} \le k$, there exists an element $t\in G$ such
that $|t|_X\le \alpha (k)$ and $f^t=g$.

Let $(p,q)$ be a minimal pair of symmetric geodesics in $\G $ with
characteristic elements $f,g$. Note that no synchronous components
of $p$ and $q$ are connected. Indeed otherwise $f$ and $g$ are
parabolic.  The rest of the proof almost coincide with the proof
of Lemma \ref{PP1}. By Lemma \ref{minpair} there are no connected
components of $p$ and $q$ at all. Further by Corollary \ref{sk3}
the length of the element $t=\overline{\phi (p)}=\overline{\phi
(q)}$ is not greater than $\rho (k)$. Note that $p$ and $q$ are
$k$--connected geodesics. By Theorem \ref{TBCP}, any component $s$
of $p$ has X--length at most $|s|_X\le \varepsilon (1,0,k) $. Thus
we have $|t|_X\le \alpha (k)$ for $\alpha (k)=\rho (k)\varepsilon
(1,0,k)$.
\end{proof}

Our results concerning algebraic properties of relatively
hyperbolic groups allow to treat 'relative versions' of some other
algorithmic problems.

\begin{defn}
The {\it order problem for hyperbolic elements} is to calculate
the order of a given hyperbolic element. The {\it root problem for
hyperbolic elements} is to decide whether for a given hyperbolic
element $g\in G$ there exists $f\in G$ and $n\in \mathbb N$,
$n>1$, such that $g=f^n$. The {\it power conjugacy problem for
hyperbolic elements} is to decide whether for given hyperbolic
elements $f,g\in G$ there exist $k,l\in \mathbb Z$, such that
$g^k$ and $f^l$ are hyperbolic and $g^k\sim f^l$.
\end{defn}

For known results about the ordinary order, power conjugacy, and
root problems in various classes of groups and relations between
these problems we refer to
\cite{Bez,Can,Com,Fine,Miller,Lar,Lip1,Lip2,Lys}.

\begin{thm}\label{OP}
Suppose that the word problem is solvable in $H_i$ for any $i=1,
\ldots , m$. Then the order problem for hyperbolic elements is
solvable in $G$.
\end{thm}

\begin{proof}
Given a hyperbolic element $g\in G$, there are only finitely many
possibilities for the order of $g$ by Theorem \ref{FinOrd}. Thus
the order problem in $G$ is reduced to the word problem.
\end{proof}

The next lemma will help us to treat the root problem.

\begin{lem} \label{RP0}
There exists a recursive function $\beta :\mathbb N\to \mathbb N$
satisfying the following condition. Let $g$ be an element of $G$,
$|g|_X=k$, and let $f$ be a hyperbolic element of $g$ such that
$f^n=g$ for some $n\in \mathbb N$. Then $f$ is conjugate to an
element $f_0\in G$ such that $|f_0|_X\le \beta (k).$
\end{lem}

\begin{proof}
Let $f_0$ be the element with minimal relative length in the
conjugacy class $f^G$, $U$ a shortest word in $X\cup\mathcal H$
representing $f_0$. By $t$ we denote an element of $G$ such that
$f_0=f^t$. Let $T$, $S$ be the shortest words in $X$ representing
$t$ and $g$ respectively. For any $j\in \mathbb N$, there is a
cycle $c_j=p_jq_j^{-1}$ in $\G $ such that $\phi (p_j)\equiv
U^{nj}$, $\phi (q_j)\equiv T^{-1}S^jT$, as $f_0^{nj}=(g^j)^t$.

If $|f_0|_{X\cup\mathcal H}\ge 8\delta +1$, where $\delta $ is the
hyperbolicity constant of $\G $, then $p_j$ is $(3, 2\delta
)$--quasi--geodesic by Lemma \ref{k-loc}. We obtain
$$|t^{-1}g^jt|_{X\cup\mathcal H} = \dxh ((p_j)_-, (p_j)_+)
\ge \frac{1}{3} (l(p_j)-2\delta )\ge \frac{j}{3}
|f_0|_{X\cup\mathcal H}-2\delta .$$ Thus $$ j|g|_{X\cup \mathcal
H}+ 2|t|_{X\cup \mathcal H} \ge \frac{j}{3}|f_0|_{X\cup\mathcal
H}-2\delta .$$ Dividing by $j$ and passing to the limit as $j\to
\infty $, we obtain $$|f_0|_{X\cup \mathcal H}\le 3|g|_{X\cup
\mathcal H}\le 2|g|_X\le 3k .$$ Thus in any case we have
\begin{equation}\label{RP1}
|f_0|_{X\cup \mathcal H} \le \max\{  8\delta +1,\, 3k\} .
\end{equation}

Suppose that for some $j\in \mathbb N$, there exist two connected
components of $p_j$, i.e., $p_j=as_1p_0s_2c$, where $s_1,s_2$ are
connected components of $p_j$. Without loss of generality we may
assume that no components of $p_0$ are connected. Note that $f$ is
hyperbolic. Repeating the arguments from the proof of Lemma
\ref{FinOrd2} and using (\ref{RP1}), we obtain the estimate
\begin{equation}\label{RP2}
|f_0|_{X} \le 2ML(\max\{  8\delta +1,\, 3k\} )^2.
\end{equation}

Further assume that for any $j\in \mathbb N$, all components of
$p_j$ are isolated in $p_j$. This yields that all components of
$c_j$ are isolated in $c_j$ since $\phi (q_j)$ is a word in the
alphabet $X$. By Lemma \ref{31}, for any $i=1, \ldots , m$, any
$H_i$--syllable $V$ of $U$ satisfies the inequality
$$
j|n||\overline V|_X \le ML l(c_j) \le ML (j|n||f_0|_{X\cup
\mathcal H} +2|t|_{X}+j|g|_{X}).
$$
Thus
$$
|\overline V|_X \le ML \left( |f_0|_{X\cup \mathcal H}
+\frac{2}{j}|t|_{X}+k\right) .
$$
Assuming $j\to \infty $, we obtain
$$
|\overline V|_X \le ML \left( |f_0|_{X\cup \mathcal H} +k\right)
\le 2ML \max\{  8\delta +1,\, 3k\} .
$$
This implies (\ref{RP2}) again. Thus we can set $$\beta
(k)=2ML(\max\{  8\delta +1,\, 3k\} )^2.$$
\end{proof}

\begin{thm}
Suppose that the word problem is solvable in $H_i$ for any $i=1,
\ldots , m$. Then the root problem for hyperbolic elements is
solvable in $G$.
\end{thm}

\begin{proof}
Let $g$ be a hyperbolic element of $G$. Since the order problem
for hyperbolic elements is solvable in $G$, we can decide whether
the order of $g$ is finite. Let us consider two cases.

1. The order of $g$ is finite. Recall that the set of powers of
hyperbolic elements of $G$ is finite. Let $n_0$ be the maximal
finite order of hyperbolic elements in $G$. If $f^n=g$ for some
positive $n$, then we may assume that $n\le n_0$.  By Lemma
\ref{RP0}, to decide whether $g$ has a non--trivial root in $G$,
it suffices to decide whether there exists an element $f_0\in G$
and a natural number $n\le n_0$ such that $|f_0|_X\le \beta
(|g|_X)$ and $f_0^n\sim g$. Thus the root problem is reduced to
the conjugacy problem for hyperbolic elements, which is solvable
by Theorem \ref{ConjProb}.

2. The order of $g$ is infinite. If $f^n=g$ for some $f\in G$,
$n\in \mathbb N$, then
$$\tau ^{rel} (f)=\frac{1}{n} \tau (g)\le \frac{1}{n} |g|_{X\cup
\mathcal H}\le \frac{1}{n}|g|_X .$$ Thus, by Theorem \ref{trd00},
$n\le |g|_X/d $, where $d$ is a constant which is independent of
$g$ and $f$. Using Lemma \ref{RP0}, we now can reduce the root
problem to the conjugacy problem for hyperbolic elements as in the
first case.
\end{proof}

To deal with the power conjugacy problem, we need the following.

\begin{lem}\label{GCP0}
There exists a recursive function $\gamma :\mathbb N\to \mathbb N$
such that if $f,g\in G$ are two hyperbolic elements of infinite
order and $f^k\sim g^l$ in $G$ for some $k,l\in \mathbb Z\setminus
\{ 0\} $, then there exist $k^\prime ,l^\prime \in \mathbb
Z\setminus \{ 0\} $ such that $f^{k^\prime }\sim g^{l^\prime }$ in
$G$ and $$ \max \{ |k^\prime |, \; |l^\prime |\} \le \gamma (\max
\{ |f|_X, \; |g|_X\} ).$$
\end{lem}

\begin{proof}
Suppose that $k$ and $l$ are numbers with minimal $\max \{ |k|, \;
|l|\} $ among all pairs $k,l$ satisfying $f^k\sim g^l$. We are
going to show that
\begin{equation} \label{GCP1}
\max \{ |k|, \; |l|\} \le \gamma (\max \{ |f|_X, \; |g|_X\} )
\end{equation}
for some recursive function $\gamma $. For definiteness we assume
that $|k|\ge |l|$.

Let $t$ be an element of $G$ such that

1) $g^l=tf^kt^{-1}$ and

2) $t$ has minimal relative length among all elements of $G$
satisfying 1).

\noindent Let $U$, $V$ (respectively $T$) be shortest words over
$X$ (respectively $X\cup \mathcal H$) representing elements $f, g$
(respectively $t$). The conditions of the lemma imply that there
exists a quadrangle
$$
rps^{-1}q^{-1}
$$
in $\G $ such that $\phi (r)\equiv \phi (s)\equiv T$, $\phi
(p)\equiv U^k$, $\phi (q)\equiv V^l$. For convenience, we assume
$r_-=1$.

Let $a_1, \ldots a_k$ denote the ending vertices of the subpaths
$p_1, \ldots , p_k$ of $p$ with $(p_i)_-=p_-$ and $\phi
(p_i)\equiv U^i$, $i=1, \ldots , k$. Let $\lambda >0$ be a number
such that there is no hyperbolic element $h$ of infinite order in
$G$ with $\tau ^{rel}(h)< \lambda $. Note that $p$ and $q$ are
$(\lambda , c)$--quasi--geodesics in $\G $, for $c=\max \{ |f|_X,
\; |g|_X\} $ since $\tau ^{rel} (g)=\inf \{ |g^n|_{X\cup \mathcal
H}/n\} $. By Corollary \ref{quad} and Theorem \ref{TBCP}, for any
vertex $a_i\in p$ there exists a vertex $b_i\in q\cup r\cup s$
such that
\begin{equation} \label{GCP2}
\dx (a_i, b_i)\le 2(\nu +\varepsilon ),
\end{equation}
where $\varepsilon =\varepsilon (\lambda , c,0)$. First of all we
wish to show that if $a_i$ is sufficiently far from the endpoints
of $p$, then $b_i\in q$.

\begin{figure}
\unitlength 0.6mm 
\begin{picture}(130,50)(-23,30)

\put(33.01,41.83){\line(0,1){39.102}}
\put(33.01,80.94){\line(1,0){88.715}}
\multiput(121.72,80.94)(.03003,-5.58595){7}{\line(0,-1){5.58595}}
\put(121.93,41.83){\line(-1,0){88.715}}

\put(33,41.92){\circle*{1.5}} \put(33.15,80.72){\circle*{1.5}}
\put(121.89,80.57){\circle*{1.5}}
\put(121.89,41.92){\circle*{1.5}} \put(71.95,41.77){\circle*{1.5}}
\put(32.85,54.26){\circle*{1.5}}
\qbezier(32.85,54.11)(66.52,55.82)(71.95,41.77)
\put(30.77,42.96){\makebox(0,0)[cc]{}}
\put(30,42.96){\makebox(0,0)[cc]{$t$}}
\put(31.07,81.91){\makebox(0,0)[cc]{$1$}}
\put(140,82){\makebox(0,0)[cc]{$g^l=tf^kt^{-1}$}}
\put(130,41.92){\makebox(0,0)[cc]{$tf^k$}}
\put(72.54,37.5){\makebox(0,0)[cc]{$a_i=tf^i$}}
\put(86.22,84){\makebox(0,0)[cc]{$q$}}
\put(100.34,438.5){\makebox(0,0)[cc]{$p$}}
\put(36.5,58){\makebox(0,0)[cc]{$b_i$}}
\put(36,69.27){\makebox(0,0)[cc]{$r$}}
\put(119.66,62.14){\makebox(0,0)[cc]{$s$}}
\put(48.01,38){\makebox(0,0)[cc]{$p_i$}}
\put(56.93,54.55){\makebox(0,0)[cc]{$o$}}

\put(32.85,70.01){\vector(0,-1){.07}}
\put(45.04,41.92){\vector(1,0){.07}}
\put(98.11,41.92){\vector(1,0){.07}}
\put(54.85,53.07){\vector(4,-1){.07}}
\put(121.89,63.62){\vector(0,-1){.07}}
\put(83.84,80.87){\vector(1,0){.07}}

\end{picture}
 \caption{}\label{glfk}
\end{figure}

For instance, suppose that $b_i\in r$. Let $o$ be a geodesic in
$\G $ such that $o_-=b_i$, $o_+=a_i$ (see Fig. \ref{glfk}).  From
(\ref{GCP2}), we obtain $l(o)\le 2(\nu
 +\varepsilon ). $ Note that the element $a_i$ satisfies
 $g^la_if^ka_i^{-1}$ since $a_i=tf^i$. By the choice of $t$ this
 means that
\begin{equation} \label{GCP3}
|a_i|_{X\cup \mathcal H} \ge |t|_{X\cup \mathcal H}.
\end{equation}
Notice that
\begin{equation} \label{GCP4}
|a_i|_{X\cup \mathcal H}\le |b_i|_{X\cup \mathcal H}+l(o)\le
|b_i|_{X\cup \mathcal H} +2(\nu +\varepsilon)
\end{equation}
and
\begin{equation} \label{GCP5}
|t|_{X\cup \mathcal H}=|b_i|_{X\cup \mathcal H}+l([b_i, t])
\end{equation}
as $r$ is geodesic. Combining (\ref{GCP3}), (\ref{GCP4}), and
(\ref{GCP5}), we obtain the following estimate on the length of
the segment $[b_i, t]$ of $r$:
$$
l([b_i, t])\le 2(\nu +\varepsilon ).
$$
Hence,
$$
\dxh (t,a_i)\le \dxh (t, b_i) +\dxh (b_i, a_i) \le 4(\nu
+\varepsilon ).
$$

Since $p$ is $(\lambda , c)$--quasi--geodesic, we have
$$
l(p_i)\le \lambda\dxh (t, a_i)+c\le 4\lambda (\nu +\varepsilon
)+c.
$$
Therefore, if $i\ge N$ for $$ N=4\lambda (\nu +\varepsilon )+c, $$
then $b_i$ can not belong to $r$ as $l(p_i)=\| U\| i\ge i$.
Similarly, if $i<k-N$, then $b_i$ can not belong to $s$.

Suppose that $k> 2N$. Then for any $N\le i\le k-N$, $b_i\in q$. We
denote by $c_j$ the terminal vertex of the subpath $q_j$ of $q$,
such that $(q_j)_-=q_-=1$ and $\phi (q_j)\equiv V^j$. Thus
$c_j=g^j$. Let $c_{j(i)}$ be the closest vertex (with respect to
the relative metric) to the vertex $b_i$. (Here the index $j(i)$
depends of $i$.) Obviously
$$
\dx (c_{j(i)}, b_i)\le \frac{1}{2}|g|_X.
$$
Hence,
$$
|a_i^{-1}c_{j(i)}|_X=\dx (a_i, c_{j(i)})\le \frac{1}{2}|g|_X+2
(\nu +\varepsilon ).
$$

Let
$$
k_0=2N +(card \; X)^{\frac{1}{2}|g|_X+2 (\nu +\varepsilon )} +1.
$$
If $k>k_0$, then there exists two pairs $a_{i_1}, c_{j(i_1)}$ and
$a_{i_2}, c_{j(i_2)}$ such that $a_{i_1}^{-1}c_{j(i_1)}=
a_{i_2}^{-1}c_{j(i_2)}$ since the number of different elements of
$G$ of $X$--length at most $n$ is less that or equal to $(card\,
X)^n+1$. Let $t_0=a_{i_1}^{-1}c_{j(i_1)}$, $j_0=j(i_1)-j(i_2)$,
and $i_0=i_1-i_2$. Then $t_0^{-1}$ conjugate $g^{j_0}$ to
$f^{i_0}$ or $f^{-i_0}$. Since $|i_0|\le |k|$, we arrive to the
contradiction with the choice of $k$ and $l$. Thus $k$ is less
than or equal to $k_0$ and the inequality (\ref{GCP1}) is true for
$$\gamma =2N +(card \; X)^{\frac{1}{2}|g|_X+2 (\nu +\varepsilon )} +1.$$
\end{proof}

\begin{thm}
Suppose that the word problem is solvable in $H_i$ for any $i=1,
\ldots , m$. Then the power conjugacy problem for hyperbolic
elements is solvable in $G$.
\end{thm}

\begin{proof}
Let $f,g$ be two hyperbolic elements of $G$. There are three cases
to consider (by Theorem \ref{OP} they can be effectively
recognized).

1. Both the elements $f,g$ have finite orders $n_1$ and $n_2$
respectively. Then it suffices to decide whether the elements
$f^k$ and $g^l$ are hyperbolic and conjugate for some $k$ and $l$
satisfying $0<k<n_1$, $0<l<n_2$. This can be done by Theorems
\ref{ParProb} and \ref{ConjProb}.

2. $g$ has finite order, $f$ has infinite order (or conversely).
Clearly $f$ is not conjugate to $g$ in this case.

3. Both the elements $f,g$ have infinite order. Then applying
Lemma \ref{GCP0}, we can reduce the question to the conjugacy
problem for hyperbolic elements in $G$ as in the first case.
\end{proof}

\setcounter{chapter}{5} \setcounter{thm}{0}

\chapter*{Open questions}

\addcontentsline{toc}{chapter}{Open questions}

Here we state some natural problems and conjectures which seem to
be important for the further studying relatively hyperbolic
groups.

If $G$ is a finitely presented group and is hyperbolic relative to
a subgroup $H$, then $H$ is finitely generated by Proposition
\ref{Hfg}. However the following important question remains open.

\begin{prob}
Let $G$ be a finitely presented group hyperbolic relative to a
subgroup $H$. Does it follow that $H$ is finitely presented?
\end{prob}

Assume that $G$ is a group generated by a finite set $X$
hyperbolic relative to a collection of subgroups $H_1, \ldots,
H_m$. As we have already mentioned in Section 4.2, Theorem
\ref{qc-hyp} does not hold without the assumption $card\; (R\cap
H_i^g)<\infty $. However, one can try to prove a similar result in
the general case. More precisely, let $R$ be a relatively
quasi--convex subgroup of $G$. We consider the set of subgroups
$$
\mbox{\eufm S} = \{ H_i^g\cap R\; | \; i=1, \ldots , m, \; g\in
G\} .
$$
The subgroup $R$ acts on {\eufm S} by conjugations. If $R$ is
quasi--convex, it is not hard to show that the number of orbits of
this action is finite. Let $P_1, \ldots , P_l$ be representatives
of the orbits.

\begin{prob}
Prove that the group $R$ is relatively hyperbolic with respect to
$\{ P_1, \ldots , P_l\} $.
\end{prob}

Note that $R$ and $P_i$ need not be finitely generated in this
case. The simplest example of this type is the pair $G=H_1\ast
H_2$ and  $R=\langle K_1, K_2\rangle $ considered in Remark
\ref{qcrem1}. The proof of the conjecture should be slightly more
complicated than one of Theorem \ref{qc-hyp}, although the
difficulties are rather technical.

Another problem about relatively quasi--convex groups is

\begin{prob}
Does the notion of relative quasi--convexity formulated in this
paper coincide with the dynamical quasi--convexity introduced by
Bowditch in \cite{B0}?
\end{prob}

Let $G$ be a group generated by a finite set $X$, $H_1, \ldots ,
H_m$ subgroups of $G$. We consider the following condition:

\begin{enumerate}
\item[($\ast $)] \it There exists a constant $\nu >0$ such that
for any geodesic triangle $\Delta = pqr$ in $\Gamma (G, X\cup
\mathcal H)$ and any vertex $v$ on $p$, there exists a vertex $u$
on the union $q\cup r$ such that $\dx (u,v)\le \nu .$
\end{enumerate}
Theorem \ref{r0} shows that if $G$ is hyperbolic relative to $\{
H_1, \ldots , H_m\} $, then ($\ast $) is satisfied. The converse
is not true. Indeed any group $G$ satisfies $(\ast )$ with respect
to any subgroup $H$ of finite index. However, by Proposition
\ref{malnorm}, $G$ is never hyperbolic relative to $H$ unless
$G=H$ or $H$ is finite. On the other hand, it is easy to see that
any group $G$ satisfying $(\ast )$ is weakly hyperbolic (or
hyperbolic in the sense of Farb) relative to $H_1, \ldots , H_m$.
And again the converse is not true. For example, if $G=H\times
\mathbb Z$, where $H$ is an infinite finitely generated group,
then $G$ is weakly hyperbolic relative to $H$, but do not satisfy
$(\ast )$.  Thus the class of groups satisfying $(\ast )$ is
intermediate between the classes of hyperbolic and weakly
hyperbolic groups.

\begin{prob}
Study the class of finitely generated groups satisfying $(\ast )$
with respect to a finite collection of subgroups.
\end{prob}

The next two problems are inspirited by well--known questions
about ordinary hyperbolic groups.

\begin{prob}\label{Hopf}
Assume that $G$ is a finitely generated group hyperbolic relative
to a collection of subgroups $\{ H_1, \ldots , H_m \} $. Suppose
that all subgroups $H_i $ are Hopfian. Does it follow that $G$ is
Hopfian?
\end{prob}

Recall that a group $G$ is said to be Hopfian if every epimorphism
$G\to G$ is an isomorphism. First examples of non--Hopfian groups
were found by Boumslag and Solitar \cite{BS62}: the group
\begin{equation} \label{BS}
BS(m,n)=\langle a,t\; | \; (a^m)^t=a^n\rangle ,
\end{equation}
is not Hopfian whenever $m$ and $n$ are relatively prime and
neither of $|m|$, $|n|$ is equal to $1$. We note that $BS(m,n)$ is
weakly relatively hyperbolic with respect to $\langle a\rangle $
(see \cite{Wh}). Sela proved that if a torsion--free hyperbolic
group does not decompose as a free product, then it is Hopfian. It
seems to be reasonable to assume that analogous result is true in
the relative case.

\begin{prob}\label{rf}
Assume that $G$ is a group hyperbolic relative to a collection of
subgroups $\{ H_\lambda \} _{\lambda \in \Lambda }$. Suppose that
all subgroups $H_\lambda $ are residually finite. Can $G$ be
non--residually finite?
\end{prob}

This question is open even for ordinary hyperbolic group $G$. It
is known \cite{Gr1,Ol95} that an infinite hyperbolic group is
never simple, but it is still unknown whenever it is always
residually finite (some speculations can also be found in
\cite{IO}, \cite{OlsBL}, and \cite{KW}). One of the results in
this direction can be found in the paper \cite{W02}, where Wise
showed that negatively curved polygons of finite groups are
residually finite. We note that in the case of CAT(0)--groups the
answer to Problem \ref{rf} is known to be positive. In the paper
\cite{W96}, Wise produced examples of compact non--positively
curved spaces whose fundamental groups are not residually finite.
Subsequently, Burger and Mozes \cite{BM} constructed compact
non--positively curved 2--complexes whose fundamental groups are
simple.

\setcounter{chapter}{6} \setcounter{thm}{0}


\chapter*{Appendix. Equivalent definitions of relative
hyperbolicity}


\addcontentsline{toc}{chapter}{Appendix. Equivalent definitions of
relative hyperbolicity}

\subsection*{The definition of Bowditch}

The original definition of Bowditch characterizes relative
hyperbolicity in dynamical terms.

\begin{defn} \label{RGDyn}
A finitely generated group $G$ is {\it hyperbolic relative to a
collection of finitely generated subgroups} $H_1, \ldots , H_m$ if
it admits a properly discontinuous isometric action on a
path--metric hyperbolic proper space $X$ such that the induced
action of $G$ on the boundary $\partial X$ satisfies the following
conditions.

(1) $G$ acts on $\partial X$ as a geometrically finite convergence
group.

(2) The maximal parabolic subgroups of $G$ are precisely the
subgroups of $G$ conjugate to $H_1, \ldots , H_m$.

\end{defn}

Convergence groups were introduced by Gehring and Martin \cite{GM}
in order to describe the dynamical properties of Kleinian groups
acting on the standard sphere in $\mathbb R ^n$ and were
generalized to groups acting on compact Hausdorff spaces by Tukia
and Freden \cite{Tuk}, \cite{Fred}. Their motivation came from the
observation that an isometry group of a hyperbolic space $X$ acts
as a convergence group on the hyperbolic boundary of $X$.

We recall that a group $G$ of homeomorphisms of a metrizable
compactum $M$ acts on $M$ as a {\it convergence group} if the
induced action on the space of distinct triples of elements of $M$
is properly discontinuous (for equivalent definitions we refer to
\cite{B0}). One says that a subgroup $H\le G$ is {\it parabolic}
if it is infinite, fixes some point of $M$ and contains no
elements $g\in G$ of infinite order such that $card (fix(g))=2$.
In this case the fixed point of $H$ is unique and is called a {\it
parabolic point}. A parabolic point $x$ is said to be {\it
bounded} if $\left ( M\setminus \{ x\} \right) /Stab_G(x) $ is
compact. A point $x\in M$ is called a {\it conical limit point} if
there is a sequence $\{ g_i\} $ and two distinct points $a,b\in M$
such that $g_ix$ converges to $a$ and $g_iy$ converges to $b$ for
any $y\in M\setminus \{ x\} $. Finally, a convergence group $G$ is
said to be {\it geometrically finite} if every point of $M$ is a
conical limit point or a bounded parabolic point.

In \cite{B1}, Bowditch proved that Definition \ref{RGDyn} is
equivalent to the following.

\begin{defn}\label{BowDef}
A finitely generated group $G$ is {\it hyperbolic relative to a
collection of finitely generated subgroups} $H_1, \ldots , H_m$ if
it admits an action on a hyperbolic graph $K$ such that the
following condition hold.

1) All edge stabilizers are finite.

2) All vertex stabilizers are finite or conjugate to one of the
subgroups $H_1, \ldots , H_m$.

3) The number of orbits of edges is finite.

4) The graph $K$ is fine, that is, for every $n\in \mathbb N$, any
edge of $K$ is contained in finitely many circuits of length $n$.
(Here circuit means a cycle without self--intersections).
\end{defn}

\subsection*{The definition of Farb}

Let $G$ be a group generated by a finite set $X$ and let $\{ H_1,
H_2, \ldots H_m\} $ be a collection of subgroups of $G$. We begin
with the Cayley graph $\Gamma (G, X)$ of $G$ and form a new graph
as follows: for each left coset $gH_i$ of $H_i$ in $G$, add a
vertex $v(gH_i)$ to $\Gamma (G, X)$, and add an edge $e(gh)$ of
length $1/2$ from each element $gh$ of $gH_i$ to the vertex
$v(gH_i)$. The new graph is called the {\it coned--off Cayley
graph} of $G$ with respect to $\{ H_1, H_2, \ldots H_m\} $, and is
denoted by $\widehat \Gamma (G, X)$. We give this graph the path
metric. Note that $\widehat \Gamma (G, X)$ is not a proper metric
space as closed balls are not necessarily compact.

\begin{defn} \label{FarbDef}
The group $G$ is {\it hyperbolic relative to $\{ H_1, H_2, \ldots
H_m\} $} if the coned--off Cayley graph $\widehat \Gamma (G, X)$
of $G$ with respect to $\{ H_1, H_2, \ldots H_m\} $ is a
hyperbolic metric space.
\end{defn}

\begin{defn}
Given a path $p$ in $\widehat \Gamma (G, X)$, we say that $p$ {\it
penetrates} the coset $gH_i$ if $p$ passes through the cone point
$v(gH_i)$; a vertex $v_1$ (respectively $v_2$) of the path $p$
which precedes to $v(gH_i)$ (respectively succeeds to $v(gH_i)$)
is called an {\it entering vertex} (respectively an {\it exiting
vertex}) of $p$ in the coset $gH_i$. Notice that entering and
exiting vertices are always vertices of $\Gamma (X,G)$. A path $p$
in $\widehat \Gamma (G, X)$ is said to be a path {\it without
backtracking} if, for every coset $gH_i$ which $p$ penetrates, $p$
never returns to $gH_i$ after leaving $gH_i$.
\end{defn}

\begin{defn}\label{BCPdef}
(Bounded coset penetration). The pair $(G, \{ H_1, H_2, \ldots
H_m\} )$ is said to satisfy the {\it Bounded Coset Penetration
property} (or BCP property for brevity) if, for every $\lambda\ge
1$, there is a constant $a=a(\lambda )>0$ such that if $p$ and $q$
are $(\lambda ,0) $--quasi--geodesics without backtracking in
$\widehat\Gamma (G, X)$ such that the endpoints of $p$ and $q$ are
in $\Gamma (G, X)$, $p_-=q_-$, and $\dx (p_+, q_+)\le 1$, then the
following conditions hold.

(1) If $p$ penetrates a coset $gH_i$ but $q$ does not penetrate
$gH_i$, then the entering vertex and the ending vertex of $p$ in
$gH_i$ are an $X$--distance of at most $a$ from each other.

(2) If both $p$ and $q$ penetrate a coset $gH_i$, then the
entering vertices of $p$ and $q$ in $gH_i$ lies an $X$--distance
of at most $a$ from each other; similarly for the exiting
vertices.
\end{defn}

\begin{ex}
The group $G=\langle a,b \; | \; [a,b]=1\rangle \cong \mathbb
Z\times \mathbb Z$ is weakly hyperbolic relative to the cyclic
subgroup $H=\langle a\rangle $. However the pair $(G,H)$ does not
have the BCP property: the paths $a^n$ and $ba^n$ are relative
geodesics without backtracking ending a distance $1$ apart in
$\Gamma $, but they clearly violate condition (1) of Definition
\ref{BCPdef} when $n$ is large enough.
\end{ex}

Dahmani \cite{Dah1} shows that $G$ satisfies Definition
\ref{BowDef} if and only if it satisfies Definition \ref{FarbDef}
and the pair $(G, \{ H_1, H_2, \ldots H_m\} )$ has the BCP
property. (The proof of this fact in \cite{Dah1} contains some gaps; the complete version of the proof is available in \cite{DahPhD}.)

Now we are going to reformulate Farb's definition in terms of the
relative Cayley graph $\G $.

\begin{defn}\label{QI}
Two metric spaces $M_1, M_2$ are said to be {\it quasi--isometric}
if there exist $\lambda >0$, $c\ge 0$, $\varepsilon \ge 0$, and a
map $\alpha : M_1\to M_2$ such that the following two condition
hold.

\begin{enumerate}
\item For any $x,y\in M_1$, we have $$ \frac{1}{\lambda }
dist_{M_1} (x,y)-c \le dist_{M_2} (\alpha (x),\alpha (y))\le
\lambda dist_{M_1} (x,y)+c. $$

\item For any $z\in M_2$ there exists $x\in M_1$ such that $$
dist_{M_2}(\alpha (x), z)\le \varepsilon .$$
\end{enumerate}
\end{defn}

Recall that hyperbolicity of metric spaces is invariant under
quasi--isometry.

\begin{lem} \label{po-chel1}
Let $G$ be a group, $\{ H_1, H_2, \ldots , H_m\} $ a finite
collection of subgroups of $G$. The conned--off Cayley graph
$\widehat \Gamma $ is quasi--isometric to $\G $ endowed with the
relative metric $\dxh $. In particular, $\widehat \Gamma $ is
hyperbolic if and only if so is $\G $.
\end{lem}

\begin{proof}
Note that the identity map on $G$ induces in isometric embedding
$\iota $ of the vertex set $V(\G )$ of $\G $ to $\widehat \Gamma $
and $\widehat \Gamma $ belongs to the closed $1$--neighborhood of
the image $\iota (V(\G ))$.
\end{proof}

Clearly, the BCP property can be rewritten as follows (see Section
2.2 for necessary definitions).

\begin{lem}\label{po-chelovecheski}
Let $G$ be a group generated by a finite set $X$, $\{ H_1, H_2,
\ldots , H_m\} $ a finite collection of finitely generated
subgroups of $G$. The pair $(G, \{ H_1, H_2, \ldots , H_m\} )$
satisfies the BCP property if and only if for any $\lambda \ge 1$,
there exists constant $a=a(\lambda )$ such that the following
conditions hold. Let $p$, $q$ be $(\lambda , 0)$--quasi--geodesics
without backtracking in $\G $ (in the sense of Definition
\ref{backtrac}) such that $p_-=q_-$, $\dx (p_+, q_+)\le 1$.

1) Suppose that for some $i$, $s$ is an $H_i$--component of $p$
such that $\dx (s_-, s_+) \ge a$; then there exists an
$H_i$--component $t$ of $q$ such that $t$ is connected to $s$.

2) Suppose that for some $i$, $s$ and $t$ are connected
$H_i$--components of $p$ and $q$ respectively. Then $\dx (s_-,
t_-)\le a $ and $\dx (s_+, t_+)\le a$.
\end{lem}

\subsection*{Proof of the main theorem}

Our main goal here is to prove

\begin{thm}\label{MainThap}
Let $G$ be a group generated by a finite set $X$,  $\{ H_1, H_2,
\ldots , H_m \}$ a collection of subgroups of $G$. Then the
following conditions are equivalent.

1) $G$ is relatively finitely presented with respect to $\{ H_1,
H_2, \ldots , H_m \}$ and the relative Dehn function of the pair
$(G, \{ H_1, H_2, \ldots , H_m \} )$ is linear.

2) $G$ is hyperbolic with respect to the collection $\{ H_1, H_2,
\ldots , H_m \}$ in the sense of Farb and satisfies the BCP
property (or, equivalently, $G$ is hyperbolic with respect to $\{
H_1, H_2, \ldots , H_m \}$ in the sense of Bowditch).
\end{thm}

\begin{proof}
Theorem \ref{TBCP} gives the implication $ 1)\Rightarrow 2)$. Let
us show that Farb's definition implies relative hyperbolicity in
the sense of our paper.

Let $p$ be a cycle in $\G $. We say that $p$ is {\it atomic } if
any subpath $q$ of $p$ of length $l(q)\le 1/2 l(p)$ is geodesic in
$\G $. Recall that $\G $ is hyperbolic by Lemma \ref{po-chel1}. In
what follows we denote by $\delta $ the hyperbolicity constant of
$\G $.

\begin{lem}\label{atomic}
Let $p$ be an atomic cycle in $\G $. Then the following conditions
hold.

1) $l(p)\le 4\delta +9$.

2) For any $i=1, \ldots , m$, any $H_i$--component of $p$ is
isolated.
\end{lem}

\begin{proof}
1) Suppose that $l(p)>4\delta +9$. We can represent $p$ as the
product $p=p_1p_2p_3$, where $l(p_1)=l(p_2)$ and $l(p_3)\le 1$.
Then $p_1p_2p_3$ is a geodesic triangle in $\G $ since $p$ is
atomic. At least one of the sides $p_1$, $p_2$ is longer than
$2\delta +4$. Assume that $l(p_1)>2\delta +4$. Let us take the
middle point $w$ of $p_1$. Since $\G $ is $\delta $--hyperbolic,
there exists a point $z\in p_2p_3$ such that $\dxh (w,z)\le \delta
$. Let $v$ (respectively $u$) be a vertex of $p_1$ (respectively
$p_2p_3$) that is closest to $w$ (respectively $z$). Then $\dxh
(v,u)\le \delta +1$. Notice that both the segments $[u,v]$ and
$[v,u]$ of $p$ are not geodesic as both of them have length at
least $l(p_1)/2-1/2>\delta +3/2$. However, at least one of these
segments has length at most $1/2l(p)$ contradictory to the
assumption that $p$ is atomic.

2) Let $p=arbs$, where $a,b$ are connected $H_i$--components of
$p$. Then both the subpaths $ar$ and $bs$ are not geodesic since
$$\dxh (a_-,r_+)=\dxh (b_-, s_+)= 1.$$ Thus we get a contradiction
as above.
\end{proof}

\begin{cor}
Let $\mathcal A$ be the set of labels of all atomic cycles in $\G
$. Then $card \; \mathcal A< \infty $.
\end{cor}

\begin{proof}
Let $p$ be an atomic cycle. We can represent $p$ as $p=cd$, where
$l(c)\le l(d)\le l(c)+1$. Since $p$ is atomic, $c$ is geodesic and
$d$ is $(2,0)$--quasi--geodesic as any proper subpath of $d$ is
geodesic. Since all components of $p$ are isolated, $c$ and $d$
are paths without backtracking. Therefore, for any $i$, the
$X$--length of every $H_i$--component of $p$ is at most $a$, where
$a=a(2)$ is the constant provided by Lemma \ref{po-chelovecheski}.
Hence there are only finitely many possibilities for labels of
$H_i$--components of atomic cycles. This fact together with the
first assertion of Lemma \ref{atomic} implies the finiteness of
$\mathcal A$.
\end{proof}

Let us return to the proof of the theorem. We are going to show
that $G$ has the relative presentation $$ \langle X, \mathcal H\;
| \; R=1,\; R\in \mathcal A\rangle $$ and the corresponding
relative Dehn function satisfies $\delta ^{rel}(n)\le 2^n$.  Let
$W$ be a word of length at most $n$ over $(X\cup \mathcal H)^\ast
$ representing $1$ in $G$. We have to show that there exists a van
Kampen diagram with boundary label $W$ having at most $2^n$ cells
labelled by words from $\mathcal A$.

We proceed by induction on $n$. Consider the cycle $p$ in $\G $
labelled $W$. If $p$ is atomic, the existence of the required
diagram is obvious. If $p$ is not atomic, then without loss of
generality we can assume that $p$ is combinatorially homotopic to
a product $q_1q_2$ of cycles $q_1, q_2$ of length $l(q_i)\le
l(p)-1$, $i=1,2$. This means that up to a cyclic shift
$W=_FV_1V_2$, where $\| V_i\|\le\| W\| -1\le n-1$, and $V_i$
represents $1$ in $G$ for $i=1,2$. By the inductive assumption,
there exist van Kampen diagrams with boundary labels $V_1$, $V_2$
and with the number of cells labelled by words from $\mathcal A$
at most $2^{\| V_1\| }$ and $2^{\| V_1\| }$ respectively. Gluing
these diagrams together in the obvious way, we obtain a van Kampen
diagram with boundary label $W$ and with the number of cells
labelled by words from $\mathcal A$ at most
$$2^{\| V_1\| }+2^{\| V_2\| }\le 2^{n-1}+2^{n-1}=2^n.$$

Thus we proved that the relative Dehn function of $G$ with respect
to $\{ H_1, \ldots , H_m\} $ is well--defined. Since $\G $ is
hyperbolic, the relative Dehn function is, in fact, linear by
Corollary \ref{Gammahyp}.
\end{proof}


\begin{thebibliography}{99}

\addcontentsline{toc}{chapter}{Bibliography}

\bibitem{Ale57}
A.D. Alexandrov, {\it Rulled surfaces in metric spaces}, Vestnik
Leningrad. Univ., {\bf 12} (1957), 5--26.


\bibitem{Ali}
E. Alibegovic, {\it A Combination Theorem for Relatively
Hyperbolic Groups Authors,} prep., 2003, available at
http://xxx.lanl.gov/abs/math.GR/0310257.


\bibitem{AloB}
J.M. Alonso, M.R. Bridson, {\it Semyhyperbolic groups}, Proc.
London Math. Soc., (3) {\bf 70} (1995), 56--114.

\bibitem{Bah}
P. Bahls, {\it Relative hyperbolicity and right-angled Coxeter
groups,}, prep., 2004, available at
http://xxx.lanl.gov/abs/math.GR/0401280.

\bibitem{BS62}
G. Baumslag, D. Solitar, {\it Some two-generator one-relator
non-Hopfian groups},  Bull. Amer. Math. Soc. {\bf 68} (1962)
199--201.


\bibitem{BerN93}
V.N. Berestovskii, I.G. Nikolaev, {\it Multidimensional
generalized Riemannian spaces,} Geometry IV (Yu. Reshetnyak, ed.),
Encyclopedia of Math. Sci., Vol. 70, Springer--Verlag, 1993,
165--243.

\bibitem{Bez}
V.N. Bezverkhnii, V.A. Grinblat, {\it The root problem in Artin
groups} (Russian), Algorithmic problems of the theory of groups
and semigroups, pp. 72--81, Tulsk. Gos. Ped. Inst., Tula, 1981

\bibitem{BoG}
O.V. Bogopolskii, V.N. Gerasimov, {\it Finite subgroups of
hyperbolic groups,} Algebra and Logic {\bf 34} (1995), 343--345.

\bibitem{Bow91}
B.H. Bowditch, {\it Notes on Gromov's hyperbolicity criterion for
path--metric spaces,} Group Theory from a Geometrical Viewpoint
(E. Ghys, A. Haefliger, A. Verjovsky, ed.), Proc. ICTP Trieste
1990, Word Scientific, Singapore, 1991, 373--464.

\bibitem{Bow95}
B.H. Bowditch, {\it A short proof that a subquadratic
isoperimetric inequality implies a linear one,} Michigan J. Math.,
{\bf 42} (1995), 103--107.

\bibitem{Bow-98}
B.H. Bowditch, {\it A topological characterization of hyperbolic
groups,} J. Amer. Math. Soc., {\bf 11} (1998), 643--667.

\bibitem{B96}
B.H. Bowditch, {\it Convergence groups and configuration spaces,}
In: Group Theory Down Under (ed. J.Cossey, C.F. Miller, W.D.
Neumann, M. Shapiro), de Gruyter, 1999, 23--54.

\bibitem{B1}
B.H. Bowditch, {\it Relatively hyperbolic groups,} prep., 1999.

\bibitem{Bow03}
B.H. Bowditch, {\it Intersection numbers and the hyperbolicity of
the complex of curves,} prep., 2002.

\bibitem{BrC}
S.G. Brick, J.M. Corson, {\it Annular Dehn functions of groups,}
Bull. Austral. Math. Soc. 58 (1998), no. 3, 453--464.


\bibitem{BC1}
S.G. Brick, J.M. Corson, {\it Dehn functions and complexes of
groups,} Glasgow Math. J., {\bf 40} (1998), 33--46.

\bibitem{BC}
S.G. Brick, J.M. Corson, {\it On Dehn functions of amalgamations
and strongly undistorted subgroups,} IJAC, {\bf 10} (2000), 5,
665--681.

\bibitem{Bri}
M. Bridson, {\it Polynomial Dehn functions and the length of
asynchronously automatic structures,} Proc. LMS, to appear.

\bibitem{BriH}
M. Bridson, A. Haefliger, Metric spaces of non--positive
curvature, Springer, 1999.

\bibitem{Bum}
I. Bumagin, {\it Conjugacy problem for relatively hyperbolic
groups}, submitted to Alg. Geom. Topology.

\bibitem{Bum1}
I. Bumagin, {\it On the definition of relatyvely hyperbolic
groups,} prep., 2004, available at
http://xxx.lanl.gov/abs/math.GR/0402072.

\bibitem{BM}
M. Burger, S. Mozes, {\it finitely presented simple groups and
products of trees}, C. R. Acad. Sci. Paris, {\bf 324} (1997),
747--752.

\bibitem{Camm}
R. Camm, {\it Simple free products,} J. London Math. Soc., {\bf
28} (1953), 66-76.

\bibitem{Can}
F.B. Cannonito, R.W. Gatterdam, {\it The word problem and power
problem in $1$-relator groups are primitive recursive}, Pacific J.
Math. {\bf 61} (1975), no. 2, 351--359.

\bibitem{Col}
D.J. Collins, {\it The word, power and order problems in finitely
presented groups,} Word problems: decision problems and the
Burnside problem in group theory (Conf., Univ. California, Irvine,
Calif., 1969; dedicated to Hanna Neumann), 401--420. Studies in
Logic and the Foundations of Math., Vol. 71, North-Holland,
Amsterdam, 1973

\bibitem{Com}
L.P.J. Comerford, {\it A note on power-conjugacy,} Houston J.
Math. {\bf 3} (1977), no. 3, 337--341.

\bibitem{CoDP}
M. Coornaert, T. Delzant, A. Papadopoulos,  Notes sur les groupes
hyperboliques de Gromov, Springer LNM 1441, 1990.

\bibitem{Co1}
J.M. Corson, {\it Groups acting on complexes and complexes of
groups,} Geometric Group Theory (Charney, Davis and Shapiro,
eds.), Walter de Gruyter, Berlin, New York, 1955, 79--97.

\bibitem{Co2}
J.M. Corson, {\it Howie diagrams and complexes of groups,} Comm.
Alg., {\bf 23 } (1995), 14, 5221--5242.

\bibitem{Dah1}
F. Dahmani, {\it Classifying space and boundary for relatively
hyperbolic groups}, Proc. of London Math. Soc., to appear.

\bibitem{DahPhD}
F. Dahmani, Les groupes relativement hyperboliques et leurs bords,
PhD thesises, 2003.

\bibitem{Dah2}
F. Dahmani, {\it Combination of convergence groups,} prep., 2002.

\bibitem{DD}
W. Dicks, M.J. Dunwoody, Groups acting on graphs, Cambridge Univ.
Press, Cambridge, 1989.

\bibitem{E}
P. Eberlein, {\it Lattices in spaces of nonpositive curvature,}
Annals of Math., {\bf 111} (1980), 435--476.

\bibitem{Ep-etal}
D.B.A. Epstein, J. Cannon, D.F. Holt, S. Levy, M.S. Paterson, W.P.
Thurston, Word processing in groups, Jones and Bartlett, 1992.

\bibitem{F0}
B. Farb, {\it The extrinsic geometry of subgroups and the
generalized word problem}, Proc. London Math. Soc. {\bf 68}
(1994), 3, 577--593.

\bibitem{F}
B. Farb, {\it Relatively hyperbolic groups,} GAFA, {\bf 8} (1998),
810--840.

\bibitem{Fine}
B. Fine, {\it On power conjugacy and ${\rm SQ}$-universality for
Fuchsian and Kleinian groups}, Modular functions in analysis and
number theory, 41--54, Lecture Notes Math. Statist., 5, Univ.
Pittsburgh, Pittsburgh, PA, 1983. .

\bibitem{Fred}
E.M. Freden, {\it Properties of convergence groups and spaces},
Conform. Geom. Dynam., {\bf 1} (1997), 13--23.

\bibitem{GM}
F.W. Gehring, G.J. Martin, {\it Discrete quasiconformal groups I},
Proc. London Math. Soc., {\bf 55} (1987), 331--358.

\bibitem{Ge1}
S. Gersten, {\it Reducible diagrams and equations over groups,}
Essays in Group Theory (S. Gersten, eds.), MSRI Publ.,
Springer--Verlag, 1987, 15--74.

\bibitem{Ger}
S.M. Gersten, {\it Subgroups of word hyperbolic groups in
dimension 2,} J. London Math. Soc., {\bf 54} (1996), 261--283.

\bibitem{GS}
S.M. Gersten, H.Short, {\it Rational subgroups of biautomatic
groups,} Ann. of Math., (2) {\bf 134} (1991), 125--158.

\bibitem{GhH}
E. Ghys, P. de la Harpe, Eds., Sur les groupes hyperboliques
d'apr\'es Mikhael Gromov, Progress in Math., 83, Birka\"user,
1990.

\bibitem{Gr1}
M. Gromov, {\it Hyperbolic groups,} Essays in Group Theory, MSRI
Series, Vol.8, (S.M. Gersten, ed.), Springer, 1987, 75--263.

\bibitem{MG1}
M. Gromov, \textit{Asymptotic invariants of infinite groups},
Geometric group theory, Vol. 2 (Sussex, 1991), 1--295, London
Math. Soc. Lecture Note Ser., 182, Cambridge Univ. Press,
Cambridge, 1993.

\bibitem{GubS}
V.S. Guba, M.V. Sapir, \textit{On Dehn functions of free products
of groups,} Proc. Amer. Math. Soc. \textbf{127} (1999), 7,
1885--1891.

\bibitem{Harv}
W.J. Harvey, {\it Boundary structure of the modular group},
Riemannian Surfaces and Related Topics: Proc. of the 1978 Stony
Brook Conference (I. Kra and B. Maskit, eds.), Ann. of Math. Stud.
97, Princeton, 1981.

\bibitem{Hru}
C. Hruska, Relative hyperbolicity and relative quasiconvexity for countable groups, \emph{Algebr. Geom. Topol.} \textbf{10} (2010), no. 3, 1807-1856.

\bibitem{IO}
S.V. Ivanov, A.Yu. Olshanskii, {\it Hyperbolic groups and their
quotients of bounded exponents}, Trans. Amer. Math. Soc. {\bf 348}
(1996), no. 6, 2091--2138.

\bibitem{Kal}
K.A. Kalorkoti, {\it Decision problems in group theory}, Proc.
London Math. Soc. (3) {\bf 44} (1982), no. 2, 312--332.

\bibitem{KW}
I. Kapovich, D.Wise, {\it The equivalence of some residual
properties of word-hyperbolic groups},  J. Algebra {\bf 223}
(2000), no. 2, 562--583.

\bibitem{Kap}
I. Kapovich, {\it Relative hyperbolicity and Artin groups}, prep.,
2002.

\bibitem{Lar}
L. Larsen, {\it On the computability of conjugate powers in
finitely generated Fuchsian groups}, Acta Math. 139 (1977), no.
3-4, 267--291.

\bibitem{Lip1}
S. Lipschutz, M. Lipschutz, {\it A note on root decision problems
in groups}, Canad. J. Math. {\bf 25} (1973),  702--705.

\bibitem{Lip2}
S. Lipschutz, Ch.F. Miller III, {\it Groups with certain solvable
and unsolvable decision problems}, Comm. Pure Appl. Math. {\bf 24}
(1971), 7--15.

\bibitem{LS}
R.C. Lyndon, P.E. Shupp, Combinatorial Group Theory,
Springer--Verlag, 1977.

\bibitem{Lys}
I.G. Lysenok, {\it On some algoritmic properties of hyperbolic
groups,} Math. USSR Izv., {\bf 35} (1990), 145--163.

\bibitem{MKS}
W. Magnus, A. Karras, D. Solitar, Combinatorial group theory,
Interscience Publ., 1966.

\bibitem{MM}
H. Masur, Y. Minsky, {\it Geometry of complex of curves I:
Hyperbolicity,} Invent. Math. 138 (1999), no. 1, 103--149.

\bibitem{Miller}
C.F. Miller III, {On group Theoretic Decision Problems and their
Classification,} Annals of Math. Stud., {\bf 68}, Princeton Univ.
Press, Princeton NJ, 1971.


\bibitem{Ols-book}
A.Yu. Ol'shanskii, Geometry of defining relations in groups,
Kluwer Academic Publisher, 1991.

\bibitem{Ol91}
A.Yu. Ol'shanskii, {\it Hyperbolicity of groups with subquadratic
isoperimetric inequalities}, Internat. J. Algebra Comput., {\bf 1}
(1991), 281--289.

\bibitem{Ol95}
A.Yu. Olshanskii, {\it {\rm SQ}-universality of hyperbolic groups}
(Russian), Mat. Sb. {\bf 186} (1995), no. 8, 119--132; translation
in Sb. Math. {\bf 186} (1995), no. 8, 1199--1211.

\bibitem{OlsBL}
A.Yu. Olshanskii, {\it On the Bass-Lubotzky question about
quotients of hyperbolic groups}, J. Algebra {\bf 226} (2000), no.
2, 807--817.

\bibitem{OlsIJAC}
A.Yu. Olshanskii, {\it On residualing homomorphisms and
$G$--subgroups of hyperbolic groups}, Int. J. Alg. Comp., {\bf 3}
(1993), 4, 365--409.

\bibitem{OS}
A.Yu. Olshanskii, M.V. Sapir, {\it The conjugacy problem for
groups, and Higman embeddings}, Electron. Res. Announc. Amer.
Math. Soc. {\bf 9} (2003), 40--50.

\bibitem{Wh}
D.V. Osin, {\it Weak hyperbolicity and free constructions},
Contemp. Math., to appear.

\bibitem{Osin}
D.V. Osin, {\it Relatrively hyperbolic groups and embedding
theorems,} preprint, 2004.

\bibitem{OS02}
D.V. Osin, M.V. Sapir, {\it Asymptotic cones of relatively
hyperbolic groups}, preprint, 2002.

\bibitem{Pan1}
A. Pankrat'ev, {\it Hyperbolic products of groups,} Vestnik Mosk.
Universiteta, Ser. 1, Mathematics, Mechanics, {\bf 2} (1999),
9--13 (in Russian).

\bibitem{Pan2}
A. Pankrat'ev, {\it On the infinite hyperbolic quotients of
hyperbolic products of groups}, Fundam. Prikl. Mat. 7 (2001), no.
2, 465--493.

\bibitem{PanPhD}
A. Pankrat'ev, {\it Hyperbolic products of groups,} PhD Thesises,
Moscow State University, 2001.

\bibitem{Paps95}
P. Papasoglou, {\it On the sub--quadratic isoperimetric
inequality,} Geometric Group Theory, (R. Charney, M. Davis, M.
Shapiro, ed.), de Gruyter, Berlin -- New--York, 1995, 149--158.


\bibitem{Rips}
I. Rips, ...

\bibitem{Sel99}
Z. Sela, {\it Endomorphisms of hyperbolic groups I: The Hopf
property,} Topology, {\bf 38} (1999), 301--322.

\bibitem{Sela}
Z. Sela, {\it Diophantine geometry over groups I: Makanin --
Razborov diagrams}, IHES Publ. Math., {\bf 93} (2001), 31--105.

\bibitem{Trees}
J-P. Serre, Trees, Springer--Verlag, 1980, Translation of "Arbres,
Amalgames, $SL_2$", Ast\'erisque, {\bf 46} 1977.

\bibitem{Sh90}
H. Short, {\it Groups and combings}, preprint, ENS Lion, 1990.

\bibitem{Short}
H. Short (editor), {\it Notes on word hyperbolic groups,} Group
Theory from a Geometric Viewpoint, (E. Ghys, A. Haefliger, A.
Verjovsky, ed.), Proc. ICTP Trieste 1990, World Scientific,
Singapore, 1991, 3--64.

\bibitem{Sz}
A. Szczepa\'nski, {\it Relatively hyperbolic groups,} Michigan
Math. J., {\bf 45} (1998), 611--618.

\bibitem{Tuk}
P. Tukia, {\it Convergence groups and Gromov's metric hyperbolic
spaces}, New Zeland J. Math, {\bf 23} (1994), 157--187.

\bibitem{W96}
D. Wise, Non--positively curved squared complexes, aperiodic
tillings, and non--residually finite groups, Ph.D. Thesis,
Princeton Univ., 1996.

\bibitem{W02}
D. Wise, {\it The residual finiteness of negatively curved
polygons of finite groups}, Invent. Math. {\bf 149} (2002), no. 3,
579--617.

\bibitem{Yaman}
A. Yaman, {\it A topplogical characterization of relatively
hyperbolic groups}, prep., 2002.
\end{thebibliography}
\end{document}